%% file: solvability2.tex
\documentclass[12pt]{article}

\input{commands}
\usepackage{epsfig}

\begin{document}

\pagenumbering{roman}

\thispagestyle{empty}

\vbox to 0pt{}

\vskip 30pt

\centerline{SOLVABILITY IN GROUPS OF PIECEWISE-LINEAR}
\centerline{HOMEOMORPHISMS OF THE UNIT INTERVAL}

\vskip 160pt

\centerline{BY}
\vskip 10pt

\centerline{COLLIN BLEAK}
\vskip 10pt

\vskip 220pt

\centerline{COPY OF DISSERTATION}
\centerline{July, 2005}
\vskip 10pt

\newpage

\thispagestyle{empty}

\vbox to 4truein{}

\centerline{copyright by}
\vskip 10pt

\centerline{Collin Bleak}
\vskip 10pt

\centerline{2005}

\newpage

\centerline{\textbf{Abstract}}
We investigate subgroups of the group $\ploi$ of piecewise-linear,
orientation preserving homeomorphisms of the unit interval with
finitely many breaks in slope, and also subgroups of Thompson's group
$F$.  We find geometric criteria determining the derived length of any
such group, and use this criteria to classify the solvable
and non-solvable subgroups of $\ploi$ and of $F$.

Let $H$ be a subgroup of $\ploi$ or $F$.  We find that $H$ is solvable if and
only if $H$ is isomorphic to a group in a well described class $\ws$ of
groups.  We also find that $H$ is non-solvable if and only if we can
embed a copy of a specific non-solvable group $W$ into $H$.

We strengthen the non-solvability classification by finding weak
geometric criteria under which we can embed other groups (all
containing $W$) into non-solvable subgroups of $\ploi$ or $F$.

\newpage

\centerline{\textbf{Acknowledgments}}

The author would like to acknowledge the patient support of the
Mathematics Department at Binghamton University, and the support of
his family and friends.

The author further wishes to thank Matthew Brin and Patricia Bleak for
their support and advice over the years.

\newpage

\tableofcontents

\newpage
\pagenumbering{arabic}
\setcounter{page}{1}

\include{solvIntro2}

\bibliographystyle{amsplain}
\bibliography{solvability2}

\end{document}

%% file: commands.tex
\usepackage{latexsym}
\usepackage{amssymb}
\usepackage{amsmath}
\usepackage{amsfonts}
\usepackage{mathrsfs}
\usepackage[all]{xy}
\usepackage{texdraw}
\usepackage{graphicx}
\usepackage{psfrag}
\usepackage{makeidx}

\newtheorem{Theory}{Theory}[section] 
\newtheorem{theorem}[Theory]{Theorem}
\newtheorem{lemma}[Theory]{Lemma}
\newtheorem{technicalLemma}[Theory]{Technical Lemma}
\newtheorem{corollary}[Theory]{Corollary}
\newtheorem{proposition}[Theory]{Proposition}
\newtheorem{fact}{Fact}  
\newtheorem{remark}[Theory]{Remark} 
\newtheorem{question}{Question} 
\newtheorem{conjecture}[question]{Conjecture}

\newtheorem{Ntn}{Description} 
\newtheorem{Dn}[Ntn]{Definition}

\newcommand{\be}{\begin{enumerate}}
\newcommand{\ee}{\end{enumerate}}
\newcommand{\bq}{\begin{question}}
\newcommand{\eq}{\end{question}}
\newcommand{\bcj}{\begin{conjecture}}
\newcommand{\ecj}{\end{conjecture}}
\newcommand{\bc}{\begin{corollary}}
\newcommand{\ec}{\end{corollary}}
\newcommand{\bl}{\begin{lemma}}
\newcommand{\el}{\end{lemma}}
\newcommand{\btl}{\begin{technicalLemma}}
\newcommand{\etl}{\end{technicalLemma}}
\newcommand{\bt}{\begin{theorem}}
\newcommand{\et}{\end{theorem}}
\newcommand{\bp}{\begin{proposition}}
\newcommand{\ep}{\end{proposition}}
\newcommand{\bft}{\begin{fact}}
\newcommand{\eft}{\end{fact}}
\newcommand{\brk}{\begin{remark}}
\newcommand{\erk}{\end{remark}}
\newcommand{\bd}{\begin{Dn}}
\newcommand{\ed}{\end{Dn}}

\newcommand{\ga}{\alpha}
\newcommand{\gb}{\beta}

\newcommand{\ploi}{PL_o(I) }

\newcommand{\Z}{ \mathbf Z }
 
\newcommand{\N}{\mathbf N }
 
\newcommand{\R}{\mathbf R }

\newcommand{\mm}{\mathscr{M}}
\newcommand{\ws}{\mathscr{R}}

\newcommand{\cg}{\mathscr{G}}
\newcommand{\bsc}{\mathcal{B}}
\newcommand{\wc}{\mathcal{W}}
\renewcommand{\sc}{\mathcal{S}}
\newcommand{\pc}{\mathcal{P}}

\newcommand{\rtr}{ \R\times\R }

\newcommand{\bSeq}[4]{\left\{#1_{#2}\right\}_{#2=#3}^{#4}}
\newcommand{\pSeq}[4]{(#1_{#2})_{#2=#3}^{#4}}

\newcommand{\numSet}[1]{\left\{1,\,\ldots\,,#1\right\}}
\newcommand{\ol}[1]{\overline{#1}}

\newcommand{\pow}[1]{\mathscr{P}(#1)}

\author{Collin Bleak}

%% file: solvIntro2.tex
\vbox to 1.5truein{}
\section{Introduction}
We find descriptions of solvable and non-solvable subgroups of
$\ploi$, the group of orientation preserving piecewise linear
homeomorphisms of the unit interval with finitely many breaks in
slope.  One description we find of the solvable subgroups of $\ploi$
is as the set of isomorphism classes of the subgroups of an easily
described, countable set $\mathscr{M}
=\left\{G_0,G_1,G_2,\ldots\right\}$ of countable, solvable subgroups
of $\ploi$.  We also find one non-solvable group $W$ so that a
subgroup $H$ of $\ploi$ is non-solvable if and only if $H$ contains an
isomorphic copy of $W$.  We will say more about these descriptions,
and others, later on in the introduction.

Much of our analysis of the solvable and non-solvable subgroups of
$\ploi$ does not rely on the piecewise linear nature of the
homeomorphisms of the unit interval; we tend to use point ``dynamics''
under homeomorphisms of the interval.  It would be interesting to
discover how close the classifications below come to classifiying the
solvable and non-solvable subgroups of $Homeo_+(I)$, the group of
orientation preserving homeomorphims of the unit interval.

Unknown to the author at the time of this work, Andr\'es Navas
(\cite{NavasSolv}, \cite{NavasPL}) has also been using dynamics to
analyze the solvable subgroups of the group $Diff^2_+(\R)$.  In
\cite{NavasPL}, he uses his results to show that a finitely generated
solvable subgroup of $\ploi$ with connected support is isomorphic to
a semi-direct product of a group $H$ with the integers $\Z$, where
$H$ is a group in a special class of groups (we also discover this
class of groups).  This result is contained in our investigations
below, but his techniques are different from our own.  Perhaps a
combination of Navas' techniques with those below may be sufficient to
classify the solvable and non-solvable subgroups of $Homeo_+(I)$.

Our investigations also have an impact on the theory of Thompson's group
$F$.  Both the class $\mathscr{M}$ of groups that we realize in
$\ploi$, and the group $W$, can be realized in the standard
realization of $F$ in $\ploi$, so that we also have descriptions of
the solvable and non-solvable subgroups of $F$.

\subsection{An example free statement}
Most of our understanding of the subgroups of $\ploi$ is derived
through analysis of important examples of subgroups of $\ploi$.
Nonetheless, we give here an example free description of the solvable
subgroups of $\ploi$.

Define $\ws$ to be the smallest
non-empty class of groups which is closed under the following three
operations: 

\be

\item Restricted wreath product with $\Z$ (ie. $H\mapsto H\wr\Z$).
\item Bounded direct sum (defined below).
\item Taking subgroups.

\ee 
Here the ``bounded direct sum'' is a countable direct sum of
groups in $\ws$ for which there is a uniform bound on the lengths of
their derived series.

\vspace{.01 in}

With this definition in place, we can give another algebraic version of
our result on the solvable groups in $\ploi$: 

\bt
\label{solveClassification}

$H$ is isomorphic to a solvable subgroup of
$\ploi$ if and only if $H$ is isomorphic to a group in the class $\ws$
of groups.

\et
\subsection{Key examples}
In this section, we will mention some key examples.  All of these
examples can be realized in $\ploi$, and even in Thompson's group $F$.
Our examples rely on an understanding of the wreath product, both as a
standard restricted wreath product as discussed by P. M. Neumann in
\cite{NeumannW}, and as a permutation wreath product which is
discussed in detail in section \ref{pwp}.  In situations such as 
\[
(\ldots((G_1\wr G_2)\wr G_3)\wr\ldots\wr G_n)
\]
where parenthesies ``accumulate on the left'' and each $G_i$ acts on
itself by right multiplication, the resulting groups are the same if
$\wr$ represents the standard restricted wreath product or the
permutation wreath product.  Until we need to distinguish the two
types of wreath product we will not do so.

Before progressing into this discussion, let us fix $\N$ as
representing the positive integers for the remainder of our
investigations.

First, let us build the class $\mathscr{M}$ of groups.  Define
\[G_0=\left\{1\right\},\]
the trivial group, and for each $n\in\N$,
inductively define
\[
G_n=\bigoplus_{i\in\Z}(G_{n-1}\wr \Z).
\]
Note, for example, that $G_1\cong\bigoplus_{i\in\Z}\Z$.  Now define
\[\mm = \left\{G_i\,|\,i\in\Z, i\geq 0\right\}.\]

The class $\mm$ of groups has nice properties, the first of
which is obvious, and the latter two of which we will prove later:

\be
\item Given non-negative $i\in\Z$, $G_i\cong \bigoplus_{j\in\Z}
G_i$\qquad(note: the subscript is not the sum index).
\item Given non-negative $n\in\Z$, $G_n$ has derived length $n$.
\item If $H$ is a subgroup of $G_k$ for some non-negative $k\in\Z$, and $H$ has
derived length $n$, then $H$ is isomorphic to a subgroup of $G_n$.
\ee

Our main interest in these groups is that they play a key role in
understanding the class $\ws$.

We will see that taking a restricted wreath product with $\Z$ is
something that is easy to realize in $\ploi$, which is why this
activity plays a key role in the definition of the class $\mm$.  In
fact, it is so natural that another collection of groups becomes
relevant. Define $W_0 = 1$, the trivial group, and for all $i\in\N$,
define $W_i = W_{i-1}\wr\Z$, so that
\[
W_i =
(\ldots(((\Z\wr\Z)\wr\Z)\wr \Z)\ldots)\wr\Z,
\]
where there are $i$ appearances of $\Z$ on the right.  Our main result
towards understanding $\ws$ can be rephrased in terms of the $W_i$
instead of the $G_i$.  However, the $W_i$ lack the corresponding first
and third properties of the $G_i$ stated above, and these are serious
deficiencies in the class, from a computational point of view.

Nonetheless, the $W_i$ are useful since an isomorphic copy of the
group
\[
W = \bigoplus_{i\in\N}W_i
\]
occurs as a subgroup of any non-solvable subgroup of $\ploi$.

We now describe some often recurring groups, all of which contain an
isomorphic copy of $W$ as a subgroup.  The first group is 
\[
(\Z\wr)^{\infty} = \ldots(((\Z\wr\Z)\wr\Z)\wr\Z)\wr\ldots
\]
and the second group is the permutation wreath product:
\[
(\wr\Z)^{\infty} = \ldots\wr(\Z\wr(\Z\wr(\Z\wr\Z)))\ldots
\]
Sapir had raised the question of whether every non-solvable subgroup
of $F$ contained an isomorphic copy of $(\Z\wr)^{\infty}$.  In
\cite{BrinEG}, Brin shows that both $(\Z\wr)^{\infty}$ and
$(\wr\Z)^{\infty}$ occur as non-solvable subgroups of Thompson's group
$F$, but that neither $(\Z\wr)^{\infty}$ nor $(\wr\Z)^{\infty}$
contain the other as a subgroup.  Brin also shows that a third group,
the permutation wreath product
\[
(\wr\Z\wr)^{\infty} = (\wr\Z)^{\infty}\wr(\Z\wr)^{\infty} 
\] 
is also realized in $F$.  The group $(\wr\Z\wr)^{\infty}$ contains
both an isomorphic copy of $(\Z\wr)^{\infty}$ and an isomorphic copy
of $(\wr\Z)^{\infty}$.  A natural follow-up question to Sapir's
initial question is whether one of $(\Z\wr)^{\infty}$ or
$(\wr\Z)^{\infty}$ is always to be found as a subgroup of any
non-solvable subgroup of $F$.  These investigations answer this
question.  See the next section.

\subsection{Formal statements of algebraic results}
There will be a later section giving geometric results which are used
to obtain the results we give here.  In this section, we will state
all of our main results, and make a few comments.  

Our two chief results are Theorem \ref{solveClassification} above, and
the following:

\bt 
\label{nonSolveClassification}

Suppose $H$ is a subgroup of $\ploi$.  $H$ is non-solvable if and
only if $H$ contains a subgroup isomorphic to $W$.  

\et

Theorem \ref{solveClassification} is further explained by the following result:

\bt
\label{RandM}

$H\in\ws$ if and only if $H$ is isomorphic to a subgroup of a group in $\mm$.
\et

This second description of the class $\ws$ of groups will greatly
assists us in doing calculations relating to the solvable subgroups of
$\ploi$ from a purely algebraic perspective.  The proofs of both
Theorem \ref{solveClassification} and Theorem \ref{RandM} will depend
on the following Lemma.

\bl
\label{solveInM}
If $G$ is a solvable subgroup of $\ploi$ with derived
length $n$, then $G$ is isomorphic to a subgroup of $G_n$.  
\el

It is immediate from construction that the groups in $\mm$ are all
countable, so the last lemma also has the following corollary.

\bc
If $H$ is a solvable subgroup of $\ploi$, then $H$ is countable.
\ec

Since the groups in $\mm$ can all be realized in Thompson's group $F$,
we further have:

\bc 

A group $H\in\ws$ if and only if $H$ is isomorphic with a solvable
subgroup of $F$.

\ec

It is easily seen that $W$ admits no finite index solvable
subgroup, so we have another corollary based on Theorem
\ref{nonSolveClassification}.

\bc
Virtually solvable subgroups of $\ploi$ are solvable.
\ec

Theorem $\ref{nonSolveClassification}$ does not tell the whole story
of the non-solvable subgroups of $\ploi$.  Later we give weak
geometric conditions, each of which implies that a subgroup of $\ploi$
contains a subgroup isomophic to $(\wr\Z\wr)^{\infty}$.

We also show the following:

\bt 

If $H$ is a finitely generated non-solvable subgroup of $\ploi$,
then $H$ contains an isomorphic copy of $(\Z\wr)^{\infty}$ or
$(\wr\Z)^{\infty}$.  

\et 

It is a question whether every finitely generated non-solvable
subgroup of $\ploi$ contains $(\wr\Z\wr)^{\infty}$.  (Note, $W$ is not
finitely generated and non-solvable, and no copy of
$(\wr\Z\wr)^{\infty}$ embeds in it.)

\subsection{Geometry}

In this section, we will realize the $W_i$ in $F\leq\ploi$, and use
these realizations to motivate some geometric definitions which will
enable us to state our main geometric result upon which both Theorem
\ref{solveClassification} and Theorem \ref{nonSolveClassification}
depend.  Here, we are thinking of $F$ as the realization of Thompson's
group in $\ploi$ which consists of all the elements of $\ploi$ which
have all slopes powers of two, and which have all breakpoints occuring
at the dyadic rationals $\Z[\frac{1}{2}]$.  See Cannon, Floyd, and Parry
\cite{CFP} for an introduction to the remarkable group $F$.
\subsubsection{Realizing the $W_i$}
Consider the two elements $\alpha_1$,
$\alpha_2\in\ploi$ defined below:
\[
\begin{array}{llcrr}
x\alpha_1=&\left\{
\begin{array}{lr}
2x&0\leq x<\frac{1}{4},
\\ 
x +\frac{1}{4} & \frac{1}{4}\leq x< \frac{1}{2},
\\ 
\frac{1}{2}x+ \frac{1}{2} &\frac{1}{2}\leq x\leq 1,
\end{array}
\right.&&
x\alpha_2 =&\left\{
\begin{array}{lr}
x&0\leq x<\frac{1}{4},
\\
2x-\frac{1}{4}&\frac{1}{4}\leq x<\frac{5}{16},
\\
x +\frac{1}{16}&\frac{5}{16}\leq x\leq \frac{3}{8},
\\
\frac{1}{2}x +\frac{1}{4}& \frac{3}{8}\leq x< \frac{1}{2},
\\
x &\frac{1}{2}\leq x\leq 1.
\end{array}
\right.
\end{array}
\]
Here are the graphs (superimposed) of these functions:

\begin{center}
\psfrag{a1}[c]{$\alpha_1$}
\psfrag{a2}[c]{$\alpha_2$}
\includegraphics[height=340pt,width=340 pt]{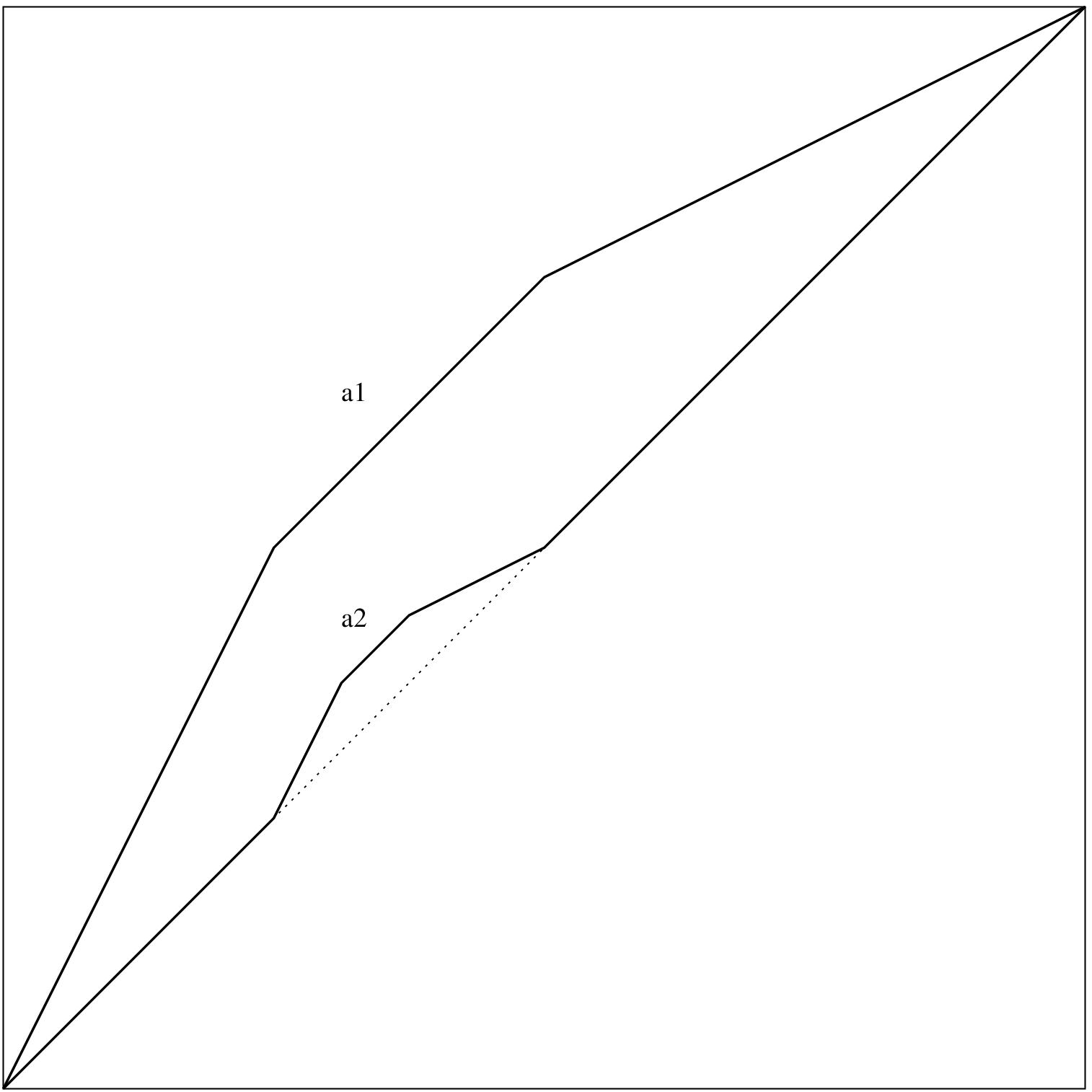}
\end{center}

Either element alone generates a group isomorphic to $\Z\cong W_1$ in
$\ploi$, but the ``action'' of $\alpha_2$ occurs in a single
fundamental domain of $\alpha_1$; that is, $\frac{1}{4}\alpha_1 =
\frac{1}{2}$, but $\alpha_2$ is the identity off of the interval
$[\frac{1}{4},\frac{1}{2}]$.  In particular, $\alpha_2^{\alpha_1} =
\alpha_1^{-1}\alpha_2\alpha_1$ has support
$(\frac{1}{2},\frac{3}{4})$, which is disjoint from the support of
$\alpha_2$. (In this discussion, following the notation in Brin's
papers \cite{BrinU} and \cite{BrinEG}, elements of $\ploi$ act on the
right on $I$, and the support of any particular element of $\ploi$ is
the open set of points in $I$ that are moved by that element.)  In
particular, any two distinct conjugates of $\alpha_2$ by powers of
$\alpha_1$ commute with each other, since their supports will
be disjoint in the interval $I$.  Now, if we consider any element
$h\in\langle \alpha_1,\alpha_2\rangle$, it is a standard algebraic
fact that we can write $h$ as a product of the form a power of
$\alpha_1$ followed by a product of conjugates of $\alpha_2$ (and
$\alpha_2^{-1}$) by various powers of $\alpha_1$.  In particular, the
element $\alpha_1$ generates a group isomorphic to $\Z$ which acts on
the normal subgroup of $\langle \alpha_1,\alpha_2\rangle$ consisting
of the direct sum of $\Z$'s generated by the conjugates of $\alpha_2$
by different powers of $\alpha_1$.  The following chain of
isomorphisms should now make sense:
\[
\langle \alpha_1,\alpha_2\rangle\cong (\bigoplus_{i\in\Z}\Z)\rtimes\Z\cong\Z\wr\Z\cong W_2.
\]
Note before we move on that $\alpha_1$ and $\alpha_2$ are both
elements of Thompson's group $F$.

All of the $W_i$ can be realized in Thompson's group $F$ in an
entirely similar way, using one-bump generators; take $\alpha_i$ to be
the conjugate $\alpha_i = \alpha_{i-1}^s = s^{-1}\alpha_{i-1}s$, where
$s$ is the ``shrinking function'' $s(x) = \frac{1}{4}x+\frac{1}{4}$,
and the conjugation takes place in $PL_o(\R)$, where we replace
$\alpha_1$ by the function in $PL_o(\R)$ which behaves as the identity outside of the unit interval $[0,1]$ for this inductive
definition.  Given $i\in\N$, the support of $\alpha_i$ is contained in
a single fundamental domain of $\alpha_{i-1}$, so that given $j\in\N$,
we have $W_j\cong \langle \alpha_1,\alpha_2,\ldots,\alpha_j\rangle$.
Since conjugating any element of $F$ by $s$ will still produce an
element of $F$, we see that the $W_i$'s can all be realized in $F$.

\subsubsection{Some geometric ideas and results}

Much of the language of this section is motivated by thinking of
subgroups of $\ploi$ as permutation groups acting on the set $I =
[0,1]$.  We will become progressively less formal in this section, as
we are only trying to indicate the nature of our geometric results.
Formal definitions will be given in later sections.

If $H$ is a subgroup of $\ploi$, then its support naturally falls into
a collection of disjoint, open intervals.  Each such we will call an
\emph{orbital}\index{orbital!group} of the group $H$.  Given an
element $\gamma\in H$, the group $\langle\gamma\rangle$ has its own
orbitals, which are the connected components of the support of
$\gamma$.  Given such an interval $A = (a,b)\subset [0,1]$, we call
$A$ an \emph{orbital}\index{orbital!element} of $\gamma$.  Since $H$
now must have infinitely many elements that also have $A$ as an
orbital, we typically use \emph{signed orbitals}, which are pairs of
the form $(A,\gamma)$ to indicate not only which element orbital we
are considering, but also the specific element which is ``owning''
that orbital in our discussion.  In the example groups $W_i$, each of
the $\alpha_k$ had an orbital whose closure was fully contained in the
orbital of $\alpha_{k-1}$, whenever $k>1$.  Such a ``stack'' of $k$
signed, nested orbitals forms what we call a ``tower''.  Observe that many
subgroups of $\ploi$ have pairs of elements with orbitals that
overlap each other on one end, and are therefore not arranged as a tower.  In
this case, we can use the respective elements to move points in $I$
across one orbital into one end of an overlapping orbital, and then
use the second element to move the points across the new orbital out
of the original orbital.  This horizontal technique, and arrangement
of orbitals, motivates our definition of a ``transition chain'', which
consists of a ``chain'' of neighbor-overlapping signed orbitals.
Typically, the elements of a tower generate easier groups to
understand, while the set of elements of a transition chain generate
highly complex groups.

We are now in a position to approach our main geometric result.  We
will say that a group has depth $n\in\N$ if and only if we can find
towers of height $n$, but no towers of height $n+1$, in the group.
Our main geometric result is as follows:

\bt
\label{geoClassification}

Suppose $G$ is a subgroup of $\ploi$ and $n\in\Z$ with $n\geq 0$.  G
is solvable with derived length $n$ if and only if $G$ has depth $n$.

\et

Theorem \ref{geoClassification} is used in both the classification of
the solvable subgroups of $\ploi$ and the classification of the
non-solvable subgroups of $\ploi$.

As hinted at above, the techniques we use can naturally be broken down
into two types: vertical techniques that correspond to working with
towers, and horizontal techniques that correspond to working with
transition chains.  Vertical techniques lead to Theorem
\ref{geoClassification}.  The horizontal techniques extract
information from the transition chains, and then work to get rid of
them.
\newpage

\vbox to 1.5truein{}

\section{Classification of solvable subgroups in $\ploi$\label{solveClassificationSection}}
In this section we will pursue the algebraic classification of the
solvable subgroups in $\ploi$.  One direction of the classification is
mostly algebraic, and requires less knowledge of the terminology of
$\ploi$.  We will engage in that direction first.

\subsection{\label{pwp}Wreath products}
In this paper the restrited wreath product is as given by Neumann in
\cite{NeumannW}.  The restricted wreath product is distinguished from
the unrestricted wreath product in that the restricted wreath product
$A\wr B$ is formed as $C\rtimes B$ where $C$ is a sum of copies of $A$
as opposed to a product of copies of $A$.  One can also think of a
restricted wreath product of two groups as the group resulting from a
permutation wreath product (definition below) where the sets that the two
groups in the product act on are their own underlying sets using right
multiplication.

In \cite{Robinson}, Robinson describes the permutation wreath product,
and proves some of the basic facts about this product.  We will
generally follow his presentation, with only minor modifications.

Let $H$ and $K$ be permutation groups acting on sets $X$ and $Y$
respectively.  We will describe how to construct a new permutation
group $H\wr K$, called the permutation wreath product of $H$ and $K$.
Our new group $H\wr K$ will act on the set $Z = X\times Y$.

If $\gamma\in H$, $\kappa\in K$, and $y\in Y$, define the permutations
$\gamma_y$ and $\kappa^*$ of $Z$ by the rules:
\[
\begin{array}{rl}
\gamma_y: & \left\{
\begin{array}{rl}
(x,y)\mapsto (x\gamma,y),&\\
(x,y')\mapsto(x,y')&\textrm{if}\,\,y'\neq y,
\end{array}\right.\\
\\
\kappa^*:&(x,y)\mapsto(x,y\kappa).
\end{array}
\]
We note that $(\gamma^{-1})_y = (\gamma_y)^{-1}$ and $(\kappa^{-1})^*
= (\kappa^*)^{-1}$, so that $\gamma_y$ and $\kappa^*$ are invertible
and therefore really are elements in $Sym(Z)$ (the group of
permutations of the set $Z$).  One can check that with $y$ fixed, the
map $\phi_y:H\to Sym(Z)$ defined by the rule $\gamma\mapsto\gamma_y$
is a monic homomorphism, and regardless, that $\phi^*: K \to Sym(Z)$
defined by the rule $\kappa\mapsto\kappa^*$ is also a monic
homomorphism. (The monic property follows from the definition of
permutation groups as subgroups of the symmetric group on the set
being acted upon, so that all elements except the identity move
something.)  If we denote the images of $\phi_y$ and $\phi^*$ by $H_y$
and $K^*$ respectively, then \emph{the wreath product of $H$ and $K$
is the subgroup of $Sym(Z)$ generated by all of the $H_y$ for $y\in Y$
and the group $K^*$}\index{product!wreath}, ie. $H\wr K= \langle
H_y,\,K^*\,|\,y\in Y\rangle$.

Because the $\gamma_y$ have disjoint support in $Z$, the $H_y$'s
generate a subgroup $B$ of $H\wr K$ with $B\cong\bigoplus_{y\in Y}H_y$.
The group $B$ is normal in $H\wr K$, since
$(\kappa^*)^{-1}\gamma_y\kappa^*$ acts by mapping $(x,y\kappa)$ to
$(x\gamma,y\kappa)$, and by fixing $(x',y')$ if $y' \neq y\kappa$,
ie., $(\kappa^*)^{-1}\gamma_y\kappa^* = \gamma_{y\kappa}\in B$.  In
particular, conjugating $B$ by an element $\kappa^*$ of $K^*$ simply
permutes the direct factors of $B$ the same way that $\kappa$ permutes
the elements of $Y$.  Now since non-trivial elements of $B$ always
move some element of $Z$, and can only effect the first coordinate of
such elements of $Z$, while non-trivial elements of $K^*$ always move
some element of $Z$, and can only effect the second coordinates of
such elements, we see that $B\cap K^* = 1$ in $H\wr K$, so that $H\wr
K\cong B\rtimes K^*$.  We typically call $B$ the ``Base group'' of
$H\wr K$, and $K^*$ (or $K$ for simplicity) the ``Top group'' of $H\wr
K$.  We will often use the coordinate representation (or implied
variants of it) of $H\wr K$ as the set of elements of the form $(b,k)$
where $b\in B$ and $k\in K$ in our arguments.  We will often denote a
$b\in\bigoplus_{y\in Y}H_y$ by $(b)_{y\in Y}$.

As mentioned in the introduction, we need a little care to specify
when we are using the restricted wreath product or when we are using
the permutation wreath product.  Here is an example: If we temporarily
let $\wr_p$ represent the permutation wreath product, and $\wr$
represent the restricted wreath product, and assume that $\Z$ is a
permutation group acting on its base set by right multiplication in
the permutation wreath products below, then observe the difference in
the following two groups:
\[
\begin{array}{l}
A =\Z\wr(\Z\wr\Z)\cong(\bigoplus_{a\in\Z\wr\Z}\Z)\rtimes\Z\wr\Z,\qquad\textrm{and}\\
B=\Z\wr_p(\Z\wr_p\Z)\cong(\bigoplus_{a\in \Z\times\Z}\Z)\rtimes\Z\wr_p\Z.
\end{array}
\]
The top groups are isomorphic, and so are the base groups, but overall
the groups are very different in the sense of how the summands of the
base group are indexed and moved by the action of the top group.  

Note that differences arise in the two wreath products (restricted and
permutation) since in the permutation wreath product, the set
associated with the group may not actually be the underlying set of
the group.  When we take wreath products on the right with $\Z$, this
will never be an issue, so we will not focus on this issue unless it
needs explicit care (such as in the construction of
$(\wr\Z)^{\infty}$).

Note that in \cite{Robinson}, Robinson shows that the permutation
wreath product of any three permutation groups is associative.  From
the last example, we see that this is not the case for restricted
wreath products.

One should consider the realization of $W_2$ in the introduction to
understand what sets play the roles of $X$, $Y$, and $Z$, and what are
the incarnations of the groups $B$ and $K$ in that context.  Brin's
paper \cite{BrinEG} goes in depth into understanding this realization.

\subsection{The class $\ws$\label{WS}}
Recall that $\ws$ represents the smallest non-empty class of groups
which is closed under the following three operations.

\be
\item Restricted wreath product with $\Z$.
\item Bounded direct sum.
\item Taking subgroups.
\ee 

Here bounded direct sum means a direct sum of groups in $\ws$, all of
whose derived lengths are less than some bounding integer $M$, and where
the index set of the direct sum is countable.

We will investigate this class and come to understand it from a second
perspective, and then by using Theorems
\ref{RandM} and \ref{geoClassification}, we will prove Theorem
\ref{solveClassification}.

We begin to understand $\ws$ via a close study of the class $\mm =
\left\{G_i\,|\, i\geq 0, i\in\Z\right\}$ defined in the introduction.
Let us gather some facts about the groups $G_k$.

\bl 
\label{MFacts}

\be

\item If $F_0$, $F_1$, $H_0$, and $H_1$ are groups, where $F_0\leq
F_1$ and $H_0\leq H_1$, then $F_0\wr H_0\leq F_1\wr H_1$.
\item For any non-negative $m$ and $n\in\Z$, with $m< n$, $G_m$ embeds
  as a normal subgroup of $G_n$.
\item For any group $G$ with derived length $n$, the groups $G\wr\Z$
and $\bigoplus_{i\in\Z}(G\wr\Z)$ have derived length $n+1$.
\item For any non-negative $n\in\Z$, $G_n$ has derived length $n$.

\ee
\el

pf: The first point is immediate by examining the following chain
of subgroup inclusions, where the inclusions are based on the
underlying sets.
\[
\begin{array}{l}
F_0\wr H_0 = \left\{((f_0)_{h\in H_0},h_0)\,|\,(f_0)_{h\in H_0}\in\bigoplus_{a\in H_0}F_0, h_0\in H_0\right\}
\leq \\
F_1\wr H_0 = \left\{((f_1)_{h\in H_0},h_0)\,|\,(f_1)_{h\in H_0}\in\bigoplus_{a\in H_0}F_1, h_0\in H_0\right\}
\leq \\
F_1\wr H_1 = \left\{((f_1)_{h\in H_1},h_1)\,|\,(f_1)_{h\in H_1}\in\bigoplus_{a\in H_1}F_1, h_1\in H_1\right\}
\end{array}
\]
To see the second point, we will demonstrate an embedding of $G_{n-1}$
into $G_n$, and thus inductively define an embedding of $G_m$ into
$G_n$, for any non-negative integers $m <n$.  Note that there are many
copies of $G_{n-1}$ in $G_n$, but we are particularly interested in
the one given in the next paragraph, which is the copy that we use to
inductively define our particular copy of $G_m$ in $G_n$.  The
normality of this embedded copy of $G_m$ in $G_n$ follows easily from
the theory of group actions, which can be checked in section
\ref{realizingGn} below where we realize $G_n$ in $\ploi$.  (Let $X$ be
the set of left hand endpoints of the components of the support of
$G_m$ in $I$.  $X$ is acted upon by $G_n$, and the kernal of this
action is $G_m$.)  A second proof is by noting that the embedded copy
of $G_{n-1}$ in $G_n$ that we demonstrate in the paragraph below is
characteristic in $G_n$, which proof can be carried out with the help
of the geometric tools established in the proof of Theorem
\ref{geoClassification}.  We will not use the normality of our
embedded copy of $G_m$ in $G_n$ later.

Now let us describe our particular embedded copy of $G_{n-1}$ in $G_n$.  First,
identify the base group of $G_{n-1}\wr\Z$ with $G_{n-1}$ using the fact
that $\bigoplus_{i\in\Z}G_{n-1}\cong G_{n-1}$, now since the direct sum
of the base groups of the $G_{n-1}\wr\Z$ summands in the definition of
$G_n$ is also a subgroup of $G_n$, we see that
$\bigoplus_{i\in\Z}G_{n-1}$ is a subgroup of $G_n$.  But now again, this
last direct sum is isomorphic with $G_{n-1}$, so that $G_{n-1}$ (as
embedded here as the direct sum of the base groups of the
$G_{n-1}\wr\Z$ summands of $G_n$) is a subgroup of $G_n$.

For the third point, it follows from Neumann \cite{NeumannW} (Theorem
4.1 and Corollary 4.5) that if $G=A\wr B$, then $G'$ is contained in a
sum of copies of $A$, and also that $G'$ surjects onto $A$ (these
facts are under the condition that $B$ is abelian).  Both facts are
easy exercises in our situation with $B=\Z$.  This immediately implies
the third point.

The fourth point follows directly from the third.

\qquad$\diamond$

Now let us examine $\ws$ for a short time.  We give (modulo work to
come later) a characterization of $\ws$ that completes one direction
of Theorem \ref{solveClassification}.  The key result (Lemma
\ref{threeClassifications} below) will not be used in the rest of the
paper.

Since $\ws$ is nonempty and closed
under subgroups, the group $1$ is a group in $\ws$.  Let $\cg$
represent the class of all groups.  Define $\pc:\pow{\cg}\to\pow\cg$ to be
the function representing the closure operation that takes a set $X$
of groups and computes the smallest class of groups which contains $X$
and is closed under the operations of restricted wreath product with
$\Z$ and bounded direct sum.  Since each $G_i$ is obtained by applying
a finite sequence of bounded direct sums and wreath products with $\Z$
to the trivial group $1$, we see that $\mathscr{M}\subset
\left\{1\right\}\pc$.  By the definition of $\ws$, it is immediate
that $\left\{1\right\}\pc\subset \ws$.  In particular, we have:

\[
\mathscr{M} \subset \left\{1\right\}\pc\subset \ws
\] 

Now let us consider three operators
\[
\begin{array}{l}
\sc:\pow{\cg}\to\pow{\cg},\\
\wc:\pow{\cg}\to\pow{\cg},\\
\bsc:\pow{\cg}\to\pow{\cg},
\end{array}
\]
that represent taking closure under the operations of taking
subgroups, taking restricted wreath products with $\Z$, and building
bounded direct sums, respectively.

We now investigate the relationship between the closure operations of
the three actions; taking subgroups, computing restricted wreath
products with $\Z$, and computing bounded direct sums.  If $\Gamma$ is
the set of all finite length words of the form $(\wc\bsc)^k$ or
$(\bsc\wc)^k$ where $k\in\left\{0,1,2,\ldots\right\}$, and
$X\in\pow{\cg}$, then the union $\Upsilon_{X} =
\cup_{\gamma\in\Gamma}X\gamma$ is the smallest closed set containing
$X$ that is closed under both operations of taking bounded direct sums
and wreath products with $\Z$, so that $\Upsilon_{\left\{1\right\}} =
\left\{1\right\}\pc$.  In particular, the
smallest class of groups which contains the trivial group and which is
closed under the operations of taking wreath products with $\Z$ and
building bounded direct sums equals $\left\{1\right\}\pc$.

\bl 
\label{BWinM}

Let $m$ be a non-negative integer.  If $G$ is a group in
$\left\{1\right\}\pc$ with derived length $m$, then $G$ is isomorphic
to a subgroup of $G_m$.  
\el

pf: We can prove this by inducting on the derived length of $G$.  If
$G$ has derived length $0$, then $G$ is the trivial group so $G =
G_0$.  Let $n\in\N$, and suppose $G$ has derived length $n$ and that
for any group $H$ in $\left\{1\right\}\pc$ which has derived length
$m$ where $0\leq m<n$, we know that $H\leq G_m$.  Now, there is a
$j\in\N$ so that $G\in\left\{1\right\}(\bsc\wc)^j$, or
$G\in\left\{1\right\}(\wc\bsc)^j$.  Note that if
$G\in\left\{1\right\}(\wc\bsc)^j$ then
$G\in\left\{1\right\}(\bsc\wc)^{j+1}$, so we will assume that
$G\in\left\{1\right\}(\bsc\wc)^k$ for some minimal non-negative
integer $k$.  There are now two cases:

\be

\item $G\in\left\{1\right\}(\bsc\wc)^{k-1}\bsc$ but $G$ is not in
$\left\{1\right\}(\bsc\wc)^{k-1}$.

In this case, The last operation required to build $G$ was a bounded
direct sum of groups, all of which groups have derived lengths
necessarily less than or equal to $n$ (note here that a finite
sequence of bounded direct sums is isomorphic to a bounded direct
sum).  The summands of this bounded direct sum are all groups in
$\left\{1\right\}(\bsc\wc)^{k-1}$.  We will argue that each of these
groups is actually a subgroup of $G_n$, and therefore, since a
countable or finite direct sum of groups isomorphic to $G_n$ is
actually isomophic to $G_n$, we will have finished this case.

Let $H$ be a summand of the final bounded direct sum which created
$G$, and assume that the derived length of $H$ is actually $n$ (at
least one summand must have this derived length), and that
$H\in\left\{1\right\}(\bsc\wc)^{k-1}$.  If $H$ is actually in
$\left\{1\right\}(\bsc\wc)^{k-2}\bsc$, then we can replace $H$
inductively by a summand (with derived length $n$) of the last bounded
sum operation used to create $H$, so that there is a
$t\in\N$ so that $H$ is now an element of
$\left\{1\right\}(\bsc\wc)^t$, but not an element of
$\left\{1\right\}(\bsc\wc)^{t-1}\bsc$.  In particular, $H$ is the
result of applying $s$ wreath products with $\Z$ to a
group $H^*$ in $\left\{1\right\}(\bsc\wc)^{t-1}\bsc$ for some $s\in\N$.  $H^*$ is therefore a
group in $\left\{1\right\}\pc$ with derived length $n-s$, and
therefore $H^*$ is a subgroup of $G_{n-s}$.  But note that for any
integer $p$ we have that $G_p\wr \Z\leq G_{p+1}$, so that $H =
H^*(\wr\Z)^s\leq G_{n-s}(\wr\Z)^s\leq G_{n-s+1}(\wr\Z)^{s-1}\leq
\ldots \leq G_n$.

If $H$ is a summand of $G$ with derived length $m$ where $m<n$, Then by
our induction hypothesis, $H\leq G_m$, but $G_m\leq G_n$, so $H\leq
G_n$.

\item $G\in\left\{1\right\}(\bsc\wc)^k$ but $G$ is not in
$\left\{1\right\}(\bsc\wc)^{k-1}\bsc$.

This case is entirely similar to the last case, except that we already
know that $G$ is the result of applying $s$ wreath products with $\Z$
to a group $H$ in $\left\{1\right\}(\bsc\wc)^{k-1}\bsc$ for some
positive integer $s$.  The derived length of $H$ must be $n-s$, and
therefore $H\leq G_{n-s}$ so that $G\leq G_n$ as in the penultimate
paragraph of the previous case.

\ee 
\qquad$\diamond$

We need one last technical lemma before we can complete our
exploration of the class $\ws$:

\bl 
\label{shallowGroupEmbeddings}

Suppose $m$ and  $n\in\Z$, with $0\leq m\leq n$. If $H\leq G_n$ and $H$ has
derived length $m$, then $H$ is isomorphic to a subgroup of $G_m$.

\el

pf: This follows from Corollary \ref{solveEmbeddings} below, and the
fact that we can realize the groups $\mm$ in $\ploi$ (see next
subsection).   \qquad$\diamond$

Our arguments after this section do not rely on the
classification of $\ws$ given in the next lemma.

Finally, we have a nice description of $\ws$.  Note that the following lemma implies Theorem \ref{RandM}.
\bl
\label{threeClassifications}

$\ws = \left\{1\right\}\pc\sc = \mathscr{M}\sc$
\el

Pf: 
We have already shown that $\left\{1\right\}\pc \subset
\mathscr{M}\sc$, so we know that $\left\{1\right\}\pc\sc \subset
\mathscr{M}\sc$, and by definition, $\mathscr{M}\subset
\left\{1\right\}\pc$, so that $\mathscr{M}\sc\subset
\left\{1\right\}\pc\sc$.  In particular, $\mathscr{M}\sc =
\left\{1\right\}\pc\sc$.

We will now show that $\mathscr{M}\sc\bsc = \mathscr{M} \sc$ (implying
that $\mathscr{M}\sc$ is already closed under the operation of taking
bounded direct sums) and that $\mathscr{M} \sc\wc = \mathscr{M} \sc$ (implying
that $\mathscr{M}\sc$ is already closed under wreath products with
$\Z$, so that we can conclude that $\mathscr{M} \sc = \ws$.

To see that $\mathscr{M}\sc\bsc = \mathscr{M}\sc$, let $G\in
\mathscr{M}\sc\bsc$.  There is an $M\in\N$ so that we can write
$G=\bigoplus_{i\in\Z}H_i$, where each $H_i$ has derived length bounded
above by $M$.  Now, each $H_i\leq G_M$ (by the Lemma
\ref{shallowGroupEmbeddings}), so we see that $G \leq
\bigoplus_{i\in\Z}G_M$, hence $G\in \mathscr{M}\sc$.

To see that $\mathscr{M}\sc\wc = \mathscr{M}\sc$, we note that if
$G\in\mathscr{M}\sc\wc$, then either $G\in\mathscr{M}\sc$ or we can
write $G = ((\ldots((H\wr \Z)\wr\Z)\ldots)\wr\Z$, where
$H\in\mathscr{M}\sc$, and where there are $k$ wreath products with
$\Z$, for some $k\in\N$.  In the first case we are done.  In the
second case $H \leq G_M$ for some non-negative integer $M$.  But now
$G=((\ldots(H\wr\Z)\wr\Z)\wr\ldots)\wr\Z\leq
((\ldots(G_M\wr\Z)\wr\Z)\wr\ldots)\wr\Z=J\in \left\{1\right\}\pc$,
where $J$ has derived length $M+k$, so that $G\leq G_{M+k}$.  Finally
we have that $G\in \mm\sc$.  
\qquad$\diamond$

\subsection{Realizing $\ws$ in $\ploi$ and $F$\label{realizingGn}}
Here we will explain how we can realize the groups in $\mathscr{M}$
inside Thompson's group $F$, as realized in $\ploi$.  This will prove
one half of Theorem \ref{solveClassification}.  To realize the other
direction of Theorem \ref{solveClassification}, we will need to take a
lengthy detour through the geometric definitions of $\ploi$.

The elements $\alpha_1$, $\alpha_2\in\ploi$ defined in the
introduction will play a major role here, as will the shrinking
conjugator $s$.

First, observe that we can realize $G_1$ fairly simply.  Let
$\beta_0$ = $\alpha_2$, and define $\beta_k$ for each $k\in\Z$ as
$\beta_0^{\alpha_1^k}$.  $G_1$ is immediately isomorphic to $\langle
\beta_k|k\in\Z\rangle$, and $G_1$ has been realized in Thompson's group $F$.

Now we will show that given any group $H$ which has been realized as a
subgroup of $\ploi$, we can realize $\bigoplus_{i\in\Z}(H\wr \Z)$ as a
subgroup of $\ploi$.  Firstly, conjugate the elements of $H$ by $s$
twice, to create a new group $H_0$ isomorphic with $H$.  The supports
of all of the elements of $H_0$ are contained in the set
$(\frac{5}{16},\frac{3}{8})$, which is contained in a single
fundamental domain of $\alpha_2 = \beta_0$.  At this juncture, the
group generated by $H_0$ and $\beta_0$ is isomophic to $H\wr\Z$, where
$\beta_0$ is the generator of the top $\Z$ factor.  In particular, we
can realize $H\wr\Z$ in $\ploi$.  But now $\bigoplus_{i\in\Z}(H\wr\Z)$
is the base group of $(H\wr\Z)\wr\Z$, which we can also realize by
repeating the previous procedure, so we are done.

  Observe that the shrinking map conjugates elements of Thompson's
group $F$ into Thompson's group $F$, so that if $H$ is realized as a
subgroup of $F$, then this construction of
$\bigoplus_{i\in\Z}(H\wr\Z)$ also produces a subgroup of $F$.  In
particular, we see inductively that each group $G_n$ can now be
realized in Thompson's group $F$, and therefore that each group of $\ws$ can be
realized in Thompson's group $F$.

\subsection{$\ploi$}
Beginning with this section, and through to the end of section
\ref{towers}, we build the geometric tools and
analysis necessary to complete the proof of Theorem
\ref{solveClassification}.

We will now build some required terminology for working with subgroups
of $\ploi$.  We will use notation similar to that in Brin's paper
$\cite{BrinU}$ on the ubiquitous nature of Thompson's group $F$ in
subgroups of $\ploi$.

We note that the set of points moved by an element $h$ of $\ploi$ is
open, by the continuity of elements of $\ploi$.  But then the support
of $H$, for any subgroup $H\leq\ploi$, is a
countable union of pairwise disjoint open intervals in $(0,1)$.  Let the
collection $\mathscr{O}_H$ always denote the countable set of open,
pairwise disjoint intervals of the support of $H$.  We call these
intervals the \emph{orbitals of $H$}\index{orbital!group}. There is a
natural total order on $\mathscr{O}_H$, where if $A$,
$B\in\mathscr{O}_H$, where $A\neq B$, we will say $A<B$ or \emph{$A$
is to the left of $B$}\index{orbitals!natural ordering of} if and only
if given any $x\in{}A$ and $y\in{}B$, we have $x<y$ under the natural
order induced by $I\subset\R$.  Since $A$ and $B$ are disjoint,
connected subsets of $\R$, this definition does not depend on the
choices of $x$ and $y$.

If the collection $\mathscr{O}_H$ is finite, we may speak of the
``first'' orbital, or ``second'' orbital, etc., where the first
orbital is the leftmost orbital under the definition given above, the
second orbital is the orbital to the left of all other orbitals in
$\mathscr{O}_H$ other than itself and the first orbital, and so on.

Given an open interval $A =(a,b)\subset\R$, where $a<b$, we will refer
to $a$ as the \emph{leading end of $A$}\index{end!leading}, and to $b$
as the \emph{trailing end of $A$}\index{end!trailing}.  If the
interval is an orbital of some group $H\in\ploi$, we will refer to
the \emph{ends of the orbital}\index{orbital!end} in the same fashion.

If $h\in{}H$ and $x\in{}Supp(h)$, we will say that \emph{$h$ moves $x$
to the left}\index{movement!left} if $xh<x$, and we will say that
\emph{$h$ moves $x$ to the right}\index{movement!right} if $xh>x$.
Furthermore, we will say that $x\in{}I$ is a \emph{breakpoint for
$h$}\index{breakpoint} if the left and right derivatives of $h$ exist
at $x$, but are not equal.  We recall that by definition, $h$ will
admit only finitely many breakpoints.  If $\mathscr{B}_h$ represents
the set of breakpoints of the element $h$, then
$(0,1)\backslash\mathscr{B}_h$ is a finite collection of open
intervals, which we will call \emph{affine components of
$h$}\index{component!affine}, which admit a natural ``left to right''
ordering as before.  We shall unhesitatingly refer to the ``first''
affine component of $h$, or the ``second'' affine component of $h$,
etc.  We sometimes will refer to the first affine component of $h$ as
the \emph{leading affine component of $h$}\index{component!leading
affine}, and to the last affine component of the domain of $h$ as the
\emph{trailing affine component of $h$}\index{component!trailing
affine}.

Given an element $h\in{}H$ the group $\left<h\right>$ generated by $h$
has its own orbitals, which we will typically call the \emph{orbitals
of $h$}\index{orbital!of an element}.  Any such collection of orbitals for
an element $h\in\ploi$ is finite by Remark \ref{finiteOrbitals}, and
we will denote the number of orbitals for such an element $h$ as
$o_h$.  We will denote the ordered (left to right, as before)
collection of orbitals of $h$ as
$\mathscr{A}^h=\bSeq{A^h}{i}{1}{o_h}$.

The following are some useful remarks.

\brk
\label{finiteOrbitals}
\label{transitiveElementOrbital}
\be

\item If $A$ is an orbital for $h\in H$, then either $xh>x$ for all
points $x$ in $A$, or $xh<x$ for all points $x$ in $A$.

\item Any element $h\in\ploi$ has only finitely many
orbitals.

\item If $h\in\ploi$ and $A=(a,b)$ is an orbital of $h$, then given any
$\epsilon>0$ and $x$ in $A$, there is an integer $n$ so that
$|xh^{-n}-a|<\epsilon$ and $|xh^n-b|<\epsilon$.
\ee
\erk

pf: 

For the first point, the difference of $x$ and $xh$ can never be zero
in an orbital, but the difference function is continuous, so the
intermediate value theorem of calculus implies the difference is
always positive, or always negative, throughout the
orbital.

For the second point, assume it is incorrect, then, $h$ will behave as
the identity function for an infinite sequence of points
$\pSeq{x}{i}{1}{\infty}$ in $I$ so that between any two points $x_i$
and $x_{i+1}$ there would be a point $y_i$ so that $y_ih\neq
y_i$. Hence, $h$ would have to admit to infinitely many affine
components, which contradicts the definition of
$\ploi$.

For the final point, suppose that $h$ moves points to the right on its orbital $A$
(if not, a symmetric argument will prove the result), and let
$\epsilon>0$.

The sequence $y_i = xh^i$ defined for all natural numbers $i$ is an
increasing sequence of points which is bounded above by $b$.  In
particular, it must have a limit $y$ (greater than all of the $x_i$)
by the completeness of the real numbers.  But now $yh = y$ by the
continuity of $h$, so $y$ must be $b$.  Hence, there is a natural
number $m_1$ so that $|y_{m_1}-b|=|xh^{m_1}-b|<\epsilon$.  A symmetric argument
shows that there is a natural number $m_2$ so that $| xh^{-m_2}-
a|<\epsilon$.  Let $n = max(m_1,m_2)$.
\qquad$\diamond$

Given an orbital $A$ of $H$ we say that \emph{$h$ realizes an end of
$A$} if some orbital of $h$ lies entirely in $A$ and shares an
end with $A$.  Note that Brin uses the word ``Approaches'' for this
concept in \cite{BrinU}, but we will use ``Approaches'' to also
indicate the direction in which $h$ moves points, as follows: we will
say that \emph{$h$ approaches the end $a$ of $A$ in $A$} if $h$ has an
orbital $B$ where $B\subset A$ and $B$ has end $a$, and $h$
moves points in $B$ towards $a$.  In particular, $h$ realizes $a$ in $A$
and $h$ moves points in its relevant orbital towards $a$.  If $A$ is
an orbital for $H$ then we say that \emph{$h\in H$ realizes $A$} if
$A$ is also an orbital for $h$.

If $g$ and $h$ are elements of $\ploi$ and there is an interval $A =
(a,b)\subset I$ so that both $g$ and $h$ have $A$ as an orbital, then
we will say that $g$ and $h$ \emph{share the orbital $A$}. 

Given $g\in\ploi$, and any $x\in(0,1)$, there is an $\epsilon>0$ so
that $(x,x+\epsilon)$ is contained in an affine component $C_x$ of
$h$.  The slope of $h$ on $C_x$ is the value of the right hand
derivative $h_+'(x)$ of $h$ at $x$.  Given another element
$h\in\ploi$, we can find the maximal subinterval $[0,a]\subset[0,1]$
so that for any $x\in[0,a]$ we have $h(x) = g(x)$.  If $a=1$, then $h
= g$.  If $a<1$, then $h_+'(a)\neq g_+'(a)$.  In particular, we can
define a total order on the set of elements of $\ploi$, by defining
$g<h$ if $g_+'(a) < h_+'(a)$, where $a$ is the largest point in
$[0,1]$ where $g$ and $h$ are identical on $[0,a]$.  We will call this
the \emph{left total order on $\ploi$}.

\subsubsection{Conjugation and orbitals}
In this section, we will establish some notation to help us to more
flexibly pass to subgroups of a group in $\ploi$.  This section will help us to
understand the orbital structure of conjugates of an element in
$\ploi$.

Let $g$, $h\in\ploi$ and let $k = g^h = h^{-1}gh$.  Suppose that
$\mathscr{O}_g=\bSeq{A}{i}{1}{n}$ are the $n = o_g$ orbitals of $g$ in left
to right order, where $n\in\N$ and $i\in\left\{1,2,\ldots,n\right\}$.
Define
\[
A_i^*=\left\{x\in{}I|xh^{-1}\in{}A_i\right\}=A_ih
\]
for all $i\in\left\{1,2,\ldots,n\right\}$.

The following is standard.

\bl 

$o_k = o_g = n$ and collection $\bSeq{A^*}{i}{1}{n}$ is the ordered
set of orbitals of $g^h$ in left to right order.  

\el

pf: Throughout this proof, let $i$ always refer to an index with
$1\leq i\leq n$.

Suppose $x\in{}I$ so that $xh^{-1}$ is not an element in any $A_i$.
Then $xg^h=xh^{-1}gh=xh^{-1}h=x$, since $xh^{-1}$ is in the fixed set
of $g$.  In particular, the support of $g^h$ is in the union of the
sets $A_i^*$.  Now suppose $x\in{}A_i^*$ for some index $i$.  Then
$xh^{-1}\in{}A_i$, now $xg^h = (xh^{-1}g)\cdot{}h\neq{}x$, since $g$
moves the point $xh^{-1}$.  In particular, each $A_i^*$ is contained
in the support of $g^h$.  Now we see that the support of $g^h$ is
precisely the union of the finite collection of sets $A_i^*$.  Now,
the collection of $A_i^*$'s is a finite collection of pairwise
disjoint, open intervals, being the homeomorphic image of a finite
collection of pairwise disjoint, open intervals. In particular, each
$A_i^*$ is a connected component of the support of $g^h$, and is
therefore an orbital of $g^h$.  Since the set of orbitals
$\bSeq{A}{i}{1}{n}$ is ordered in left to right fashion for $g$, and
since $h$ is order preserving, the indexed collection
$\bSeq{A^*}{i}{1}{n}$ of the orbitals of $g^h$ is also ordered in left
to right fashion.  \qquad$\diamond$

In the setting of the above lemma, we will say that the $A_ih$  are the
\emph{induced orbitals}\index{orbital!induced} of $k$ from $g$ by the
action of $h$.  We might also say that the orbitals of $k$ are induced
from the orbitals of $g$ by the action of $h$.

The following two unrelated points are worth pointing out:
\brk 
\label{breakpoints}
\be
\item Suppose $g$, $h\in\ploi$ and $f = gh$.  If $b$ is a breakpoint of $f$
then $b$ is a breakpoint of $g$ or $bg$ is a breakpoint of $h$.
\item Let $\ga$, $\gb\in\ploi$, then $\ga^{\gb}$ has the same leading
and trailing slopes on its orbitals as\, $\ga$ has on each of its
corresponding orbitals.
\ee
\erk

pf: 

For the first point, suppose $b$ is a breakpoint of $f$, but $b$ is not a
breakpoint of $g$.  The lefthand slope of $f$
at $b$ is the product of the lefthand slopes of $g$ at $b$
and $h$ at $bg$. Also, the righthand slope of $f$
at $b$ is the product of the righthand slopes of $g$ at $b$
and $h$ at $bg$.  Since the lefthand and righthand slopes
of $g$ at $b$ are the same, the lefthand and righthand slopes of $g$
at $bh$ must be different.  In particular, $bg$ is a breakpoint for $h$.

For the second point, let $A = (a,b)$ be an orbital of $\ga^{\gb}$,
corresponding to an orbital $A' = (a',b')$ of $\ga$.  There is
$c\in(a,b)$ so that $(a,c)$ is contained in an affine component of
$\gb^{-1}$, so that $(a',c') = (a,c)\gb^{-1}$ is in an affine
component of $\ga$, and so that $(a',c'') = (a',c')\ga$ is contained
in an affine component of $\gb$.  By Remark \ref{breakpoints}, $(a,c)$
is therefore contained in an affine component of $\ga^{\gb}$, and
hence the slope of $\ga^{\gb}$ on $(a,c)$ is the leading slope of
$\ga^{\gb}$ on $A$, since the leading affine component of $\ga^{\gb}$
with non-trivial intersection with $A$ must now contain $(a,c)$.  The
slope of $\ga^{\gb}$ on $(a,c)$ is the product of the slope of
$\gb^{-1}$ on $(a,c)$, the slope of $\ga$ on $(a',c')$, and the slope
of $\gb$ on $(a',c'')$.  But the functions $\gb$ and $\gb^{-1}$ are
inverse, so if $x$ is in an affine component of $\gb^{-1}$ with slope
$s$, then $y=x\gb$ is in an affine component of $\gb$ with slope
$1/s$.  Now then the image of $(a,c)$ under $\gb^{-1}$ must be in an
affine component $(a',d)$ of $\gb$, and the slopes of $\gb^{-1}$ on
$(a,c)$ and $\gb$ on $(a',d)$ are multiplicative inverses.  But now we
see that $(a',c'')\subset(a',d)$, so in particular, the lead slope of
$\ga^{\gb}$ on $(a,b)$ is the lead slope of $\ga$ on $(a',b')$.

The argument for the trailing slopes is similar. 
\qquad$\diamond$

\subsubsection{Signed orbitals}
We give this discussion its own section, since the definitions we
discuss here are new, and are of extreme importance to understanding
the arguments to follow.  Otherwise, this could easily be included at
the end of the previous section.

Let $1\neq g\in G$ and $A$ be an orbital of $g$.  It is immediate that $A$
is realized by infinitely many elements of $G$ (for instance, non-zero
powers of $g$ will always realize $A$).  Often, we want to associate
some specific element of $G$ which realizes the orbital $A$ with $A$.

Given a group $G\leq\ploi$, we define a \emph{signed
orbital}\index{orbital!signed} (in $G$) to be a pair $(A,g)$ where $A$ is an
orbital of the element $g\in G$.

Now we will build up some convenient language around signed orbitals,
and note some trivial facts.  Given $\gamma = (A,g)$, a signed orbital
of $G$, we will call $A$ the \emph{orbital of $\gamma$}, and $g$ the
\emph{signature of $\gamma$}\index{signature!signed orbital}.  We
observe, as before, that if $n$ is a non-zero integer then $(A,g^n)$
is another signed orbital with orbital $A$.  In particular, if $(A,g)$
is a signed orbital, then so is $(A,g^{-1})$.  Note in this case that
if $g$ moves points to the right on $A$, then $g^{-1}$ moves points to
the left on $A$, and vice-versa, so that we will often replace a
signed orbital by the signed orbital consisting of the same orbital,
but with the inverse of the original signature, so that the new
signature moves points in a desired direction on the orbital.

Given a set $P$ of signed orbitals,  we will often refer to two
induced sets from $P$; the set of orbitals of $P$, denoted $O_P$, and
the set of signatures of $P$, denoted $S_P$.  For completeness,
\[
\begin{array}{l}
O_P=\{A\subset[0,1]\,|\,\textrm{there is }g\in G \textrm{ with } (A,g)\in P\},
\\
S_P=\{g\in G\,|\,\textrm{there is }A \textrm{ an orbital of }g\textrm{ with }(A,g)\in P\}.
\end{array}
\]
Noting the standard (left) total order on elements of $\ploi$, we
observe that the set of signed orbitals of $\ploi$ is partially
ordered by using the lexicographical ordering, where we use the partial order
on subsets of $I$ (induced via inclusion) for the first coordinate,
and the left total order on elements of $\ploi$ for the second
coordinate.  In particular, the set of signed orbitals of any subgroup
of $\ploi$ is also a poset.

Suppose $H$ is a subgroup in $\ploi$, and $H$ has two elements
$\alpha$ and $\beta$, which then generate some subgroup $G$ of $H$.
If we restrict our attention to the orbitals of $G$, one orbital at a
time, we discover that the orbitals of $\alpha$ and $\beta$ must be
arranged so that the orbitals of $\alpha$ cover the fixed set of
$\beta$ within each orbital of $G$, and vice-versa.  This arrangement
of orbitals allows us to move points from one side of an orbital
of $G$ to the other, by repeated application of the elements $\alpha$
and $\beta$, and of their inverses, at the appropriate moments.  A
detailed argument that we can move any point inside an orbital $A$ of
$H$ arbitrarily close to either end of $A$ will be given shortly.

Let $\mathscr{C}=\bSeq{O}{i}{1}{n}$ be a set of signed orbitals of
$G$, so that $O_i = (A_i,g_i)$ for each $i\in \numSet{n}$.  Suppose
that $\mathscr{C}$ satisfies the following properites.

\be
\item For any $k\in\numSet{n-1}$, $A_k\bigcap A_{k+1}\neq\emptyset$.
\item For any $k\in\numSet{n}$, there is a point $p_k\in A_k$,
  which is not in any $A_j$ for $j\neq k$.
\item For any $k\in\numSet{n-1}$, $p_k<p_{k+1}$. 

\ee 

then we will call $\mathscr{C}$ a \emph{transition
chain for the group $G$}\index{transition chain}.  The integer
$n$ will be called the \emph{length of the transition chain
$\mathscr{C}$}\index{transition chain!length}.

Given a transition chain $\mathscr{C}$ of length $n$,
$\mathscr{C}=\{(A_i,g_i)\,|\,i\in\numSet{n}\}$, we fix the notation
$A_{\mathscr{C}} = \bigcup_{i = 1}^n A_i$, as it will often be convenient to refer to the set $A_{\mathscr{C}}$.  

If $\mathscr{C}$ is a transition chain with length $n$ as above, and
$1\leq j\leq k\leq n$, then $\mathscr{C}' = \left\{O_i\right\}_{i =
j}^k$ will be called a \emph{contiguous subset of
$\mathscr{C}$}\index{transition chain!contiguous subset}.  Then it is
an immediate observation that any contiguous subset of a transition
chain is also a transition chain.  Given two transition chains
$\mathscr{C'}$ and $\mathscr{C}$, we will say
$\mathscr{C'}\leq\mathscr{C}$ whenever $\mathscr{C'}$ is a contiguous
subset of $\mathscr{C}$.  Given a group $G$, its set of transition
chains forms a partially ordered set under this binary relation.

Here is a useful lemma:
\bl
\label{transitiveOrbital}
If $H\leq\ploi$ and $A=(a,b)$ is an orbital for $H$, then given any points
$c$, $d\in A$, with $c<d$, there is an element $g\in{}H$ so that
\mbox{$cg>d$}.
\el

pf: $[c,d]$ is contained in an orbital of $H$, and therefore it is
contained in the union of the orbitals of the elements of $H$.  Since
$[c,d]$ is compact, it follows that it is covered by a finite
subcollection $\mathscr{C}'$ of the orbitals of the elements of $H$.
This implies that there is a smallest positive integer $n$ with a
transition chain $\mathscr{C} = \left\{(A_i,g_i)|1\leq i\leq n,
i\in\N\right\}$ for $H$ whose orbitals cover $[c,d]$.  If $n =1$ we
are done by Remark \ref{transitiveElementOrbital}, so assume $n>1$.
We note in passing that $c$ is an element of $A_1$ but not of $A_2$
and $d$ is an element of $A_n$ but $d$ is not an element of $A_{n-1}$.

Improve $\mathscr{C}$ by supposing we chose signatures intelligently,
so that each $(A_i,g_i)\in\mathscr{C}$ satisfies the property that
$g_i$ moves points to the right on $A_i$.

Now, by the definition of transition chain, we know that $B_i =
A_i\cap A_{i+1}$ is actually a non-empty open interval for each
integer $i$ where $1\leq i <n$.  Let $\delta_1$ be the minimum length
of these intervals $B_i$.  By the definition of transition chain, each
orbital $A_i$ has length at least $\delta_1$ (and $2\delta_1$ for
interior orbitals in the chain), since each such orbital has a
distiguished point which is outside of the intersection of this
orbital with the other orbitals of the transition chain.  Let
$\delta_2$ be the minimum of the two distances, one from the left end
of $A_1$ to $c$, and the other from $d$ to the right end of $A_n$.
Now let $\delta = min(\delta_1,\delta_2)/2$.  It is immediate by
construction that each orbital of $\mathscr{C}$ has length greater
than $2\delta$, and also that each interval $B_i$ created above
has length greater than $2\delta$.  For each
$k\in\left\{1,2,\ldots,n\right\}$, let $x_k$ be a point in $A_k$ a
distance less than $\delta$ from the left side of $A_k$.  By Remark
\ref{transitiveElementOrbital}, for each integer $k$ where $1\leq
k\leq n$, there is a positive integer $m_k$ so that $g_k^{m_k}$ will
take $x_k$ to a point $y_k$ within $\delta$ of the right end of $A_k$.
Observing that for all integers $k$ where $1\leq k<n$ we have that
$y_k>x_{k+1}$, we see that $g = g_1^{m_1}g_2^{m_2}\cdots g_n^{m_n}$
moves $x_1 $ to the right of $y_n$.  But now $x_1<c<d<y_n$ by
construction, so $cg>d$.\qquad$\diamond$

\subsection{Thompson's group $F$ and balanced subgroups of $\ploi$}
Brin showed in \cite{BrinU} the following theorem:

\bt [Ubiquitous F]
\label{UbiquitousF}
If a group $H \leq \ploi{}$ has an orbital $A$ so that some element
$h\in{}H$ realizes one end of $A$, but not the other, then $H$ will
contain a subgroup isomorphic to Thompson's group $F$.
\et

The condition on the orbital is weak enough that one readily observes
that ``$\ploi{}$ is riddled with copies of $F$'', quoting Brin.  Hence
it becomes a natural question to ask what can be said about subgroups
of $\ploi{}$ which have no orbitals satisfying the ubiquity condition.

  We will say that an orbital $A$ of a group $H \leq \ploi{}$ is
\emph{imbalanced}\index{orbital!imbalanced} if some element $h \in H$
realizes one end of $A$, but not the other, and we will say $A$ is
\emph{balanced}\index{orbital!balanced} if whenever an element $h \in
H$ realizes one end of $A$, then $h$ also realizes the other end of
$A$ (note that $h$ might do this with two distinct orbitals).
Extrapolating, given a group $H \leq \ploi{}$, we will say that
\emph{$H$ is balanced}\index{group!balanced} if given any subgroup
$G\leq H$, and any orbital $A$ of $G$, every element of $G$ which
realizes one end of $A$ also realizes the other end of $A$.
Informally, $H$ has no subgroup $G$ which has an orbital that is
``heavy'' on one side.  In the case where $H$ has a subgroup $G$ with
an imbalanced orbital, then we will say that \emph{$H$ is
imbalanced}\index{group!imbalanced}.

\brk 

If $H\leq\ploi$ and $H$ is imbalanced, then $H$ has a subgroup
isomorphic to Thompson's group $F$.

\erk

Since $F'$ is non-trivial and simple (\cite{CFP}, Theorem 4.3), $F$ is
not solvable.  Thus imbalanced groups are not solvable.

Possibly less obvious is that the dynamics of balanced groups are much
easier to understand than those of imbalanced groups.  We will trade
heavily on this in the remainder, so the next few subsections will
establish some of the common tools that we will have available to us
when we are analyzing balanced groups.  In the meantime, we build some
tools to help us ``find'' imbalanced groups.

\subsubsection{A useful homomorphism}
\label{phi}
Let us suppose that $H \leq \ploi$ and $A=(a,b)$ is an orbital of $H$.
To simplify the arguments for now, let us suppose that $A$ is the only
orbital of $H$.  We can define a map $\phi:H \to\R\times\R$ defined by
$h\mapsto(h_a,h_b)$ where $h_a=\ln(h_+'(x))$ and $h_b=\ln(h_-'(x))$.
Ie., we take the logs of the slopes of $h$ at the ends of $A$.  Since
$h$ is a p.l.  orientation-preserving homeomorphism of $I$, we see
that the derivatives exist and are positive, and so $\phi$ is well
defined for all $h\in H$.  If $h$ does not realize $a$ (resp. $b$)
then we see that $h$ behaves as the identity near $a$ ($b$) in $A$,
and so $h_a = \ln(1)=0$ ($h_b = 0$).  If $h,g \in H$ then
$hg\phi=h\phi+g\phi$ in $\R\times\R$ by the chain rule. In particular,
we see the following remark:

\brk
$\phi$ is a homomorphism of groups.
\erk

Now the image of $\phi$ is quite interesting, it carries a small
amount of the complexity of $H$, but still enough to allow us to find
out if $A$ is an imbalanced orbital.

\bl
The orbital $A$ is imbalanced if and only if $Im(\phi)$ contains an
element of the form $(\alpha,0)$ or $(0,\alpha)$ where $\alpha \neq
0$.
\el 

pf: Suppose that $h \in H$ approaches one end of $A$ but not the
  other.  Then the slope of $h$ near one end is $1$, but on the other
  end is not $1$.  In particular, $h\phi = (h_a,h_b)$ is either
  $(0,\alpha)$ or $(\alpha,0)$ for some non-zero $\alpha$.  Hence, if
  $A$ is imbalanced, then $Im(\phi)$ contains an element in this form.  Conversely, if $Im(\phi)$ has an element
  $r_\alpha$ of the form $(0,\alpha)$ or $(\alpha,0)$, where $\alpha$
  is not zero, then there is an $h\in H$ so that $h\phi=r_\alpha$.
  For this $h$, either $h_a = 0$ and $h_b = \alpha$, in which case $h$
  realizes $b$ but does not realize $a$, or, in the other case,
  $h_a = \alpha$ and $h_b = 0$, where we see that $h$ realizes $a$
  but not $b$.  In either case, we see that $A$ is an imbalanced
  orbital.\qquad$\diamond$

This next technical lemma will help with the lemma that follows it:
\btl 
\label{longLeadSlope}

Suppose $H\leq\ploi$ and $H$ has an orbital $A=(a,b)$, and that $H$
has a sequence of elements $(g_n)_{n = 1}^{\infty}$ in $H$ which
satisfies the properties below.

\be

\item For each $i\in\N$, the lead slope of $g_{i+1}$ in $A$ is less than the
lead slope of $g_i$ in $A$.

\item Given any real number $q>1$, there is an $i\in\N$ so that the
lead slope of $g_i$ in $A$ is $p$ where $1<p<q$.

\ee

Then there is $c\in(a,b)$ so that given any real number $s>1$, $H$ has an
element $\alpha$ which has an affine component $\Gamma$ containing
$(a,c)$, and with slope $r$ on $\Gamma$ where $1<r<s$.

\etl

pf:

To simplify the language of this argument, we will restrict our
attention to the orbital $A$, treating it as the domain of the
elements of $H$, so that the phrase ``The first affine component of
[$h\in H$]'' will really mean the open interval $(a,u)$ where $h$ has
an affine component of the form $(v,u)$ where $v\leq a$.  We will also
refer to this as the ``leading (or lead) affine component of $h$''.
Similarly, we will refer to the slope of $h$ on its leading affine
component as the ``lead slope of $h$.''

Note that the second condition on $(g_n)_{n=1}^{\infty}$ implies that
every element of $(g_n)_{n=1}^{\infty}$ has lead slope greater than
one.

For each $i\in\N$, let $(a,b_i)$ be the first affine component of
$g_i$, and let $s_i$ represent the lead slope of $g_i$.

We are now in a position to define a new sequence of
functions $\pSeq{h}{i}{1}{\infty}$ which satisfy the following conditions:

\be

\item For each $i\in\N$, the lead slope of $h_i$ is $s_i$.

\item For each $i\in\N$, the leading affine component of $h_i$
contains $(a,b_1)$.

\ee

At this point, by taking $c = b_1$, we see that given any $s>1$, by
the hypothesies on the $g_i$ there will be an $N\in\N$ so that for all
$n>N$ we will have that $s>s_n>1$, and that $h_n$ has leading affine
component containing $(a,c)$.

For each $i\in\N$, define $n_i$ to be the smallest non-negative
integer so that $b_1g_1^{-n_i}<b_i$, and define $h_i =
g_1^{-n_i}g_ig_1^{n_i}$.  Since $g_1$ moves points to the right on its first
affine component, $n_i$ is well defined, and therefore $h_i$ is well
defined.  We now check that $h_i$ satisfies the two conditions, for
each $i$.  Firstly, we observe that the lead slope of each $h_i$ is
the product of the lead slopes of the elements of the product
$g_1^{-n_i}g_ig_1^{n_i}$, which is $(\frac{1}{s_1})^{n_i}s_is_1^{n_i}
= s_i$.  Secondly, by Remark \ref{breakpoints}, we know that if
$x\in(a,b)$ is a breakpoint of the product $g_1^{-n_1}g_ig_1^{n_i}$,
then the image of $x$ under application of some initial partial
product (possibly empty) must be a breakpoint of the next term of the
overall product.  However, the first breakpoint of $g_1^{-1}$ is $d_1
= b_1g_1>b_1$, and $g_1^{-1}$ moves points left in its first affine
component so if $x\in (a,b_1)$, then $xg_1^{-k}<d_1$ for all
non-negative integers $k$.  In particular, if $x\in(a,b_1)$, then $x$
is in the first affine component of $g_1^{-n_i}$ for any $i\in\N$.
Given an $i\in \N$, the first breakpoint of $g_i$ in $A$ is $b_i$, but
$n_i$ was chosen so that $b_1g_1^{-n_i}<b_i$, so if $x\in(a,b_1)$,
then the image of $x$ under $g_1^{-n_i}$ is in the first affine
component of $g_i$, in particular, $(a,b_1)$ is contained in the first
affine component of $g^{-n_1}g_i$.  Finally, the first breakpoint of
$g_1^{n_1}$ is the image of $b_1$ under $g_1^{-n_i + 1}$, but
$b_1g_1^{-n_i}g_1\geq b_1g_1^{-n_i}g_i$ because the leading slope of $g_i$
is less or equal to the leading slope of $g_1$.  In particular, the whole
interval $(a,b_1)$ is in the first affine component of $h_i$.
\qquad$\diamond$

The following lemma is the achievement of the section, as it will
enable us to find the remarkable ``controllers''; elements that
control the global behavior of a balanced group on its
orbitals.  

\bl
\label{phiImageLemma}

Suppose $H$ is a balanced group with single orbital $A = (a,b)$, and
$\phi$ is the log-slope homomorphism defined before, then
$\phi(H)\cong \mathbb{Z}$ or $\phi(H)$ is trivial in $\rtr$.

\el 

pf: Let $H$ be a balanced subgroup of $\ploi$.  Each element of $H$
either realizes both ends of $A$, or neither.  In particular, the
group homomorphism $\rho_1:\R\times\R\to\R$ which is projection on the
first factor has the property that $ker(\rho_1)\cap{}Im(\phi) =
\{(0,0)\}$, the trivial subgroup of $\R\times\R$.  This implies that
the image of $\phi$ in $\R\times\R$ is isomorphic to the group
$H_1\leq\R$ obtained by considering only the first factors of elements
of $Im(\phi)$.  If no element realizes the ends of $A$, then $H_1 =
\left\{0\right\}$, the trivial (additive) group, and we are done.
Therefore, let us suppose instead that some elements in $H$ realize
the left end of $A$ (and therefore also the right) so that $H_1$
cannot be the trivial subgroup of $\R$.

If $H_1$ is discrete in $\R$ then $H_1$ is either trivial, or
isomorphic to $\Z$, but by assumption, $H_1$ is not the trivial group,
hence in this case $H_1\cong\Z$.  Hence, we shall suppose that $H_1$
is not discrete in $\R$.  In this case, by taking the difference of
two elements in the image which are very near each other, we see that
we can find an element of $H_1$ which is as close to zero as we like.
This implies there are elements of $H$ whose leading slopes are as
close to one as we like, without actually being one.  If $h$ is an
element with leading slope $s\neq 1$, then one of $h$ or $h^{-1}$ has
slope greater than one, since the leading slope of $h$ is $s$, but the
leading slope of $h^{-1}$ is $1/s$ (note that $s$ cannot be zero,
since no element of $\ploi$ has an affine component with slope zero).

Now suppose that $H$ is abelian.  By a result of Brin and Squier
\cite{picric}, if two elements in $\ploi$ commute and share a common
orbital, then their projections on that orbital have a common root.
If two elements have non-disjoint support and commute, it is easy to
see that the intersections of their supports actually is a set of
commonly shared orbitals. Now, let $h$ be some element of $H$ with
leading slope $s>1$.  For each positive integer $n$, let $g_n$ be an
element of $H$ with leading slope $s_n'$ where $1<s_n'<\frac{n+1}{n}$.
Now for each $g_n$, the pair $h$ and $g_n$ has a common root $h_n$ (on
their leading orbital), but infinitely many of the roots $h_n$ have
pairwise distinct leading slopes, since these slopes are always less
than or equal to the slopes of the $g_n$, and in particular, $h$ must
then have infinitely many distinct roots in $H$ on its leading
orbital.  By another result in \cite{picric}, no element of $\ploi$
has infinitely many distinct roots, so we must conclude that $H$ is
not abelian.

Note that by the details of the previous paragraph, it is easy to construct a
countably infinite sequence of elements $(g_n)_{n = 1}^{\infty}$ in $H$ which
satisfies the properties below: 
\be

\item For each $i\in\N$, the lead slope of $g_{i+1}$ is less than the
lead slope of $g_i$.

\item Given any real number $r>1$, there is an $i\in\N$ so that the
lead slope of $g_i$ is $s$ where $1<s<r$.

\item given $i$, $j\in\N$, we will have $[g_i,g_j]=1$ implies $i = j$.

\ee

Therefore, by Lemma \ref{longLeadSlope} there is $c\in(a,b)$ and
elements of $H$ with lead slopes that are greater than but arbitrarily
close to one, and whose leading affine components contain $(a,c)$.

Let $f$ and $g$ be two elements of $H$ with $h = [f,g]\neq 1$.  The
fixed set of $h$ in $I$ is disconnected, and contains two components of the
form $[0,u']$ and $[v',1]$ for some numbers $u'$, $v'\in (a,b)$.  In
particular $\inf(Supp(h)) = u'$ and $\sup(Supp(h)) = v'$.

By Lemma \ref{transitiveOrbital} there is $\alpha\in H$ so that
$v'\alpha<c$, so that $j=h^{\alpha}$ has all of its orbitals inside
$(a,c)$.  Let $u = \inf(Supp(j))$ and $v = \sup(Supp(j))$.  In
particular, $v<c$.

  We will now perturb $j$ slightly via a conjugation which wil move
the orbitals of $j$ to the right by a distance less than $L$, so that
$j$ and the new element together will generate a group with an
imbalanced orbital.  

Suppose $L>0$ is smaller than two particular lengths.  The first
length is the length of the second component of the fixed set of $j$
in $A$ which has non-zero length (note, this component might just be
$[v,b)$, if $j$ has only two such), and the second length is the
length of the first orbital of $j$.

Choose an element $\beta\in H$ whose leading affine component in $A$
contains $(a,c)$ and whose lead slope is greater than one, but so near
one that no point in $(a,c)$ will move to the right a distance greater
than $L$.  Now the elements $j$ and $j^{\beta}$ will generate a group
$G$ with leading orbital $(u,w)$ where $j$ realizes $u$ and possibly
$w$ (if the right ends of the appropriate orbitals of $j$ and
$j^{\beta}$ are aligned), while $j^{\beta}$ will achieve $w$ but not
$u$.  To see this, note that the left end of the first orbital of
$j^{\beta}$ is in the first orbital of $j$, so that $u$ is the left
end of the leading orbital of $G$, and only $j$ realizes it.
Meanwhile, the right end of the first orbital of $j^{\beta}$ is to the
right of the right end of the first orbital of $j$, so that the
leading orbital of $j^{\beta}$ contains the right end of the leading
orbital of $j$.  As we progress to the right, if the solitary fixed
point sets of $j$ and $j^{\beta}$ align before we reach the second
component of the fixed set of $j$ in $A$ with non-zero length, then
the right end of the first orbital of $G$ will be achieved by both $j$
and $j^{\beta}$.  Otherwise, the first orbital of $G$ will extend
rightward into the interior of the second fixed set of $j$ with a
non-empty interior, so that the right end of this orbital of $G$ will
be realized only by $j^{\beta}$.  In all cases, $G$ will be
imbalanced, and hence $H$ will also be imbalanced.  This contradicts
the hypothesies of the lemma, and so we see that $H_1$ must be
discrete in $\R$, and therefore the image of $\phi$ is isomorphic to
$\Z$ or the trivial group in $\R\times\R$.  \qquad$\diamond$

The kernal of the homomorphism $\phi$ is naturally very important as
well, it is the subgroup of $H$ which consists of elements which are
the identity near the ends of $A$.  Typically, we will refer to this
normal subgroup as ${}H^{\!\!\!\!^{\circ}}$.

\subsubsection{Controllers}
\label{controllers}

A consequence of section \ref{phi} is that the structure of a balanced
group with one orbital is very special.  In this section we will
explore this idea further.

\bl [Balanced Generator Existence]
\label{balancedGenerator}

Suppose that $H$ is a balanced subgroup of $\ploi$ with single orbital
$A$, that there is some element in $H$ which realizes an end of $A$,
and that ${}H^{\!\!\!\!^{\circ}}$ is the subgroup of $H$ which
consists of all elements in $H$ which are the identity near the ends
of $A$. Then there is an element $g$ of $H$ so that
$H=\left<\left<g\right>,{}H^{\!\!\!\!^{\circ}}\right>$, where $g$
realizes both ends of $A$.  
\el

pf: 

Let $\Gamma_A$ be the set of elements of $H$ which realize both ends
of $A$.  By our assumptions, $\Gamma_A$ is not empty.  Now observe
that $H = 
\left<\left<\Gamma_A\right>,\,{}H^{\!\!\!\!^{\circ}}\right>$.

By lemma \ref{phiImageLemma} the image $\phi(H)$ is cyclic in
$\R\times\R$.  Let $\gamma$ be a generator of the image of $\phi$.
Let $g$ be an element of $H$ so that $g\phi =\gamma$.  We observe that
since $\gamma$ is non-trivial in both components, $g$ realizes both
ends of $A$.  Since $\gamma$ generates the image of $\phi$, if
$\hat{g}\in{}\Gamma_A$, then $\hat{g}\phi = g^k\phi$ for some $k\in\Z$.
Hence, $g^{-k}\cdot\hat{g}\in{}H^{\!\!\!\!^{\circ}}$.  This now implies that
$\Gamma_A\subset \left<g,{}H^{\!\!\!\!^{\circ}}\right>$, so that
$\left<g,{}H^{\!\!\!\!^{\circ}}\right> = \left<\Gamma_A,{}H^{\!\!\!\!^{\circ}}\right> =
H$.\qquad$\diamond$

We will call any element $c$ of a balanced group $H$ with one orbital
$A$, which satisfies the rule $H =
\left<\left<c\right>,\,{}H^{\!\!\!\!^{\circ}}\right>$, a
\emph{controller}\index{controller} of $H$.  A controller of $H$ is
clearly a special element.

Given a controller $c$ of a balanced group $H$ with one orbital $A$,
we can write any element $h$ of $H$ uniquely in the form
$c^k\cdot g^{\!\!\!^{\circ}}$, where $k$ is some integer, and $g^{\!\!\!^{\circ}}\in{}H^{\!\!\!\!^{\circ}}$.
We will call this the \emph{$c$-form of $h$}\index{$c$-form}.

We will say that a controller $c$ of the group $H$ is
\emph{consistent}\index{controller!consistent} if and only if its
image $(\alpha,\,\beta) = c\phi$ satisfies the property that
$sign(\alpha) = -sign(\beta)$.  Otherwise we will say the controller
is \emph{inconsistent}\index{controller!inconsistent}.  The idea
behind this definition is that a one bump controller should be
consistent, since it is either moving points to the right everywhere
on its support, or it is moving points to the left everywhere on its
support.  An inconsistent controller must have a fixed point set, and
has at least one bump where the controller moves points to the right,
and one bump where the controller moves points to the left.  It turns
out that a consistent controller actually is a one-bump element of
$H$.

\bl
\label{cFullSupport}

Suppose $H$ is a balanced subgroup of $\ploi$ and $H$ has unique
orbital $A$.  Further suppose that $H$ has a consistent controller
$c$, then $c$ realizes $A$.

\el

pf:

Since $c$ and $c^{-1}$ are both controllers, and either both satisfy
or both fail the conclusion of the statement of the lemma, we will
assume $c$ moves points to the right on its orbitals near the ends of
$A$.  Suppose $c$ has a non-trivial fixed set $K$ in $A$.  $K$ is
closed and bounded and hence compact.  Let $u = \inf{K}$ and $v =
\sup{K}$, so that $(a,u)$ is the first orbital of $c$ and $(v,b)$ is
the last orbital of $c$.  Now there are points $x\in (a,u)$, and
$y'\in(v,b)$ so that we have $a<x<u\leq v<y'<b$.  By Lemma
\ref{transitiveOrbital} there is an element $g\in H$ so that $xg>y'$.
Writing $g$ in its $c$-form, we have that $g = c^kg^{\!\!\!^{\circ}}$ for some
integer $k$ and element $g^{\!\!\!^{\circ}}\in{}H^{\!\!\!\!^{\circ}}$.  In particular, the
element $h = gc^{-k}$ is trivial near the ends of $A$, but still
satisfies $xh=y>v$. The element $h$ therefore has an orbital $D=(r,s)$
which spans the fixed set $K$ of $c$.  Suppose $e=\inf(Supp(h))$, so
that $a<e\leq r$.  By Lemma \ref{transitiveElementOrbital} there is a
positive integer $N_1$ so that for any integer $n_1>N_1$ we will have
$r<ec^{n_1}<u$. Suppose $f=\sup(Supp(h))$, so that $s\leq f<b$.  By
Lemma \ref{transitiveElementOrbital} there is a positive integer $N_2$
so that for any integer $n_2>N_2$ we will have $f<sc^{n_2}<b$.  Let $n
= max(N_1,N_2)$, then, the element $j = h^{(c^n)}$ has its first
orbital starting at some interior point of $(r,s)$, and its orbital
induced from $(r,s)$ has right end $t$ which is strictly to the right
of $f$.  In particular the group $H_1=\left<j,h\right>$ has an orbital
$B = (r,t)$, where we note that $t>s$ by construction.  Now, $h$
realizes the left end of $B$ in $B$, but not the right, hence $H_1$,
and therefore $H$, is imbalanced.  But this contradicts our
assumptions, therefore $c$ must have orbital $A$.\qquad$\diamond$

\bc 
\label{consistentRealization}
Suppose $H$ is a balanced subgroup of $\ploi$ and $H$ has a unique
orbital \mbox{$A= (a,b)$}.  If $H$ has a consistent controller, $c$,
then any element $g$ of $H$ which realizes both ends of $A$ actually
realizes $A$.  
\ec 

pf: Let $g = c^kg^{\!\!\!^{\circ}}$ be the $c$-form of $g$, and suppose that $g$
has a non-empty fixed set in $(a,b)$.  The behavior of $g$ near the
ends of $A$ (beyond the support of $g^{\!\!\!^{\circ}}$) depends entirely on $c$
and the integer $k$, and therefore $g$ is either moving points to the
right on both of its leading and trailing orbitals, or moving points
to the left on both of its leading and trailing orbitals, since $k$
cannot be zero.  Now the details of the argument of Lemma
\ref{cFullSupport} show that $H$ must be imbalanced, which contradicts
our assumptions, hence $g$ must have orbital $A$.  \qquad$\diamond$

We will now consider the case where $H$ has an inconsistent
controller.

\brk
\label{iAntiRealization}
Suppose $H$ is a balanced subgroup of $\ploi$ and $H$ has unique
orbital \mbox{$A = (a,b)$}.  If $H$ has an inconsistent
controller $c$, then no element of $H$ realizes $A$.
\erk

pf:
Any element $h$ which realizes both ends must do so inconsistently.  Hence
the difference function $D:[a,b]\to\R$ defined by the rule $xD = xh-x$ is
a continuous function which is positive near one end of $A$ and
negative near the other, and so by the intermediate value theorem,
there is some point $x\in A$ so that $xD = 0$.  The point $x$
corresponds to a fixed point of $h$.\qquad$\diamond$

In the case where $H$ is a one-orbital balanced group with an
inconsistent controller, by Remark \ref{iAntiRealization} any element
$h$ which realizes both ends of $A$ must have a non-trivial fixed set
$F_h$ in $A$.  Suppose there exists an element $g^{\!\!\!^{\circ}}$ of ${}H^{\!\!\!\!^{\circ}}$
with an orbital $A_{g^{\!\!\!^{\circ}}}$, where the ends of $A_{g^{\!\!\!^{\circ}}}$ are in
differing orbitals of $h$.  The element $g^{\!\!\!^{\circ}}$ enables the transfer
of a point $x\in{}Supp(h)\cap{}A_{g^{\!\!\!^{\circ}}}$ from its own orbital in
$h$ to another orbital of $h$.  In particular, the point $xg^k$ will
be in a different orbital of $h$ than the orbital which contains $x$,
for some integer power $k$.  Such an element $g^{\!\!\!^{\circ}}$ is very useful
to have in hand, so we will study questions about the existence and
structure of such elements in a one-orbital balanced group with
inconsistent controller.  

\bl
\label{fixedCoverLemma}
Suppose $H$ is a balanced subgroup of $\ploi$ and $H$ has a unique
orbital \mbox{$A = (a,b)$}.  Suppose that $H$ has an inconsistent
controller $c$ and let ${}H^{\!\!\!\!^{\circ}}$ represent the normal
subgroup of $H$ consisting of the elements in $H$ which achieve
neither end of $A$.  Now, the fixed set $F_c=fix(c)\cap A$ of $c$ in
$A$ is contained in an orbital of some element
$g^{\!\!\!^{\circ}}\in{}H^{\!\!\!\!^{\circ}}$.  \el

pf: Let $x = \inf{F_c}$ and let $y=\sup{F_c}$.  Note that $x$, $y\in
F_c$\!.  There is an element $g\in{}H$ so that $xg>y$ by lemma
\ref{transitiveOrbital}.  We can write $g=c^k\cdot g^{\!\!\!^{\circ}}$ for some
integer $k$ and some element $g^{\!\!\!^{\circ}}\in{}H^{\!\!\!\!^{\circ}}$.  Clearly, $c^{-k}$
does not move the fixed set of $c$, but $g$ moves a closed interval
containing the fixed set of $c$ completely off of itself.  In
particular, we have that $xc^{-k}g>y$, so that $F_c\subset [x,y]
\subset Supp(c^{-k}g)$.  But $c^{-k}g = g^{\!\!\!^{\circ}}\in{}H^{\!\!\!\!^{\circ}}$, so
$F_c\subset{}B$ for some orbital $B$ of $g^{\!\!\!^{\circ}}$.  \qquad$\diamond$

\bc [Transfer Existence]
\label{transferExistence}

Suppose $H$ is a balanced subgroup of $\ploi$ and $H$ has a unique
orbital \mbox{$A = (a,b)$}.  Suppose that $H$ has an inconsistent
controller $c$ and let $K$ be any compact set in $(a,b)$.  Then $K$
is contained in an orbital of some element
$g^{\!\!\!^{\circ}}\in{}H^{\!\!\!\!^{\circ}}$.  \ec

pf:

Suppose that $c$ moves points to the left on its first orbital, and
moves points to the right on its last, by replacing $c$ with it's
inverse, if necessary.

Let $x=\inf{K}$ and $y = \sup{K}$.  Let
$\hat{g}\in{}H^{\!\!\!\!^{\circ}}$ be an element which contains
the fixed set of $c$ in some orbital $\hat{A}_c$.  By the
constraints that $c$ moves points to the left on its first orbital,
that $c$ moves points to the right on its last orbital, and that the
conjugate $c^{-1}\hat{g}c$ has the induced orbital corresponding
to $\hat{A}_c$ larger than $\hat{A}_c$.  By conjugating
repeatedly with $c$, we can make the ends of this orbital approach the
ends of $A$, so that after some conjugation the resulting orbital will
contain the set $[x,y]$.  The resulting function is
$g^{\!\!\!^{\circ}}\in{}H^{\!\!\!\!^{\circ}}$.  \qquad$\diamond$

Suppose that we know that $A$ is an orbital of a balanced group $H$,
and that the ends of $A$ are achieved by some element $g$ in $H$.  Let
$H_A$ be the set of elements of $\ploi$ each of which is equal to the
restricition of some element of $H$ on the orbital $A$, and behaves as
the identity off of $A$.  $H_A$ is trivially a group with unique
orbital $A$, and is a quotient of $H$.  We will call $H_A$ the
projection of $H$ on $A$.

We will now generalize our language somewhat.  Let $H$ be a subgroup
of $\ploi$ with an orbital $A$, and let $H_A$ be the projection of $H$
on $A$.  $H_A$ has a controller $\tilde{c}$ for $A$.  Let
$\rho:H\to H_A$ be the projection homomorphism on the orbital $A$.
Let $T = \left\{c\in H\,|\,c\rho \textrm{ is a controller for $H_A$ on
$A$}\right\}$.  We will call any element of $T$ \emph{a controller of
$H$ on $A$}\index{controller!on an orbital}.  Again, given an element
$c\in H$ which is a controller of $H$ on $A$, we can write elements of
$H$ in a unique $c$-form.  I.e., if $g\in H$, and $c$ is a controller
for $H$ on $A$, then there is an integer $k$ so that $g =
c^kg^{\!\!\!^{\circ}}$, where $g^{\!\!\!^{\circ}}$ will not realize
either end of $A$.

\subsection{Towers\label{towers}}
We define and study objects called towers which will give us a
geometric criterion on derived length.

Given a group $G\leq\ploi$ and a set $T$ of signed orbitals of $G$, we
will say $T$ is a \emph{tower of $G$}\index{tower!of group} if $T$
satisfies the following properties:

\be

\item $T$ is a chain in the partial order on the signed orbitals of $G$.

\item For any orbital $A\in O_T$, $T$ has exactly one element of the
  form $(A,g)$.

\ee

We note that every group $G\in\ploi$ has the empty tower as one of its
towers.  We will also sometimes create a tower from a given chain in
the orbitals of the elements of $G$, where we will implicitly use a
choice function from the set of orbitals of elements of $G$ to $G$, so
that each orbital maps to some element of $G$ which has that orbital.
If a non-empty tower is finite, then it will admit an order preserving
bijection from a set $\left\{1,2,\ldots,n\right\}$ for some positive
integer $n$.  We will therefore refer to a finite tower's ``$i$-th''
element, by which we mean the image of the integer
$i\in\{1,2,\ldots,n\}$ in the tower, extending this, we will also
sometimes refer to the tower's ``smallest'' element, and its ``largest''
element.

It is a trivial observation that if $T$ is a tower of $G$ for some
subgroup $G\in\ploi$ with element $(A,g)$, then $T$ has no other
element of the form $(B,g)$.

With the definition of tower in place, we have a natural measure of
one form of complexity for a group.  Given a group $G\leq\ploi$, we
define the \emph{depth of $G$}\index{group!depth} to be the supremum
of the set of cardinalities of the towers of $G$.  Note that groups
have depth while towers have height.
 
We observe that the depth of a group $G\leq \ploi$ is well defined.
The set of all towers of $G$ is a nonempty (it contains the empty
tower) subset of the power set of the cartesian product of the set of
all orbitals of elements of $G$ with the underlying set of the group
$G$.  The set of all cardinalities of the set of all towers of $G$ is
therefore a nonempty subset of the ordinals, and this set has a
supremum.  We observe that the depth of the trivial group is zero,
and that the depth of any non-trivial group is greater than zero.  In
fact, the depth of all abelian groups is one, as we will see shortly.

The following is a good exercise to help the reader become familiar
with these concepts.  In exploring this, it is helpful to recall that
$W_i$ is a balanced group (it is solvable with derived length $i$, but
$F$ is known to be non-solvable, so $W_i$ cannot contain a copy of
$F$) and to observe that $\alpha_1$ is actually a consistent
controller for $W_2$ on the orbital $(0,1)$.

\brk 
The depth of $W_2$ is two.  
\erk 

where here we mean the realization of $W_2$ given in the introduction.

We extend the notions of depth and height to other objects.  

Let $T$ be a tower of $\ploi$.  We will call the cardinality of $T$
its \emph{height}\index{tower!height}, using the simple descriptive
\emph{infinite}\index{tower!infinite} if $T$ has an infinite
cardinality.  If there is an order preserving injection from $\N$ to
$T$, then we will say $T$ is \emph{tall}\index{tower!tall}, and if
there is an order preserving injection from the negative integers to
$T$, then we will say $T$ is \emph{deep}\index{tower!deep}.  If $T$ is
both deep and tall then we will say $T$ is a \emph{bi-infinite
tower}\index{tower!bi-infinite}, and we note that there will be an
order preserving injection from the integers to the tower.  We will
occassionally replace an infinite tower of one of these three types
with the image of the implied injection without comment, so that we
might refer to a tall, deep, or bi-infinite tower as countable, and
refer to the ``next'' element, etc., when this will not effect the
result of an argument.

Given a group $G\leq \ploi$, and an orbital $A$ of an element $g\in
G$.  We will define the \emph{depth of $A$ in
$G$}\index{orbital!depth} to be the supremum of the heights of the
finite towers which have their smallest element having the form
$(A,h)$ for some element $h\in G$.  If the depth of $A$ is an infinite
ordinal, we will simply say that $A$ is \emph{deep in
$G$}\index{orbital!deep}.  Symmetrically, we define the \emph{height
of $A$ in $G$}\index{orbital!height} to be the supremum of the heights
of the finite towers which have their largest element having the form
$(A,h)$ for some element $h\in G$.  If the height of $A$ is an
infinite ordinal, we will simply say that $A$ is \emph{high in
$G$}\index{orbital!high}.

The following are immediate from the definitions.
\brk
\label{trivialTower}

\be

\item Any subset of a tower is a tower.
\item If $T$ is a tower for some group $G\leq\ploi$, and $(A,g)\in T$,
and if $h$ is an element of $G$ with orbital $A$, then the set
$(T\backslash\left\{(A,g)\right\})\cup\left\{(A,h)\right\}$ is also a
tower for $G$.
\item If $T$ is a tower for some group $G\leq\ploi$, and $(A,g)\in T$,
and if $n$ is a non-zero integer, then
$(T\backslash\left\{(A,g)\right\})\cup\left\{(A,g^n)\right\}$ is also a
tower for $G$.
\item If $g\in G\leq\ploi$, and $T$ is a tower of $G$, and if $(A,g)$,
$(B,g)\in T$, then $A=B$.  That is, no signature appears twice in a tower.
\item Given a tower $T$ of a group $G\leq\ploi$, the group $H\leq G$
generated by the signatures of $T$ has an orbital $A$ which contains
all the orbitals of $T$.  
\item \label{conjugateTower} Given an element $k\in G\leq\ploi$ and any
tower $T$ for $G$, the set of signed orbitals
$T^k=\left\{(Ak,g^k)|(A,g)\in T\right\}$ is also a tower for $G$,
where the natural order of the signed orbitals of $T^k$ is equal to
the induced order from the signed orbitals of $T$.  
\ee 
\erk 

Given $k\in G\leq\ploi$ and a tower $T$ for $G$, the tower $T^k$
induced from the tower $T$ by the action of $k$ as discussed in item
\ref{conjugateTower} of Remark \ref{trivialTower} will be called the
\emph{tower conjugate to $T$ by the action of
$k$}\index{tower!conjugate}.  We will also say that the towers are
conjugate towers.  Conjugacy of towers for $G$ is an
equivalence relation on the set of towers for $G$.

Towers, as defined, are easy to find, but can be difficult to work
with.  For an arbitrary tower $T$, there are no guarantees about how
other orbitals of signatures of the elements of $T$ cooperate with the
orbitals of the tower.  We say a tower $T$ is an \emph{exemplary
tower}\index{tower!exemplary} if the following two additional
properties hold: \be

\item Whenever $(A,g)$, $(B,h)\in T$ then $(A,g)\leq(B,h)$ implies the
orbitals of $g$ are disjoint from both ends of the orbital $B$.

\item Whenever $(A,g)$, $(B,h)\in T$ then
$(A,g)\leq(B,h)$ implies no orbital of $g$ in $B$ shares an end with $B$.

\ee

The following three lemmas are an indication of the plethora of
exemplary towers in $\ploi$, and will be used repeatedly.

\bl
\label{imbalancedTower}

Suppose $H\leq\ploi$, and that $G\leq H$ has imbalanced orbital $A =
(a,b)$.  Then $G$ admits an exemplary bi-infinite tower $E$ whose
orbitals are all in $A$.

\el

pf:

A short incorrect proof is that $H$ contains a copy of $F$ and $F$
contains an exemplary bi-infinite tower.  The problem is that Theorem
\ref{UbiquitousF} guarantees a copy of $F$, but does not guarantee
that its generators have one orbital each.

If we find an exemplary bi-infinite tower $E'$ for the projection
$G_A$ of $G$ on $A$, then by replacing the signatures of $E'$ with
elements of $G$ which agree with the signatures of $E'$ on $A$, we can
build a new exemplary tower $E$ for $G$ with the same orbitals.  Thus,
we may assume for the purposes of this argument that $G$ only has
orbital $A$.

Since $A$ is imbalanced for $G$, there is $g_0\in G$ so that $g_0$ has
an orbital $B_0$ which shares an end with $A$, but $g_0$ does not
realize the other end of $A$.  We will assume that $B_0$ shares its
right end with the right end of $A$, in particular, $B_0 = (a_0,b)$
for some $a_0\in A$, and there is $w\in A$ so that if $x\in Supp(g_0)$
then $w<x$.  We will further assume that $g_0$ moves points to the
right on the orbital $B_0$, so that all conjugates of $g_0$ move
points to the right on their corresponding orbitals.

By Lemma \ref{transitiveOrbital} there is an element $\alpha\in G$ so
that $r = w\alpha>a_0$.  Let $g_k = g_0^{(\alpha^{-k})}$ for all
integers $k$.  By construction, given any integer $k$, the element
$g_k$ has rightmost orbital of the form $(a_k,b)$, where
$a_{k+1}<\inf(Supp(g_k))$.

We will construct a new sequence from the $g_k$ so that the righthand
orbitals of the new sequence do not get arbitrarily small, while
preserving all of the nesting properties we will need later.  Fix
$l\in (a_0, r)$.  

We now have $a<w<a_0<l<r<b$.  Now for each non-negative integer $k$,
we already have that $a_k<l$, so define $h_k = g_k$ when $k\geq 0$.
Since $h_0$ is already defined, we can define $h_{-k}$ inductively for
each $k\in\N$ as the first conjugate of $g_{-k}$ by $h_{-k+1}^{-1}$
which has the property that the rightmost orbital $(c_{-k},b)$ of
$h_{-k}$ contains $(l,b)$.  For each integer $k$, let $c_k$ represent
the left end of the rightmost orbital of $h_k$.  In particular, we
have now defined a bi-infinite sequence of functions $(h_k)_{k \in
\Z}$ that satisfies the following properties:

\be

\item Given $k\in\Z$, the rightmost orbital of $h_k$ is $(c_k,b)$.

\item Given $k\in\Z$, $(l,b)\subset (c_k,b)$.

\item Given $k\in\Z$, $c_{k+1} < (\inf(Supp(h_k))$.

\ee

Note that the support of the element $g_{-1}$ lies to the right of
$r$.  For all $k\in\Z$, define $u_k = h_kg_{-1}^{-1}$.  Since the
$h_i$ are all conjugates of $g_0$, and $g_{-1}$ is a conjugate of
$g_0$, we see that the trailing slopes of all of the $h_k$ are the
same.  In particular, for each integer $k$, the fixed set of $u_k$ has
a component of the form $[e_k,b)$ where $l<r<e_k<b$.  Furthermore, for
each integer $k$, $u_k$ has an orbital $C_k$ of the form $(c_k,d_k)$
where $r<d_k\leq e_k<b$.  There may be other orbitals of $u_k$ to the
right of $d_k$ and the $d_i$ may not be ordered with respect to $i$.

Let $v_0 = u_0$, and inductively define, for each positive integer
$k$, an element $v_k = u_k^{h_k^{m_k}}$ where $m_k$ is the smallest
positive integer so that the support of $v_{k-1}$ is fully contained
in the orbital $D_k = (c_k,r_k)$ of $v_k$ induced from
$C_k=(c_k,d_k)$.  Note that such an integer $m_k$ will always exist,
since $h_k$ moves points to the right on its rightmost orbital
$(c_k,b)$, and by reference to lemma \ref{transitiveElementOrbital}.

Now, for each positive integer $k$, inductively define elements
$v_{-k}$ by following the following two step process.

First, find the negative integer $n_k$ of smallest absolute value
which has the property that the closure of the support of $v_{-k}' =
u_{-k}^{(h_{-k}^{n_k})}$ is fully contained in the orbital $D_{-k+1}$
of $v_{-k+1}$ induced from $C_{-k+1}=(c_{-k+1},d_{-k+1})$.  At this
point $v_{-k-1}$ cannot be created since the orbital of $v_{-k}'$
induced from $C_{-k}$ might not contain the left end of the orbital
$(c_{-k-1},b)$ of $h_{-k-1}$.  

Now replace the elements of the sequences $(u_j)_{j\in\Z}$ and
$(h_j)_{j\in\Z}$ which have indices less than or equal to $-k$ by the
conjugate of each such element by $h_{-k}^{n_k}$.  Note $n_k<0$.
This will do nothing to the element $h_{-k}$ of the sequence
$(h_j)_{j\in\Z}$, but all terms of $(h_j)_{j\in\Z}$ with $j<-k$ will
have their orbitals extended leftward by $|n_k|$ iterates of
$h_k^{-1}$.  Now define $v_{-k} = u_{-k}$.  At any inductive stage
$k$, note that by construction, all of the left ends of the orbitals
of the $h_j$, for indices $j<-k$, are inside the orbital $D_{-k}$, so
this inductive definition makes sense.  For each integer $k$, the
closure of the support of $v_{k-1}$ is a subset of this orbital $D_k$.  In
particular, $E = \left\{(D_k,v_k)|k\in\Z\right\}$ is an exemplary
bi-infinite tower for $G$, with all orbitals $D_k$ in $A$, and where
the index of $E$ respects the natural order on the signed orbitals of
$E$.  \qquad$\diamond$

In the balanced group case, we have another way to sometimes find
bi-infinite, exemplary towers: 

\bl
\label{inconsistentTower}

Suppose $H$ is a balanced subgroup of $\ploi$ and $H$ has orbital $A$.
If $h$ is an element of $H$ so that $h$ realizes the ends of $A$
inconsistently, then $H$ admits an exemplary bi-infinite tower $T$ with
the orbitals of $T$ all contained in $A$.

\el
pf:

We will find such a tower for the projection of $H$ on $A$, then
replacing the signatures of our tower with elements of $H$ which agree
with our signatures over $A$ will create an exemplary tower for $H$
whose orbitals are all in $A$.  Thus, we shall assume that
$H$ has only the orbital $A$ for our discussion below.

Let $h\in H$ so that $h$ realizes the ends of $A$ inconsistently.  By
replacing $h$ with its inverse, if necessary, we can assume that $h$
moves points to the left on its first orbital and moves points to the
right on its last.  Let $F_h=fix(h)\cap A$ represent the (non-empty)
fixed set of $h$ in $A$.

By lemma \ref{transferExistence}, there is an element $g_0$ in
${}G^{\!\!\!\!^{\circ}}$ which has an orbital $B_0$ which fully
contains the set $F_h$.  Let $r = \inf(F_h)$ and $s=\sup(F_h)$, so
that $F_h\subset[r,s]\subset B_0\subset\bar{B}_0\subset A$.  We may
assume $g_0$ moves points to the right on $B_0$.  There is a smallest
positive integer $n_1$ so that the orbital $B_1$ of $g_1 =
g_0^{(h^{n_1})}$ induced from $B_0$ by the action of $h^{n_1}$
contains the closure of the support of $g_0$, since repeated
conjugation of $g_0$ by $h$ increases the size of the orbital $B_0$ so
that the ends of $B_0$ approach the ends of $A$.  In particular, we
can inductively define a bi-infinite sequence $((B_i,g_i))_{i\in\Z}$
of signed orbitals of $H$ by the property that $g_i =
g_{i-1}^{(h^{n_i})}$ where $n_i$ is the smallest positive integer so
that the closure of the support of $g_{i-1}$ is fully contained in the
orbital $B_i$ induced from the orbital $B_{i-1}$ by the action of
$h^{n_i}$ (more formally, for negative integers $i$, one defines $g_i$
from $g_{i+1}$ by saying that $g_i = g_{i+1}^{(h^{-n_{i+1}})}$ where
$-n_{i+1}$ is the largest negative integer so that the closure of the
support of $g_i$ is contained in $B_{i+1}$ where this last makes sense
since repeated conjugation by $h^{-1}$ moves points in the first and
last orbitals of $h$ arbitrarily close to the closed interval
$[r,s]$).

$T= \left\{(B_i,g_i)|i\in\Z\right\}$ is now an exemplary bi-infinite
tower for $H$ whose indexing follows the natural ordering of the
signed orbitals of $T$.
\qquad$\diamond$

In the following recall that $S_T$ is the set of signatures of a tower $T$.
\bl
\label{finiteExemplaryTower}
Suppose $G$ is a balanced group and $G$ has a tower $T$
of height $n$ for some positive integer $n$.  If $B_n$ is the orbital of
$H_n=\left<S_T\right>$ which contains the orbitals of $T$, then $G$ has
an exemplary tower $E$ of height $n$ whose orbitals are contained in $B_n$.

\el

pf: 

Let us write $T = \left\{(A_i,g_i)|1\leq i \leq n, i\in\N\right\}$
where the indexing of $T$ respects the natural ordering of the signed
orbitals of $T$.

For each positive integer $k$ in $\left\{1,2,\ldots,n\right\}$, define
the following objects:
\[
\begin{array}{c}
T_k = \left\{(A_i,g_i)|1\leq i\leq k, i\in\N\right\},
\\
H_k =\left<S_{T_k}\right>.
\end{array}
\]
Furthermore, for each such $k$, let $B_k$ be the orbital of $H_k$
which contains the orbitals of the tower $T_k$ for $H_k$, and note
that each end of $B_k$ is realized by some signatures
from $T_k$, since $H_k$ is finitely generated.  Further, since $G$
is balanced, each signature of $T_k$ which realizes one end of $B_k$
must also realize the other, and so for each $k$, some signature of
$T_k$ must realize both ends of $B_k$ in $B_k$.  But now, by Lemma
\ref{balancedGenerator} we know that the projection of the group $H_k$
on the orbital $B_k$ has a controller $c_k$ for the orbital $B_k$.  

There are now two cases to analyze. 

First, if $c_k$ is inconsistent, then the signature of $T_k$ which
realizes both ends of the orbital $B_k$ must do so inconsistently, and
so, by Corollary \ref{inconsistentTower}, $H_k$ will admit an
exemplary bi-infinite tower $F$ whose orbitals all are contained in
the orbital $B_k\subset B_n$.  $F$ will admit an exemplary subtower
$E$ of height $n$ all of whose orbitals will be contained in $B_n$,
and we would be finished.

Second, if $c_k$ is consistent, then the signatures of $T_k$ which
realize both ends of $B_k$ do so in a consistent manner.  But now, by
Lemma \ref{cFullSupport}, each signature of $T_k$ that realizes both
ends of $B_k$ actually realizes $B_k$, and at least one signature does
so, so that signature must be $g_k$, since $A_k$ properly contains the
orbitals $A_i$ with $i<k$.  In particular, $A_k = B_k$ for each
integer $k\in\left\{1,2,\ldots,n\right\}$.

Given any orbital $B$ of some signature $g_i$ of $T_n$, where $B\neq
A_i$, we see that there is no index $k$ where $B$ contains an end of
any orbital $A_k$ of $T$.  Otherwise, we see that $k > i$, so that the
orbital $B_k$ of $H_k$ will be realized inconsistently by some signature $g_r$
with $1\leq r\leq k$, which we have already assumed does not happen.
In a similar fashion, $B$ cannot share an end with any orbital $A_k$
in the orbital $A_k$. Otherwise, we again have $k>i$.  Now $g_i$ and
$g_k$ generate a balanced group with orbital $C$ containing $A_k$ and
$B$ and having their commonly shared end as an end, and therefore both
$g_i$ and $g_k$ share the other end of $C$ as well.  At the same time,
since $g_i$ cannot realize $C$ ($A_i\cup B\subset C$), we must have
that $g_i$ realizes the two ends inconsistently, which we have already
ruled out.  Therefore, $E=T = T_n$ is already exemplary.
\qquad$\diamond$

The following indicates that subgroups of $\ploi$ without transition
chains of length greater than one are structurally much less complex
than general subgroups of $\ploi$.

\bl
\label{transitionTower}

If $G$ admits a transition chain of length two, then $G$ admits an
exemplary bi-infinite tower.

\el

pf: 

If $G$ is imbalanced, then $G$ admits an exemplary bi-infinite tower
by Lemma \ref{imbalancedTower}, and we are finished, so let us
suppose instead that $G$ is balanced.

Suppose $G$ admits a transition chain $\mathscr{C}'
=\left\{(B_1,g_1),(B_2,g_2)\right\}$.  $G$ now must admit a maximal
transition chain $\mathscr{C} = \left\{(A_i,h_i)|1\leq i\leq n,
i\in\N\right\}$ using only the signatures $g_1$ and $g_2$, where $n$
is some integer greater than one.  Define $H =
\left<S_{\mathscr{C}}\right>$.  $H$ is a finitely generated subgroup
of $G$ with orbital $A$ containing the orbitals $A_i$.  Now one of
$g_1$ and $g_2$ must realize the left end of $A$, and since $H$ is
balanced, it must also realize the right end of $A$.  We will assume
that $g_1$ realizes the ends of $A$ for purposes of
discussion, since the other case is completely symmetric.  But now,
since $g_1$ does not realize $A$, it must realize both ends of $A$
inconsistently, by Corollary \ref{cFullSupport}, and therefore by
Lemma \ref{inconsistentTower}, $H$ (and therefore $G$) admits an
exemplary bi-infinite tower all of whose orbitals are in $A$.
\qquad$\diamond$

\bc
\label{infiniteExemplaryTower}

If $G$ is a balanced subgroup of $\ploi$ and $G$ admits a tall tower in
some orbital $A$, or $G$ admits a deep tower in some orbital $A$, or $G$
admits a bi-infinite tower in some orbital $A$, then $G$ admits an
exemplary tall tower in $A$, or $G$ admits an exemplary deep tower in
$A$, or $G$ admits an exemplary bi-infinite tower in $A$,
respectively.

\ec 

pf: 

This follows easily from the details of the proof of Lemma
\ref{finiteExemplaryTower}.  Suppose $T$ is an infinite non-exemplary
tower for a balanced subgroup $G$, where all the orbitals of $T$ are
contained in the orbital $A$ of $G$.  Since $T$ is not exemplary, then
we can produce a non-exemplary subtower $P = \left\{(A_1,g_1),
(A_2,g_2)\right\}$ of $T$ where $A_1\subset A_2$.  Suppose $B$ is the
orbital of $G_P = \left<g_1,g_2\right>$ that contains $A_1$ (and
therefore $A_2$, also note that $B\subset A$).  Since $P$ is not
exemplary, we must have that some orbital of $g_1$ contains an end of
$A_2$ or shares an end of $A_2$ in $A_2$.  In the first case, $G_P$
admits a bi-infinite exemplary tower in $B$ by Lemma
\ref{transitionTower}, which has a subtower of the appropriate type.
In the other case, $A_2 = B$, so that $A_2$ is actually an orbital of
$G_P$.  Now since we are in the second case, we must have that $g_1$
has some orbital $C$ in $A_2$ that shares one end of $A_2$ in $A_2$.
Since $G_P$ is balanced, we have that $g_1$ realizes both ends of
$A_2$ from within $A_2$.  But now $g_1$ realizes $A_2$ since $G$ is
balanced and the orbital $A_2$ is realized consistently by $g_2$.  But
this means that $A_2=A_1$, which contradicts the fact that $A_1\neq
A_2$.  \qquad$\diamond$

We note in passing that it is an open question as to whether there are
any finitely generated, balanced subgroups of $\ploi$ which do not
admit transition chains of length two, but which are non-solvable.
This can be considered one of the main remaining geometric questions
in the theory of solvability of subgroups of $\ploi$.  Such a group
would be very interesting to examine.

\subsection{Derived groups and towers}
The following lemma represents the key concept for understanding why
solvability and depth are connected.  Note that here and in the
remainder we use $[a,b] = a^{-1}b^{-1}ab$.

\bl
\label{shortTower}

Suppose $G$ is a subgroup of $\ploi$, $n$ is a positive integer, and
$A$ is an orbital of $G$.  If $G$ has a tower of height $n$ whose
orbitals are contained in $A$, then $G'$ has a tower of height $n-1$
whose orbitals are also contained in $A$.

\el 

pf: Suppose $T=\left\{(A_i,g_i)|1\leq i \leq n, i\in\N\right\}$ is a
tower of height $n$ all of whose orbitals are contained in the orbital
$A$ of $G$ and whose indexing respects the order of the
elements of $T$.  By Lemmas \ref{imbalancedTower} and
\ref{finiteExemplaryTower}, we can assume that $T$ is an exemplary
tower.

For each integer $i$ with $2\leq i \leq n$ there is a smallest
positive integer $n_i$ so that the subset of the support of $g_{i-1}$
which is in $A_i$ is fully contained in a fundamental domain of
$g_i^{n_i}$ in $A_i$, by Lemma \ref{transitiveElementOrbital}.  Define
$h_1 = g_1$, and for each integer $i$ with $2\leq i \leq n$, define
$h_i = g_i^{n_i}$.  Now define the set $E = \left\{(A_i,h_i)|1\leq
i\leq n, i\in\N\right\}$, which is a exemplary tower.  $E$ has a nice
property: if $i$ is an integer in $1\leq i < n$, then the supports of
$h_i$ and the supports of $h_i^{h_{i+1}}$ are disjoint in $A_{i+1}$, so
that $A_i$ is an orbital of $[h_i,h_{i+1}]$.  For
each integer $i$ in $1\leq i <n$ define $v_i =
[h_i,h_{i+1}]$.  Noting that $v_i\in G'$, we see that
$\left\{(A_i,v_i)|1\leq i \leq n-1, i\in\Z\right\}$ is a tower of
height $n-1$ for $G'$ whose orbitals are all in $A$.  \qquad$\diamond$

\subsection{Geometric classification of solvable groups in $\ploi$}
We are now in a position to produce algebraic results by using our
geometric tools.  

The following is an immediate consequence of Lemma \ref{shortTower}

\bc

Suppose $G$ is a subgroup of $\ploi$.  If $G$ has towers of arbitrary
height then $G$ is non-solvable.

\ec 

The following is the key lemma towards building a geometric understanding of
derived groups in solvable groups of $\ploi$

\bl
\label{solvableClassification}

If $G$ is a subgroup of $\ploi$ of depth $n$ for some positive
integer $n$, then $G'$ is a subgroup of $\ploi$ of depth $n-1$.

\el 

pf: 

Suppose $G$ is a subgroup of $\ploi$ with depth $n$.  $G$ must
be balanced by Lemma \ref{imbalancedTower}, and $G$ must have no
transition chains of length two by Lemma
\ref{transitionTower}.

Since $n>0$, $G$ is not the trivial group, and so $G$ has at least one
orbital.  Let $A = (a,b)$ be an orbital of $G$.  Note that if 
\[
\Upsilon =
\left\{B\,|\,B\textrm { is an orbital of some element of G}, B\subset
A\right\},
\]
then $A = \cup_{B\in\Upsilon}B$.  But $G$ has no transition chains of
length two or the depth of $G$ would be infinite, so $A$ can be
written as a union of a chain of properly nested orbitals of elements
of $G$.  Taking these orbitals, paired with appropriate signatures, we
create a tower $T$ whose height is bounded above by $n$.  Let the
height of $T$ be $m$, and let $T = \left\{(A_i,g_i)\,|\, 1\leq i\leq
m, i\in\Z\right\}$.  But now, by construction, $A$ is the union of the
orbitals $A_i$, all of which are contained in $A_m$, so that $A_m$
must actually be $A$, and no element $g\in G$ has an orbital $B$
properly containing $A_m$ (or $A = A_m$ would not be an orbital of
$G$).  Thus, $(A_m,g_m)$ is a signed orbital of depth one, and the
orbitals of $G$ are precisely the orbitals of depth one for $G$. (Note
that since $G$ admits no transition chains, any signed orbital of
depth one for $G$ is automatically an orbital of $G$.)

Let $g$, $h\in G$, and consider the element $[g,h] \in
G'$.  Let $\Gamma$ be the set of all orbitals of $g$ and $h$.  Suppose
$[g,h]$ has $m$ orbitals, where $m$ is some positive integer, and let
$\left\{A_i|1\leq i\leq m, i\in\N\right\}$ be the orbitals of $[g,h]$
in left to right order.  Both $g$ and $h$ fix $a$, so the slope of the
leftmost affine component of $[g,h]$ that intersects $A$ is the
product of the slopes of the leftmost affine components of $g^{-1}$,
$h^{-1}$, $g$, and $h$ in $A$, which product is one.  In particular,
$[g,h]$ cannot realize $A$, so no orbital of $[g,h]$ is an orbital of
$G$, and so no tower for $G'$ contains an orbital of depth one for
$G$, and thusly, all towers of $G'$ can have height at most $n-1$.

By the last paragraph, we see that the depth of $G'$ is at most $n-1$.
By Lemma \ref{shortTower} $G'$ has a tower of depth $n-1$, so the
depth of $G'$ is actually $n-1$.
\qquad$\diamond$

We are now ready to prove our main geometric result, that given a
non-negative integer $n$, a subgroup $H$ in $\ploi$ is solvable with
derived length $n$ if and only if $H$ has depth $n$.  

\emph{Proof of Theorem \ref{geoClassification}:}

If the derived length of a group is $n$, then it must contain a tower
of height at least $n$, otherwise by Lemma
\ref{solvableClassification} and the fact that a depth zero group is
trivial, the derived series will terminate too soon.  But if $G$ has a
tower of height greater than $n$, then again by Lemma
\ref{solvableClassification}, the $n$-th derived group $G^{(n)}$ of
$G$ will admit a tower of height at least one, and so $G^{(n)}$ will
not be trivial.  \qquad$\diamond$

The next lemma is a technical lemma that we will use in completing our
proof of Theorem \ref{solveClassification}.

\bl
\label{oneBumpGenerators}

If $G$ is a solvable subgroup of $\ploi$ with derived length $n$,
generated by a collection $\Gamma$ of elements of $\ploi$ which each
admit exactly one orbital, and where no generator can be conjugated by
an element of $G$ to share an orbital with a different generator, then
$G$ is isomorphic to a group in the class $\left\{1\right\}\pc$ with
derived length $n$.

\el

pf: 
Before getting into the main body of the proof, note that the
hypothesies imply that each generator is the only generator with that
orbital, and that no element orbital is a union of element orbitals
that do not realize the original orbital.

We now enter the main body of the proof.  We will proceed by induction
on $n$.

If $n = 0$ then $G$ is the trivial group, and $G\in
\left\{1\right\}\pc$.  If $n = 1$ then $G$ is abelian, and in
particular, there can be at most countably many generators in
$\Gamma$, all of which have disjoint support, so that $G$ is
isomorphic with a countable (or finite) direct sum of $\Z$ factors.

Now let us suppose that $n > 1$ and that the statement of the lemma is
correct for any such solvable group with derived length $n-1$.  Let
$X$ represent the generators in $\Gamma$ whose orbitals are all depth
$2$ or deeper.  Let $Y$ be the set of elements in $\Gamma$ whose
orbitals have depth $1$.  We note in passing that the cardinality of
$Y$ is at most countably infinite, and that the collection of orbitals
of the elements of $Y$ actually form the orbitals of the group $G$.
We will assume that all of the elements in $Y$ move points to the
right on their orbitals.  We can partition the elements in $X$ into
sets $P_y$, where the $y$ index runs over the elements in $Y$, and
where an element of $X$ is in $P_y$ if and only if that element's
orbital is contained inside the orbital of $y$.  Given $y\in Y$, if
$P_y$ is empty, then define $H_y = \langle y \rangle$.  Otherwise, let
$\gamma\in P_y$, and suppose that $\gamma$ has smallest depth possible
for the elements in $P_y$, and that $\gamma$ has orbital $A = (a,b)$.
$y$ has a fundamental domain $D_y = [a,ay)$, and each element of $P_y$
may be conjugated by some power of $y$ so that the resultant element's
orbital lies in the fundamental domain $D_y$ (if some element,
$\beta$, conjugates to contain $a$ in its orbital, then either that
conjugate has that its orbital fully contains the orbital $A$, which
is impossible by our choice of $\gamma$ as having a minimal depth
orbital of the orbitals of all the elements in $P_y$, or the signed
orbitals of $\gamma$ and of the conjugate of $\beta$ form a transition
chain of length two, which is impossible since $G$ is solvable).  We
can now replace $P_y$ by the conjugates of the original $P_y$ found
above, and the group generated by the new $P_y$ with $y$ will be
identical to the group generated by the old $P_y$ with $y$.  However,
now that all of the elements of $P_y$ have supports in the same
fundamental domain of $y$, we have that $H_y = \langle P_y, y\rangle$
is isomorphic to $K_y\wr\Z$, where $K_y = \langle P_y \rangle$.  But
$K_y$ is a solvable group of precisely the type mentioned in the
hypothesies of the lemma, with derived length $k$ less than $n$, so
that $K_y$ is isomorphic to a group in $\left\{1\right\}\pc$ with
derived length $k<n$, which implies that $H_y$ is isomorphic to a
group in $\left\{1\right\}\pc$ (being the result of a group in
$\left\{1\right\}\pc$ being wreathed with a $\Z$ factor on the right)
with derived length $k+1\leq n$.  But this argument holds for every
$y$ in $Y$ so that
\[
G\cong \bigoplus_{y\in Y}H_y
\]
and since all of the groups in this countable direct sum have derived
length less than or equal to $n$ (and at least one of them has derived
length $n$), we see that $G$ is isomorphic to a group in
$\left\{1\right\}\pc$ with derived length $n$.  \qquad$\diamond$

The following is commonly used without comment in the remainder, since
we often work in the situation where a group in $\ploi$ is balanced
and admits no transition chains of length two.

\brk

\label{nestedOrbitals}
Suppose $G\leq\ploi$ and $G$ is balanced and does not admit transition
chains of length two.  If $(A,g)$ and $(B,h)$ are signed
orbitals of $G$ and $A\cap B\neq \emptyset$ then either
$\overline{A}\subset B$, $\overline{B}\subset A$, or $A=B$.
\erk
pf:

Suppose $(A,g)$ and $(B,h)$ are signed orbitals of $G$ and $A\cap
B\neq\emptyset$.  Suppose $A=(a,b)$ and $B=(c,d)$.  Let $H =
\left<g,h\right>$, and let $C = (e,f)$ be the orbital of $H$ which
contains $A\cup B$.  We will assume without loss of generality that
$a\leq c$.  Since $a$ is not interior to any orbital of $h$, we must
have that $e = a$.  If $b< d$ then in order for
$\left\{(A,g),(B,h)\right\}$ not to form a transition chain of length
two, we must have that $c = a = e$.  In this case $f = d$, since $d$
cannot be in any orbital of $g$ without creating a different
transition chain of length two.  Now $h$ will realize the orbital $C$
of $H$ consistently, and $g$ realizes the left end of $C$ in $C$, so
$g$ must then also realize the right end of $C$.  This means that $g$
realizes both ends of $C$ in $C$ and must be consistent on the ends of
$C$ since $h$ is consistent.  By Lemma \ref{cFullSupport} we see
that $g$ cannot have any fixed points in $C$.  But then $b=d$, which
contradicts our assumption that $b<d$.  In particular, we must have
that $d\leq b$.  If $d = b$ then by reasoning similar to the previous
case, we must have that $a = c$, so that $A = B$.  In particular let
us assume that $d<b$.  Again, if $a=c$ we will have that $A$ is an
imbalanced orbital of $H$, so we must have that $c>a$.  In particular,
we have shown that if $a\leq c$ then either $\overline{B}\subset A$ or
$A = B$.  By a symmetric argument, if $c\leq a$ then either
$\overline{A}\subset B$ or $A = B$.  In particular, either $A = B$,
$\overline{A}\subset B$, or $\overline{B}\subset A$.\qquad$\diamond$

One great tool for technical analysis of a subgroup $G$ of $\ploi$ is
the \emph{split group of $G$}\index{group!split}.  It is motivated by
the hypothesies of Lemma \ref{oneBumpGenerators}.  Suppose $G$ is a
subgroup of $\ploi$, and let $\Gamma$ be the maximal set of elements
of $\ploi$ which all have single orbitals, where if $\gamma\in\Gamma$,
then $\gamma$ is identical to an element $h$ of $G$ on $\gamma$'s
orbital.  The group $\langle\Gamma\rangle$ is the split group of $G$.
Note that $G\leq \langle \Gamma\rangle$.

Here is another technical lemma, which we use for our main result below.

\btl
\label{splitStable}

Suppose $G$ is a balanced subgroup of $\ploi$ that admits no
transition chains of length two, that has derived length $t$ for
some $t\in\N$, and that has split group $H$.  If $(A,h)$ is
a signed orbital of $H$, then there is $g\in G$ so that $(A,g)$ is a
signed orbital of $G$.

\etl

pf: 

Let us denote by $\Gamma$ the maximal set of elements of $\ploi$ which
all have one orbital, so that if $\gamma\in\Gamma$ then $\gamma$ is
identical to an element of $G$ over $\gamma$'s orbital, so that $H =
\langle\Gamma\rangle$.

For each non-negative $k\in\Z$, let $P(k)$ be the statement that if
$\gamma = \gamma_1\gamma_2\cdots\gamma_k$ is a product of elements of
$\Gamma$ of length $k$ for some non-negative $k\in\Z$, and if $A$ is
an orbital of $\gamma$, then there is an element $\alpha$ of $\Gamma$
which also has orbital $A$.

We will prove that $P(k)$ is true for all non-negative $k\in\Z$ via
induction, at which point we will have proven the lemma.

If $k = 0$ the statement $P(0)$ is vacuously true.  If $k = 1$ the
statement $P(1)$ is trivially true.  The the key statement of the
induction proof occurs when $k = 2$.  Suppose $\gamma =
\gamma_1\gamma_2$, where $\gamma_1$, $\gamma_2\in\Gamma$, and let $A$
be an orbital of $\gamma$.  We have two possible cases, which will be
exhaustive as a consequence of Remark \ref{nestedOrbitals}.

\be

\item $\gamma_1$ and $\gamma_2$ share an orbital $C$ with
$\bar{A}\subset C$, or 
\item at least one of $\gamma_1$ and $\gamma_2$ has the orbital $A$,
while the other has orbital $B$ with $B\subset A$ or $B\cap
A=\emptyset$.

\ee

Let us suppose it is the first case.  Then there are elements $f$ and
$g$ in $G$ with $f$ identical to $\gamma_1$ on $C$, and $g$ identical
to $\gamma_2$ on $C$.  The product $fg$ is an element in $G$ with
orbital $A$ so that $fg$ is identical to $\gamma$ on $A$, so in this case,
$P(2)$ is true.

The second case trivially satisfies the statement $P(2)$, so we have
shown $P(2)$ is true in all cases.

Therefore let us assume that $k>2$ and that we know
that $P(m)$ is true for all integers $m$ with $0\leq m<k$.

Let $\gamma = \gamma_1\gamma_2\cdots\gamma_m$ be a product of $m$
elements from $\Gamma$ with an orbital $A$.  Let $g_{m-1} =
\gamma_1\gamma_2\cdots\gamma_{m-1}$.  The product $g_{m-1}\gamma_m$
has orbital $A$, so either both $g_{m-1}$ and $\gamma$ share an
orbital $C$ with $A\subset C$, or at least one of $g_{m-1}$ and
$\gamma_{m}$ has orbital $A$.  If either of these two elements has
orbital $A$, we are done, since we have found a shorter product of
elements of $\Gamma$ that produces the orbital $A$ (possibly length
one).  So let us assume that instead that $\bar{A}\subset C$, and $C$
is the orbital of both $g_{m-1}$ and $\gamma_m$.  In this case, there
is a smallest index $j$ with $1\leq j \leq m$ so that $\gamma_j$ has
orbital $D$ with $A\subset D$ and all the orbitals of the $\gamma_i$
are contained in or disjoint from $D$.  We note that we can permute
the order of the products, moving $\gamma_j$ all the way to the left
and replacing the $\gamma_i$ that $\gamma_j$ moves past with a
conjugate element, which we can still find in $\Gamma$ (since the
orbital of $\gamma_j$ always contains or is disjoint from the orbitals
of the elements we are passing across, the new conjugate elements are
guaranteed to be in $\Gamma$, which is not the case if the orbital of
$\gamma_j$ was interior to the orbital we were trying to pass across),
so that the total number of elements in the product stays constant,
and the resulting product still is the same element of $H$.  Repeat
this procedure with the smallest index in the set of indices
$\left\{2,3,\ldots,m\right\}$ so that we place this new element in the
second slot of the product that produces $\gamma$.  $\gamma_1$ and
$\gamma_2$ now have the largest orbitals of the orbitals of the
$\gamma_i$ which contain $A$.  It is immediate therefore that $D$ is
also the orbital of $\gamma_2$, else $\gamma$ will have orbital $D$
which contains $\bar{A}$.  There are elements $g_1$ and $g_2$ in $G$
which agree with $\gamma_1$ and $\gamma_2$ over $D$, so that the
element $g1g2$ in $G$ has every orbital of $\gamma_1\gamma_2$ as an
orbital of $g_1g_2$, so that these orbitals are all realized by
elements in $\Gamma$.  Now, if $\gamma_1\gamma_2$ has orbital $D$,
then by our induction hypothesis, we are finished, since $g_1g_2$
would be an element of $G$ which behaves as the product
$\gamma_1\gamma_2$ over $D$, so that we could write $\gamma$ as a
shorter product of elements in $\Gamma$.  Therefore let us assume that
the product $\gamma_1\gamma_2$ does not have $D$ as an orbital.  In
this case all of the orbitals of $\gamma_1\gamma_2$ are interior
orbitals of $D$.  Replace the product $\gamma_1\gamma_2$ by the
product $\tau_1\tau_2\cdots\tau_r$ where the $\tau_i$ are elements in
$\Gamma$ that agree with $\gamma_1\gamma_2$ on the $r$ element
orbitals of $\gamma_1\gamma_2$.  Since these orbitals are disjoint,
the ordering in this product is immaterial.  We can even remove any
$\tau_i$ with an orbital that does not contain the boundary of $A$,
since such a $\tau_i$ will not effect whether $A$ is a resultant
orbital of the overall product.  In particular, we will retain at most
two $\tau_i$, so that the total length of the product did not change.
However, there are fewer elements in the product with orbital $D$.  In
particular, we can now repeat this process.  At each stage, the number
of elements in the product that realize the largest orbital containing
$A$, or the size of that largest orbital, is decreasing, while the
number of elements in the product does not go up.  Since $G$ is
solvable, this process must eventually halt, since the largest element
orbital of the new $\gamma_i$ cannot keep shrinking using orbitals in
$\Gamma$, or $G$ would possess an infinite tower.  \qquad$\diamond$

\bc 
\label{splitLength}

Suppose $G$ is a subgroup of $\ploi$ and let $H$ be the split
group of $G$, then the derived length of $G$ equals the derived length
of $H$.  
\ec 

pf: 

If $G$ is solvable, then $G$ must be balanced and admits no
transition chains of length two, so the Lemma \ref{splitStable}
applies, and therefore given any tower $T$ of $H$ we can find a tower
with the same orbitals in $G$.  The other direction is immediate.

If $G$ is non-solvable, then both $G$ and $H$ have towers of arbitrary
height.  
\qquad$\diamond$

The following lemma, and its corollary, complete our proof of both
Lemma \ref{shallowGroupEmbeddings} and Theorem
\ref{solveClassification}.  Note that the corollary is simply a
restatement of Lemma \ref{solveInM}.

\bl 
\label{HInOneP}

If $G$ is a solvable subgroup of $\ploi$ with derived length $n$, then
$G$ is isomorphic to a subgroup of a group $H$ in
$\left\{1\right\}\pc$ so that $H$ has derived length $n$.

\el 

pf: 

Suppose that $n\in\Z$ with $n\geq 0$ and that $G$ is a solvable
subgroup of $\ploi$ with derived length $n$.  We see immediately
that $G$ is balanced and admits no transition chains of length $2$.

Let $H$ be the split group of $G$, and let $\Gamma$ be the collection
of one-bump generators of $H$, that is, $\Gamma$ is the largest
collection of one-bump elements in $\ploi$ which satisfies that if
$g_A\in\Gamma$ with orbital $A$, then there is a $g\in G$ with orbital
$A$, and the element $g_A$ equals $g$ on $A$.  We know from lemma
\ref{splitLength} that $H$ is also solvable with derived length $n$.

Let $X_1$ represent the set of signed orbitals of $H$ with depth $1$.
For each orbital $A$ in $O_{X_1}$, there is a non-empty set of
controllers of $H$ for $A$ in the set $\Gamma$.  Let
$\phi_1:O_{X_1}\to \Gamma$ represent a function that associates to
each orbital in $O_{X_1}$ a controller of $H$ for that orbital which
moves points to the right on the orbital.  Let $Y_1 = O_{X_1}\phi_1$
be the image of $\phi_1$.  We note that each pair of elements in $Y_1$
have disjoint support, and trivially, that no element of
$Y_1$ can be conjugated by an element of $H$ to share an orbital with
a different element of $Y_1$.  Now, by the definition of controller,
$Y_1$ consists of a set of generators in $\Gamma$ sufficient so that
the set $\Gamma_1 = Y_1 \cup (\Gamma \backslash S_{X_1})$ generates
$H$.

For each $y\in Y_1$, let $A_1^y$ represent the orbital of $y$.  We may
partition the elements of $\Gamma_1$ into sets $P^y_1$ indexed by the
set $Y_1$ so that $\gamma\in P^y_1$ if the orbital of $\gamma$ is
contained in $A_1^y$.  Now let $X^y_2$ represent the set of signed
orbitals of elements in $P^y_1$ with depth two in $H$.  If $X^y_2$ is not
empty, let $(A^y_2,\gamma^y_2)$ be an element in $X^y_2$, and let
$a_y$ be the left end of the orbital $A^y_2$.  Since the orbitals of
elements of $H$ are always orbitals of $G$ by Lemma \ref{splitStable},
we see that $H$ admits no transition chains of length two.  In
particular, we can use $y$ to conjugate every signed orbital of
$X^y_2$ into the fundamental domain $[a_y,a_yy)$ to produce the set
$D^y_2$, since all of these orbitals can be conjugated to fit in the
fundamental domain (else we can get an orbital to have $a_y$ in its
interior, creating a transition chain of length two).  And likewise we
can conjugate every signed orbital in $P^y_1$ of depth greater than
two into the fundamental domain as well, producing the set
$T^y_2$. The collections of conjugates in the fundamental domain have
nice properties:

\be
\item If two elements of $D^y_2$ have orbitals that non-trivially
intersect each other, then they actually have identical orbitals.
\item $H_{A^y_1}=\langle y,S_{D^y_2},S_{T^y_2}\rangle$.  \ee 

Let $\phi^y_2:O_{D^y_2}\to S_{D^y_2}$ be a function that picks for
each orbital of depth two in $O_{D^y_2}$ a controller that moves
points to the right for that orbital as before.  Let $Y_2$ be the
union of all the images of the functions $\phi^y_2$ across the $Y$
index set, so that we have now picked an controller that moves points
to the right on its orbital for every conjugacy class of depth two
orbitals of $H$ (note that a conjugate of a controller is also a
controller), so that if $\Gamma_2 = (\Gamma_1\backslash (\cup_{y\in
Y}S_{X^y_2}))\cup Y_2$ then $H = \langle \Gamma_2\rangle$.

In a like fashion we can inductively proceed to pick sets of
controllers, one for each conjugacy class of element orbital of depth
$i$, where $i$ is an index less than or equal to $n$, in exactly the
same fashion as discussed for forming the set $Y_2$ above.  This
process will steadily improve the sets of generators $\Gamma_i$, so
that finally $\Gamma_n$ will be a set of generators for $H$ where each
generator in $\Gamma_n$ has exactly one orbital, and where no
generator can be conjugated in $G$ to share an orbital with another
generator in $\Gamma_n$.  In particular, by Lemma
\ref{oneBumpGenerators}, $H$ is isomorphic to a group in
$\left\{1\right\}\pc$ with derived length $n$, and therefore $G$ is
isomorphic to a subgroup of a group in the class $\left\{1\right\}\pc$
with derived length $n$.  \qquad$\diamond$

\bc
\label{solveEmbeddings}
If $G$ is solvable in $\ploi$ of derived length $n$, then $G$ embeds
in $G_n$.  
\ec 

pf: 

Suppose $G$ is a subgroup of $\ploi$ with derived length $n$.  Lemma
\ref{HInOneP} guarantees that $G$ embeds in a group
$H\in\left\{1\right\}\pc$ with derived length $n$.  But now, Lemma
\ref{BWinM} guarantees that $H$ embeds in $G_n$, so that $G$ must
embed in $G_n$.  \qquad$\diamond$

\newpage

\vbox to 1.5truein{}

\section{Algebraic classification of non-solvable subgroups of $\ploi$}
\newtheorem{construction}{Construction}
\newtheorem{process}{Process}

Our primary result is Theorem \ref{nonSolveClassification}, which
states that $W$ is a subgroup of any non-solvable subgroup of $\ploi$.
This result is easier than what we actually prove in this section.
See the next subsection for a list of the main results.

It is immediate from Theorem \ref{geoClassification} that a subgroup
$H$ of $\ploi$ is non-solvable if and only if $H$ admits towers of
arbitrary height.  If $H$ admits towers of arbitrary height then we
know from section \ref{solveClassificationSection} that $H$ admits
exemplary towers of arbitrary height.  If we consider the realization
of the groups $W_i$ in the introduction, then we might suspect that
the signatures of an exemplary tower of height $n$ will generate a
group isomorphic to $W_n$, and that the signatures of a countably
infinite exemplary tower will generate one of the three groups
$(\wr\Z\wr)^{\infty}$, $(\wr\Z)^{\infty}$, and $(\Z\wr)^{\infty}$.
This is false, as individual signatures of an exemplary tower may have
multiple orbitals; we are only guaranteed that the ends of the
orbitals of the exemplary tower are arranged nicely with respect to
the set of all orbitals of the signatures of the tower.  Nonetheless,
this suspicion has still some kernal of truth in it, and so it allows
us a toe-hold on what work needs to be done (cleaning up the non-tower
orbitals) in order to find towers so that their signatures generate
groups that we can analyze.

Suppose that $H$ is a non-solvable subgroup of $\ploi$, then we can
informally outline our approach in this section in terms of $H$ as follows.
\be

\item Show that given an exemplary infinite tower for $H$, we can pass
to an even better infinite tower so that the group generated by the
signatures of the new tower is one we can recognize.  This takes two
generic steps.  

\be

\item Clean the tower further so that the orbitals of any signature in
the tower which are contained in an orbital of the tower are arranged
nicely in that orbital.

\item Clean the tower further still so that the orbitals of the
signatures which are contained away from the orbital of $H$ which
supports the tower are arranged nicely with respect to each other.

\ee

\item Show that even in circumstances where we only know that we can
find exemplary towers of arbitrary height in $H$, with reasonable
extra information we can conclude there are actually exemplary
infinite towers in $H$.

\item Show that $W$ embeds in any of the groups generated by the process
outlined above when $H$ admits an infinite exemplary tower.
\item Show that if $H$ does not admit infinite towers, but is
non-solvable, then $W$ embeds in $H$ as a subgroup.

\ee

In order to carry out the two types (a) and (b) of tower cleaning
mentioned above, we will need to carry out a technical analysis
describing the results of building two types of commutators, so recall that in these investigations we use the commutator symbol as follows:
\[
[a,b] = a^{-1}b^{-1}ab = a^{-1}a^b.
\]

The main points of the outline above will be visible in the statements
of results below, although one should not be mislead into believing that the
statements below precisely follow the development of the proof of the
primary result.

\subsection{Statement of non-solvablity results}

First we find conditions under which we can guarantee that we can find
a subgroup isomorphic to one of the wreath products mentioned in the
introduction.

\bt
\label{tallWreath}

Any subgroup $H$ of $\ploi$ which admits a tall tower contains a
subgroup isomorphic to $(\Z\wr)^{\infty}$.

\et

\bt
\label{deepWreath}

Any subgroup $H$ of $\ploi$ which admits a deep tower contains a
subgroup isomorphic to $(\wr\Z)^{\infty}$.

\et

\bc

If $H$ is a subgroup of $\ploi$ and $H$ admits transition chains of
length two then $H$ has subgroups isomorphic to both
$(\wr\Z)^{\infty}$ and $(\Z\wr)^{\infty}$.

\ec

We also can find one of these wreath products in any finitely
generated non-solvable subgroup of $\ploi$.

\bt
\label{fgNonSolve}

Any finitely generated subgroup of $\ploi$ with
towers of arbitrary height contains one of $(\Z\wr)^{\infty}$ or
$(\wr\Z)^{\infty}$ as a subgroup.

\et

However, there are non-solvable subgroups of $\ploi$ which do not
contain any of the three main infinite wreath products mentioned in
the introduction.

\bl 
\label{stuffInW}

Neither $(\wr\Z)^{\infty}$, nor $(\Z\wr)^{\infty}$, nor
$(\wr\Z\wr)^{\infty}$ embed in $W$.  

\el

Although converse statements are false, as stated below.

\bl
\label{WInStuff}
$W$ embeds in both $(\wr\Z)^{\infty}$ and $(\Z\wr)^{\infty}$.

\el

Therefore, the following lemma completes a proof of Theorem
\ref{nonSolveClassification}.

\bl
\label{arbTowersW}

If $H$ is a non-solvable subgroup of $\ploi$ which does not admit
infinite towers then $H$ contains a subgroup isomorphic to $W$.

\el

and from this we have an immediate consequence, based on the easy fact
that $W$ admits no finite index solvable subgroups.

\bc
Virtually solvable subgroups of $\ploi$ are solvable.
\ec

\subsection{Improving exemplary towers}

In this section we build technical results guaranteeing that we can
often find towers whose signatures generate groups that we can
analyze.

The first result is a simple refinement of Lemma
\ref{transferExistence} using Lemma \ref{transitiveOrbital}.

\btl
\label{spanningConjugate} 

Suppose $H$ is a balanced subgroup of $\ploi$ and $H = \langle\alpha,
\beta\rangle$ for some two elements $\alpha$, $\beta\in\ploi$.
Suppose further that $A$ is an inconsistent orbital of $H$ and
$\alpha$ realizes both ends of $A$ while $\beta$ realizes neither.
There is a conjugate $\gamma$ of $\beta$ in $H$ which has an orbital
$B\subset A$ so that the fixed set of $\alpha$ in $A$ is contained in
$B$.

\etl

Pf: Let $F_{\alpha}$ represent the fixed set of $\alpha$ in $A$, and
let $x = \inf(F_{\alpha})$ and $y = \sup(F_{\alpha})$.  By Lemma
\ref{transitiveOrbital}, since $A$ is an orbital of $H$, there is
$\theta\in H$ so that $x\theta >y$.  By the continuity of $\theta$,
there is $x_1<x$ so that $x_1>y$ as well.  Let $z = x_1\theta^{-1}$ so
we have
\[
z<z\theta = x_1<x<y<x_1\theta<x\theta
\]
Now since $F_{\alpha}$ in the orbitals of $\beta$, we see that $\beta$
has an orbital $C = (r,s)$ so that $r<x<s$.  There is a power $k\in\Z$
so that $r\alpha^k=q<z$.  Now, $\beta_1 = \beta^{\alpha^k}$ has
orbital $D=(q,t)$ induced from $C$ by the action of $\alpha^k$, and
$D$ satisfies $q<z<x<t$.  Set $\gamma = \beta_1^{\theta}$.  $\gamma$
has orbital $B =(u,v)$ induced from $D$ by the action of $\theta$ on
$\beta_1$, and $u = q\theta<x<y<t\theta = v$.

\qquad$\diamond$

The following lemma is more involved, and plays a key role in
the proof of the lemma following immediately after.  

\btl
\label{chainSplitting}

Suppose $H$ is a balanced subgroup of $\ploi$ and $H = \langle
\alpha,\beta\rangle$ for some two elements $\alpha$, $\beta\in\ploi$.  If
$H$ has an inconsistent orbital $A$, and $\beta$ realizes the ends of no
orbitals of $H$, then there are elements $\alpha_1$ and $\beta_1$ in
$H$ so that if $H_1 = \langle\alpha_1,\beta_1\rangle$ then $H_1$,
$\alpha_1$, and $\beta_1$ satisfy the following properties:

\be

\item $A$ is an inconsistent orbital of $H_1$.
\item $\beta_1$ realizes no ends of any orbital of $H_1$.
\item Every inconsistent orbital of $H_1$ is the union of the orbitals
of a transition chain of length three whose first and last orbitals
are orbitals of $\alpha_1$, and whose second orbital is an orbital of
$\beta_1$.
\item $\alpha_1$ moves points to the left on its leading orbital in each of the
inconsistent orbitals of $H_1$.

\ee
\etl

pf: 

Set $\alpha_1$ to be either $\alpha$ or $\alpha^{-1}$, so that
$\alpha_1$ moves points to the left on its leading orbital $B$
contained in $A$.

Suppose $n\in\N$ and $H$ has $n$ inconsistent orbitals.  There are
non-negative integers $r$, $s\in\Z$ so that $n = r+s$, with $r>0$,
where $r$ is the number of the inconsistent orbitals of $H$ that have
$\alpha_1$ moving points to the left on its leading interior orbitals
in these inconsistent orbitals.  Let $\mathscr{B}=\left\{B_i\,|\,1\leq
i\leq r\right\}$ represent the collection of inconsistent orbitals of
$H$ where $\alpha_1$ moves points to the left on its leading orbital
in each of these orbitals, indexed from left to right.  Let
$\mathscr{C}=\left\{C_j\,|\,1\leq j\leq s\right\}$ represent the other
inconsistent orbitals of $H$, and note that if $s = 0$, then this
could be an empty collection.

By Technical Lemma \ref{spanningConjugate}, for each orbital $B_i$ in
$\mathscr{B}$ there is an element $\gamma_i$ in $H$, which is a
conjugate of $\beta$, so that the fixed set of $\alpha_1$ in $B_i$ is
contained in a single orbital of $\gamma_i$.  Likewise, for each
orbital $C_j$ in $\mathscr{C}$ there is an element $\theta_j$ in $H$,
which is a conjugate of $\beta$, so that the fixed set of $\alpha_1$
in $C_j$ is contained in a single orbital of $\theta_j$.

Firstly, inductively replace each element $\gamma_i$, for $i>1$, by a
conjugate of $\gamma_i$ by a high negative power of $\alpha_1$ so that
whenever $k$, $i\in\N$, with $1\leq i<k\leq r$, we have that the
closure of the union of all of the orbitals of $\gamma_k$ that
intersect the orbitals of $\gamma_i$ nontrivially in $B_i$ is
actually fully contained in the single orbital of $\gamma_i$ that
contains the fixed set of $\alpha_1$ in $B_i$.  Note that conjugating
by high negative powers of $\alpha_1$, we are pushing the orbitals of
the conjugates closer to the exterior boundaries of the fixed set of
$\alpha_1$ in each $B_i$, so that the last sentence is possible.

For each $i\in\left\{1,2,\ldots,r\right\}$, $\gamma_i$ has an orbital
$D_i$ that contains the fixed set of $\alpha_1$ in $B_i$, as well as the
closure of all of the orbitals of $\gamma_k$ in $B_i$ for all $i<k\leq
r$.  

We will inductively define a sequence of elements $(\rho_i)_{i =
1}^r$.  Then, modulo replacing some of the $\gamma_i$ by their
conjugates by more negative powers of $\alpha_i$, the $\rho_i$ will
have the following properties: 

\be

\item $\rho_1 =\gamma_1$.
\item For all indices $i>1$, $\rho_i$ will be either a conjugate
of $\rho_{i-1}$ by some power of $\gamma_i$, or
$\rho_i=\rho_{i-1}\gamma_i$.
\item For all indices $i$, $\rho_i$ will have an orbital $E_i$ in
$B_i$ that fully contains the fixed set of $\alpha_1$ in $B_i$.
\item If $i<r$, the orbital $E_i$ of $\rho_i$ will contain the closures of the
orbitals of $\gamma_j$ in $B_i$ for all integers $j$ with $i<j\leq r$.
\item If $i>1$, for each integer $j$ with $1\leq j <i$, $\rho_i$ will
have $E_j$ as one of its orbitals.

\ee 

Firstly, set $\rho_1 =\gamma_1$, and $E_1 = D_1$.  By construction,
$\rho_1$ satisfies the five inductive properties.  If $r = 1$, we are
done.  If not, suppose that $k$ is an integer so that $1<k\leq r$ and
for all $i\in\left\{1,2,\ldots,k-1\right\}$ we have that $\rho_i$ is
defined and satisfies the five defining properties of the induction.
Our analysis now breaks into two cases.

If $\rho_{k-1}$ has an orbital $F_k$ containing either end of $D_k$,
then there is some integer $j$ so that $\rho_k =
\rho_{k-1}^{\gamma_k^j}$ will have orbital $E_k$ induced from $F_k$ by
the action of $\gamma_k^j$ so that $E_k$ will contain the fixed set of
$\alpha_1$ in $B_k$, as well as the closure of all of the orbitals of
$\gamma_j$ in $B_k$ for integers $j$ with $i<j\leq r$ (if $k<r$) .  

If $\rho_{k-1}$ does not have an orbital $F_k$ containing either end
of $D_k$, then we have to handle the case where $\rho_{k-1}$ has
orbitals in $D_k$ that share ends with $D_k$ separately before
continuing.

If $\rho_{k-1}$ has orbitals in $D_k$ that share ends with $D_k$ then
replace $\gamma_k$ and all later $\gamma_j$ with conjugates of these
elements by a high negative power of $\alpha_1$ so that $D_k$ either
has an end contained in an orbital of $\rho_{k-1}$, or shares no end
with an orbital of $\rho_{k-1}$, and repeat the whole inductive
definition of $\rho_k$.

If $\rho_k$ is still undefined, then set $\rho_k =
\rho_{k-1}\gamma_k$.  Note that since $\rho_{k-1}$ has no orbitals in
$D_k$ that share ends with $D_k$, the product $\rho_k =
\rho_{k-1}\gamma_k$ realizes both ends of $D_k$ consistently, and
therefore realizes $D_k$ consistently since $H$ is balanced.
Therefore define $E_k = D_k$ and note that $\rho_k$ actually has $E_k$
as an orbital.

At this point, $\rho_k$ and $E_k$ are both defined, and we can
continue with our main argument.  Note that $E_k$ contains the closure
of all of the orbitals of all of the $\gamma_i$ for $i>k$, and that
for each integer $j$ in $1\leq j\leq {k-1}$, the closure of the
orbitals of $\gamma_k$ in $B_j$ are fully contained in the orbital
$E_j$, so that $\rho_k$ will have $E_j$ as an orbital as well.  Now by
construction, $\rho_k$ satisfies the five defining properties of the
induction.

We now examine the element $\rho_r$.  Observe that the element
$\rho_r$ contains an orbital $E_k$ in each $B_k$ where the fixed set
of $\alpha_1$ in $B_k$ is fully contained in $E_k$.  $\rho_r$ is
constructed as a sequence of products using various $\gamma_i$'s and
conjugates of $\gamma_i$'s so $\rho_r$ realizes no end of any orbital
of $H$, but is an element of $H$.

In an entirely analogous fashion, if $s>0$, then we can find one
element $\psi_s$ in $H$ which realizes no end of any orbital of $H$ and
which contains an orbital $F_i$ in each $C_i\in\mathscr{C}$ which contains the fixed
set of $\alpha_1$ in that $C_i$.

There is a positive integer $p$ so that $\rho= \rho_r^{\alpha_1^p}$
has the properties that follow:

\be

\item For each integer $i\in\left\{1,2,\ldots,r\right\}$, the closure
of the orbitals of $\psi_s$ in $B_i$ is actually contained in the
orbital $G_i$ of $\rho$ induced from $E_i$ by the action of
$\alpha_1^p$.
\item For each integer $i\in\left\{1,2,\ldots,s\right\}$, the closure
of the orbitals of $\rho$ in $C_i$ is actually contained in the
orbital $F_i$ of $\psi_s$.

\ee

This follows since for each orbital $B_i$ of $\mathscr{B}$, the lead
orbital of $\alpha_1$ in $B_i$ has the property that $\alpha_1$ is
moves points to the left there (and therefore moves points to the
right on the trailing orbital of $\alpha_1$ in $B_i$), and for each
orbital $C_i$ in $\mathscr{C}$, the lead orbital of $\alpha_1$ in
$C_i$ has the property that $\alpha_1$ is moving points to the right
there (and therefore $\alpha_1$ moves points to the left on its
trailing orbital in $C_i$).

We note in passing that the orbitals $G_i$ of $\rho$ contain the
orbitals $E_i$ of $\rho_r$, and therefore the fixed set of $\alpha_1$ in
the $B_i$.

Now there is a power $q$ of $\psi_s$ so that the element $\beta_1 =
\rho^{\psi_s^q}$ will have the following nice properties: 

\be

\item For each integer $i\in\left\{1,2,\ldots,s\right\}$, the orbitals
of $\beta_1$ in $C_i$ have trivial intersection with the fixed set of
$\alpha_1$ in $C_i$.

\item For each integer $i\in\left\{1,2\ldots,r\right\}$, $\beta_1$ will
have the orbital $G_i$ which contains the fixed set of $\alpha_1$ in
$B_i$.

\ee 

The first property follows since the orbitals of $\rho$ in the $C_i$
are contained in the orbitals $F_i$ of $\psi_s$, and so the
conjugation of $\rho$ by a high power of $\psi_s$ will throw these
orbitals off of the fixed set of $\alpha_1$ in the $C_i$.  The second
property follows since the orbitals of $\psi_s$ are fully contained in
the orbitals $G_i$ of $\rho$ in the $B_i$, so that conjugation of
$\rho$ by $\psi_s$ to any power will not change these orbitals.

We observe that by generating the group $H_1 = \langle
\alpha_1,\beta_1\rangle$, we succeed in producing the promised group,
since $A = B_k$ of some integer $k\in\left\{1,2,\ldots,r\right\}$ is
an inconsistent orbital of $H_1$, $\beta_1$ realizes no ends of any
orbital of $H_1$, and on each inconsistent orbital $D$ of $H_1$, $D$
is the union of a transition chain of length three of the form
$\left\{(L,\alpha_1),\,(G_k,\beta_1),\,(R,\alpha_1)\right\}$ where $L$
is an orbital of $\alpha_1$ where $\alpha_1$ moves points to the left,
and $L$ contains the left end of $G_k$, and $R$ is an orbital of
$\alpha_1$ where $\alpha_1$ moves points to the right, and $R$
contains the right end of $G_k$, and where $G_k$ is one of the
orbitals $G_i$ of the definition of $\beta_1$.

\qquad$\diamond$ 

We already know that imbalanced groups contain copies of Thompson's
group $F$, which itself contains many copies of $(\wr\Z\wr)^{\infty}$,
(see \cite{BrinEG} for details of this last point), so in our
investigations, we will generally work under the assumption that the
groups we are examining are balanced.

The next result is an example of what we are aiming for when we
``improve'' exemplary towers, and depends on the previous technical
lemma.  We will use this result in the proofs of Lemma's
\ref{tallZWreath} and \ref{deepWreathZ} to immediately reduce to the
special case of examining exemplary towers whose signatures generate a
group which admits no transition chains of length two.

\bl
\label{messyWreath}

If $G$ is balanced and $G$ contains transition chains of length two,
then there is an exemplary bi-infinite tower
$E=\left\{(B_i,\beta_i)\,|\, i\in\Z\right\}$ for $G$ where the group
$H = \langle S_E\rangle$ admits no transition chains of length two.

\el

pf:

Since $G$ admits transition chains of length two, we can
find two elements $\alpha$ and $\beta$ which form the set of
signatures of a transition chain of length two.  Let $K = \langle \alpha,\beta\rangle$.

The orbitals of $K$ are the components of the union of the orbitals of
$\alpha$ and $\beta$.  Some of these orbitals may be consistent
orbitals for $K$, so that at least one of $\alpha$ or $\beta$ realize
these orbitals.  The other orbitals are inconsistent, and are formed
by the union of a subcollection of orbitals of $\alpha$ and orbitals of
$\beta$.  Each such subcollection admits a transition chain of length
two.  A chief feature of the inconsistent orbitals is that
one of $\alpha$ or $\beta$ must realize both ends of any particular
such orbital, since $K$ is balanced.  Since $K$ admits transition
chains of length two, at least one of the orbitals of $K$ is
inconsistent.

 We are going to analyze the orbitals of $K$ still further.  Any
 particular orbital of $K$ has one of six types, the first three are
 consistent, and the last three are inconsistent:

\be

\item (Type AB) Both $\alpha$ and $\beta$ consistently realize this orbital.
\item (Type Ab) $\alpha$ consistently realizes this orbital, but not $\beta$.
\item (Type aB) $\beta$ consistently realizes this orbital, but not $\alpha$.
\item (Type $\underline{a}\underline{b}$) Both $\alpha$ and $\beta$
inconsistently realize both ends of this orbital.
\item (Type $\underline{a}$b) $\alpha$ inconsistently realizes both
ends of this orbital, but $\beta$ realizes
neither end of this orbital.
\item (Type a$\underline{b}$) $\beta$ inconsistently realizes both
ends of this orbital, but $\alpha$ realizes
neither end of this orbital.

\ee

We know that $K$ has at least one orbital, let us call it $A$, of type
$\underline{a}\underline{b}$, $\underline{a}b$, or $a\underline{b}$,
and we will assume without meaningful loss of generality that $A$ has
one of the first two types.  Let $F_a$ represent the union of the
fixed sets of $\alpha$ that are contained in the orbitals of $K$ of
type $\underline{a}\underline{b}$, type $\underline{a}b$, and
$a\underline{b}$.  $F_a$ is non-empty, and is entirely contained in
the orbitals of $\beta$.  In particular, by Remark
\ref{transitiveElementOrbital} there is a $N_1\in\N$ so that for all
$k\in\N$ with $k\geq N_1$ we have $F_a\beta^k\cap F_a=\emptyset$ in
orbitals of type $\underline{a}\underline{b}$ and $\underline{a}b$ (in
orbitals of type $a\underline{b}$, the interior componenents of $F_a$
are moved off of themselves).  Similarly, let $S$ represent the
support of $\alpha$ in the orbitals of $K$ of type $aB$, then there is
$N_2\in\N$ so that for all $k\geq N_2$, we have $S\beta^k\cap
S=\emptyset$.  Let $N = max(N_1,N_2)$.  Considering the other
direction, let $F_b$ represent the fixed set of $\beta$ in the
orbitals of $K$ of type $a\underline{b}$.  Since $F_b$ is contained in
the support of $\alpha$ by definition, there is $M\in\N$ so that for
all $j\geq M$ we have that $F_b\alpha^j\cap F_b=\emptyset$.  Now let
$j\geq M$, and let $k\geq N$, and define $\beta_1=[\alpha^j,\beta^k]$.
We observe that the fixed set of $\alpha$ in the orbitals of $K$ of
type $\underline{a}\underline{b}$ and $\underline{a}b$ is contained in
the orbitals of $\beta_1$.  Note also that the components of $F_a$ in
the orbitals of $K$ of type $a\underline{b}$ which do not realize any
end of an orbital of $K$ are all contained in the support of
$\beta_1$.  The fixed set of $\beta$ contained in the orbitals of $K$
of type $a\underline{b}$ is also contained in the support of
$\beta_1$, since any such point is moved off of $F_b$ by
$\alpha^{-j}$, then moved by $\beta^{-k}$, then moved to someplace
different (from its start) by $\alpha^j$, and finally, $\beta^k$
cannot move the resultant point to its original location in the fixed
set of $\beta$.  Now observe that the orbitals of $\beta_1$ are
either disjoint from $S$, or else are components of $S$ where
$\alpha^j$ behaves as the inverse of $\beta_1$.

We now consider
the group $K_1=\langle \alpha, \beta_1\rangle$, and we will consider
the orbitals of $K_1$ under the same classification as the orbitals of
$K$, where we replace $\beta$ by $\beta_1$ in that classification.

It is immediate to see that $K_1$ still has all the orbitals of $K$ of
type $\underline{a}b$, and that the type of these orbitals is
unchanged.  It is also immediate by construction that the orbitals of
$K$ of type $\underline{a}\underline{b}$ are also orbitals of $K_1$,
although they are now of type $\underline{a}b$.  The orbitals of $K$
of type $Ab$ are also orbitals of $K_1$ of type $Ab$, but the orbitals
of $K$ of type $aB$ are now replaced by a collection of interior
orbitals (all lying properly in the union of the orbitals of $K$ of
type $aB$), each of which is an orbital of type $aB$ that is actually
disjoint from the support of $\alpha$, or else of type $AB$, where
$\alpha^j$ and $\beta_1$ behave as inverses on these orbitals.  The
orbitals of $K$ of type $AB$ are now of type $Ab$, and may have
trivial intersection with the support of $\beta_1$ (if, in fact,
$\alpha$ and $\beta$ commuted on these orbitals).

If $B_1$ is an orbital of $K$ of type $a\underline{b}$, then $B_1$ is
not an orbital of $K_1$.  In this case $K_1$ admits a new collection
of orbitals properly contained in $B_1$.  

We first consider the case where $\beta$ is moving points to the right
on its leading orbital in $B_1$ (and therefore is moving points to the
left on its trailing such orbital). We will suppose $k$ was chosen
large enough so that the closure of the union of the orbitals of
$\beta^{-k}\alpha^j\beta^k$ that are contained in orbitals of $\beta$
in $B_1$ is actually contained in the orbitals of $\alpha$ (and
therefore of $\alpha^j$) which contain components of the fixed set of
$\beta$.  Note that any interior orbital of $\beta$ in $B_1$ is
contained in the union of the orbitals of $\alpha^j$ and
$\beta^{-k}\alpha^j\beta^k$.  Therefore, there are three possible
varieties of resulting orbitals of $K_1$ in $B_1$: firstly, of type
$AB$, where $\beta_1$ actually behaves as $\alpha^{-j}$ on these
orbitals (there may be several of these), secondly, of type $AB$,
where there is only one such orbital, and it contains the fixed set of
$\beta$, or thirdly, of type $\underline{a}\underline{b}$, where there
is one of these if the previous variety did not occur, and it contains the
fixed set of $\beta$ in this case.  We will assume $k$ was chosen
large enough so that these properties of transformation are preserved
over all orbitals of $K$ of type $a\underline{b}$ where $\beta$ is
moving points to the right on its leading relevant orbitals.

In the case of the orbitals of $K$ of type $a\underline{b}$ where
$\beta$ moves points to the left on its leading relevant orbitals. The
results depend heavily on the nature of $\alpha$ in these individual
orbitals.  To clarify the discussion, let us suppose that $B$ is such
an orbital, and discuss the possibilities that arise from the
behaviour of $\alpha$ and $\beta$ on $B$.  

Firstly, let us suppose that $\alpha$ has an orbital that contains the
fixed set of $\beta$ in $B$.  In this case, let us suppose $k$ and $j$
were chosen large enough so that the entire support of $\alpha$ is
contained inside a single fundamental domain of the single orbital of
$\beta^{-k}\alpha^j\beta^k$ that contains the fixed set of $\beta$ in
$B$.  In this case, the group $K_1$ possibly has several orbitals in
$B$, all of type $aB$.  One of these orbitals contains all of the
support of $\alpha$ in $B$, and all of the rest are orbitals of
$\beta_1$ which contain no orbitals of $\alpha$ and are therefore of
type $aB$ with trivial intersection with orbitals of $\alpha$.

Now let us suppose that $\alpha$ has more than one orbital in $B$ that
contains a component of the fixed set of $\beta$.  The first and last
such orbitals of $\alpha$ in $B$ must have that $\alpha$ behaves
inconsistently on these orbitals, otherwise it is easy to create an
imbalanced subgroup of $K_1$.  So now there are two further cases.  

Let us suppose that $\alpha$ moves points to the right on its first
orbital in $B$ which contains a component of the fixed set of $\beta$,
and therefore moves points to the left on the last such.  In this case
$K_1$ has only one orbital in the domain $B$, call it $C$, which is
again of type $a\underline{b}$.  The closure of $C$ is contained in
$B$, and $\beta_1$ moves points to the right on its leading orbital in
$C$ and moves points to the left on its trailing orbital there.  

Now let us suppose $\alpha$ moves points to the left on its leading
orbital in $B$ that contains a component of the fixed set of $\beta$,
and therefore moves points to the right on its trailing orbital in $B$
which contains a component of the fixed set of $\beta$, the group
$K_1$ again has some pure orbitals (type $aB$) plus precisely one
orbital $C$ in $B$, which is again of type $a\underline{b}$, and this
time, $\beta_1$ will move points to the left on its leading orbital in
$C$ and will move points to the right on its trailing orbital in $C$.

The result of all of this analysis is the following, we can choose $j$
and $k$ so that the group $K_1$ has orbitals of the following types:

\be
\item $AB$

Note that in this case $\alpha$ and $\beta_1$ commute on this orbital,
except in the case possibly generated from orbitals of type
$a\underline{b}$ where $b$ moves points right on its leading orbital.

\item $Ab$
\item $aB$

Note here that the behavior of $\alpha$ on this orbital is as the
identity, unless this orbital is contained in an orbital of $K$ of
type $a\underline{b}$, in which case $\alpha$ may have non-trivial
support in this orbital.

\item $\underline{a}\underline{b}$ 

Note that orbitals of this type are always contained in orbitals of
$K$ of type $a\underline{b}$ where $\beta$ moves points to the
right on its first relevant orbital.

\item $\underline{a}b$

Since these are the certain result of an orbital of type
$\underline{a}\underline{b}$ or of an orbital of type $\underline{a}b$
of $K$, we see that $K_1$ will have at least one of these.

\item $a\underline{b}$

These orbitals all have the property that whenever $\beta_1$ moves
points to the left on its leading orbital in these orbitals, then
$\alpha$ moves points to the left on its leading orbital of the
orbitals that contain a component of the fixed set of $\beta_1$.

\ee

In particular, we can repeat this process to create a new element
$\beta_2$ using $\alpha$ and $\beta_1$, and therefore a new group $K_2
= \langle\alpha,\beta_2\rangle$.  $K_2$ improves on $K_1$ since all of
its orbitals of type $a\underline{b}$ have both $\beta_2$ and $\alpha$
moving points to the left on their important leading orbitals.  In
particular, $K_2$ may still have orbitals of type
$\underline{a}\underline{b}$, and of type $AB$ (although here $\alpha$
and $\beta_2$ will commute on these orbitals). $K_2$ may have orbitals
of type $Ab$, but its orbitals of type $aB$ will all have the property
that $\alpha$ is the identity over these orbitals, while $K_2$ will
certainly have orbitals of type $\underline{a}b$.  Repeating the
process one more time to create an element $\beta_3$ and a subgroup
$K_3=\langle \alpha,\beta_3\rangle$ produces a group whose orbitals
are much easier to describe.  $K_3$ will have no orbitals of type $AB$
since $K_2$ had no orbitals of type $aB$ or $a\underline{b}$ that
could produce these orbitals (the types exist, but not with the right
subflavors of $\alpha$ and $\beta_2$ to generate these offspring).
$K_3$ may have orbitals of type $Ab$, but it will have no orbitals of
type $aB$, since the orbitals of type $aB$ in $K_2$ had $\alpha$
behaving as the identity there, and $K_2$ had no orbitals of type
$a\underline{b}$ with $\beta_2$ moving points to the left on its first
sub-orbital $D$ while $\alpha$ was moving points to the right on its
orbital containing the right end of $D$.  $K_3$ will have no orbitals
of type $\underline{a}\underline{b}$, since $K_2$ had no orbitals of
type $a\underline{b}$ with $\beta_2$ moving points to the right on its
first orbital in the orbitals of $K_2$ of this type. $K_3$ will have
at least one orbital of type $\underline{a}b$, and may have several
orbitals of the type $a\underline{b}$, but all of these last will have
$\beta_3$ moving points to the left on its leftmost orbitals in these
orbitals, and $\alpha$ will also move points to the left on its first
orbitals containing the right ends of $\beta_3$'s leftmost orbitals in
these orbitals of type $a\underline{b}$ of $K_3$.

Now, the orbitals of $K_3$ are well understood, and the behaviors of
$\beta_3$ and $\alpha$ on these orbitals are also well understood.  We
now consider the subgroup $K_4$ generated by $\alpha$ and $\beta_4 =
[\alpha^{-j},\beta_3^{k}]$, where $j$ and $k$ are chosen as in the
previous process (note the negative index on $\alpha$).  The point of
this is that now the orbitals of $K_4$ will admit no orbital of type
$a\underline{b}$ with $\beta_4$ moving points to the left on its first
orbital.  Now replacing $K_4$ with $K_5 = \langle
\alpha,\,\beta_5\rangle$ where $\beta_5 = [\alpha^j,\beta^k]$ where
$j$ and $k$ are chosen as before produces a group with no orbitals of
type $a\underline{b}$, repeating one more time to generate $\beta_6$
and $K_6$ in the same fashion that we generated $K_1$ from $K$
produces a group whose orbitals are only of types $Ab$ and
$\underline{a}b$.

Let us consider the orbital $A$ of $K$.  $A$ is also an orbital of
$K_6$, and it is of type $\underline{a}b$. We will now replace $K$ by
$K_6$ and $\beta$ by $\beta_6$ so that $K$ has an orbital of type
$\underline{a}b$ and all of its orbitals are of type $\underline{a}b$
and $Ab$.

Suppose $K$ has $n$ orbitals of type $\underline{a}b$, and let
$\mathscr{O}=\left\{A_i\,|\,1\leq i\leq n, i\in\N\right\}$ represent
this collection, where the indices respect the left to right order of
the orbitals.  By construction we know that $n\geq 1$.  Apply
Technical Lemma \ref{chainSplitting} (above) to replace $\alpha$ and
$\beta$ by new elements, and replace $K$ by the new group generated by
the new $\alpha$ and $\beta$ so that $\beta$ still realizes no end of
any orbital of $K$, and $A_1$ is still an orbital of type
$\underline{a}b$, but where every maximal transition chain (of length
greater than one) formable using $\alpha$ and $\beta$ has length three
(naturally $\alpha$ provides the leading and trailing orbitals for any
such chain), and where $\alpha$ moves points to the left on all of its
leading orbitals in orbitals of type $\underline{a}b$ for $K$ (and
therefore moving points to the right on its trailing such intervals).

Define $\gamma_0 = \beta$. For each $i\in\N$, inductively define
$\gamma_i = \gamma_{i-1}^{\alpha^{k_i}}$ where $k_i$ is chosen large
enough so that in each orbital of $K$ of type $\underline{a}b$, the
closure of all of the orbitals of $\gamma_{i-1}$ is contained in the
single orbital of $\gamma_i$ that contains the fixed set of $\alpha$
in that particular orbital of $K$ of type $\underline{a}b$, and so
that in any particular orbital $X$ of $K$ of type $Ab$, the set of
orbitals of $\gamma_i$ lie between an end of $X$ (already chosen by
the direction $\alpha$ moves points) and the orbitals of
$\gamma_j$ in $X$ for all $j$ with $0\leq j <i$.  Note that this $k_i$
exists, since $\alpha$ moves points to the left on all of its leading
orbitals in orbitals of $K$ of type $\underline{a}b$.  

Similarly, define the $\gamma_i$ for negative integers $i$ by setting
$\gamma_i=\gamma_{i+1}^{\alpha^{k_i}}$ where $k_i$ is chosen large
enough (in the negative direction) so that in any orbital of $K$ of
type $\underline{a}b$ the closure of all of the orbitals of $\gamma_i$
is contained in the orbital of $\gamma_{i+1}$ that contains the fixed
set of $\alpha$ in that orbital of $K$, and so that the orbitals of
$\gamma_i$ in orbitals of $K$ of type $Ab$ are disjoint from all of
the orbitals of the $\gamma_k$ for $k>i$.  Note that this last is
possible in these orbitals of $K$ of type $Ab$, despite the fact that
there are infinitely many such orbitals to avoid, since we are now
conjugating by $\alpha$ to negative powers, so that orbitals of the
$\gamma_i$ for negative $i$ all lie on the other side of the orbitals
of $\gamma_0$ compared to the orbitals of $\gamma_i$ with $i>0$.  

By construction, if $i$ and $j$ are integers, and $i<j$, then any
orbital $C$ of $\gamma_i$ that intersects an orbital $D$ of $\gamma_j$
has the property that $\overline{C}<D$.  Furthermore, $\overline{D}$
is then contained in an orbital of $\gamma_k$ for every integer $k>j$.
In particular, for each integer $i$ there is $m_i\in\N$ so that
whenever $B$ is an orbital of $\gamma_i$, all of the orbitals of
$\gamma_{i-1}$ in $B$ are actually contained in a single fundamental
domain $f_{m_i}$ of $\gamma_i^{m_i}$ in $B$.  Inductively, we see that
all of the orbitals of $\gamma_k$ in $B$ for any integer $k<i$ are
actually also contained in $f_{m_i}$.  Define $\beta_i =
\gamma_i^{m_i}$.  

The group $H = \langle \left\{\beta_i\,|\, i\in\Z\right\}\rangle$ now
admits no transition chains of length two.  For every integer $i$, let
$B_i$ represent the orbital of $\beta_i$ that contains the fixed set
of $\alpha$ in $A_1$.  By construction, the $B_i$ all exist, and are
nested (so that $B_i\subset B_j$ whenever $i<j$) so that $E =
\left\{(B_i,\beta_i)\,|\,i\in\Z\right\}$ is an exemplary bi-infinite
tower for $H$ whose indexing respects the ordering on the signed
orbitals of $E$.  \qquad$\diamond$

Suppose $h$, $k\in\ploi$ and they satisfy the properties that whenever
$A$ is an orbital of $h$ and $B$ is an orbital of $k$, and $A\cap
B\neq \emptyset$, then either $A = B$, or $\bar{A}\subset B$, or
$\bar{B}\subset A$.  (Note that if $H$ is a balanced subgroup of
$\ploi$ which admits no transition chains, then any two elements of
$H$ will satisfy these conditions.)  

Under these conditions, we will say that $h$ and $k$ satisfy the
\emph{mutual efficiency condition}\index{condition!mutual efficiency}
if given any orbital of $C$ of $h$ that properly contains an orbital
of $k$, then the support of $k$ in $C$ is contained in a single
fundamental domain of $h$ in $C$, and the symmetric condition that
whenever $D$ is an orbital of $k$ that properly contains an orbital of
$h$, then the support of $h$ in $D$ is contained in a single
fundametal domain of $k$ in $D$.

Note that the intitial containment conditions on the orbitals of $h$
and $k$ above occur for any two elements in any balanced subgroup of
$\ploi$ with no transition chains of length two.  

The following remark follows easily from Lemma
\ref{transitiveElementOrbital} and the nature of the orbital
alignments of the elements of the hypothesies, and is left to the
reader.

\brk

If $h$ and $k$ are elements of a balanced subgroup of $\ploi$
that admits no transition chains of length two, then there are
positive integers $m$ and $n$ so that $h^m$ and $k^n$ satisfy the
mutual efficiency condition.  

\erk

We will use this fact heavily in the remainder. 

The following remark is a simple exercises in the calculus that
partially determines the orbitals of a product of two elements of a
balanced subgroup of $\ploi$ without transition chains of length two.

\brk
\label{productOrbitals} 

Suppose $h$ and $k$ are elements of a subgroup $H$ of $\ploi$,
where $H$ is balanced and admits no transition chains of length two.

\be

\item Suppose $h$ has an orbital $A$ and $k$ has an orbital $B$, and
$A\subsetneq B$, then $hk$ has orbital $B$.

\item Suppose $h$ has an orbital $A$ and $k$ has an orbital $B$, and
$B\subsetneq A$, then $hk$ has orbital $A$.

\item If $C$ is an orbital of $hk$, then either $C$
is an orbital of $h$, and $k$ does not have an orbital properly
containing $C$, or $C$ is an orbital of $k$ and $h$ does not have an
orbital properly containing $C$, or both $h$ and $k$ have orbital $D$ which
contains $C$. 
\ee 
\erk

pf: 

To see the first point, since $H$ admits no transition chains of
length two, no end of $B$ can be contained in an orbital of $h$, so
$B$ is an orbital of $\langle h,k\rangle$.  But $k$ realizes $B$
consistently, so that any element of $\langle h,k\rangle$ which
realizes one end of $B$ must realize $B$.  Therefore since $A$ is
properly contained in $B$, $h$ must not realize any end of $B$.  But
therefore the product $hk$ is non-trivial near the ends of $B$ in $B$,
so $hk$ must actually realize $B$.

The argument for the second point is similar to the argument for the
first point, and will not be given here.

To see the third point, if neither $h$ nor $k$ have an orbital
containing $C$, then some point $x$ in $C$ is not in the support of
either $h$ or $k$ (since $H$ admits no transition chains), so that the
product $hk$ cannot have orbital $C$.  Therefore one of $h$ and $k$
must have an orbital $D$ containing $C$.  In which case the first two points
imply all three cases of the conclusion to the third point.

\qquad$\diamond$

We will also use the
following construction.

\begin{construction}[Double Commutator Operation]
Given two elements $h$, $g\in\ploi$ we can construct a third element
$[g,_2 h] = [[g,h],h]$, which we will refer to as the \emph{double
commutator of $g$ and $h$}\index{commutator!double}\index{[g,_2 h]}.

\end{construction}

Double commutators are nice to understand in the setting of a balanced
subgroup of $\ploi$ with no transition chains of length two.

\bl
\label{dcFacts}

Let $g$, $h\in H$, where $H$ is a balanced group in $\ploi$ with no
transition chains of length two.  Suppose further that $g$ and $h$
satisfy the mutual efficiency condition.  If $f = [g,_2 h]$, then $f$
has the following properties:

\be

\item Every orbital of $g$ properly contained in an orbital of $h$ is
an orbital of $f$.

\item Every orbital of $f$ is properly contained in an orbital of $h$ that
contains (perhaps not properly) an orbital of $g$.

\ee

\el 

pf: 

Suppose $A$ is an orbital of $g$ properly contained in some orbital
$B$ of $h$.  The element $[g,h] = g^{-1}h^{-1}gh$ can be thought of as
the product of $g^{-1}$ with the conjugate $g^h$.  Since the support
of $g$ in the orbital $B$ of $h$ is fully contained in a single
fundamental domain of $h$ in $B$, the orbitals of $g^h$ are completely
disjoint from the orbitals of $g$ in $B$.  In particular the product $[g,h]$
has all the orbitals of $g$, where $[g,h]$ acts as $g^{-1}$ on these
orbitals, and a full disjoint copy of these orbitals where the action
of $[g,h]$ on these orbitals is as $g^h$.  Now the element
$[g,_2 h] = ([g,h])^{-1}h^{-1}[g,h]h$ can be thought of as
the product of $([g,h])^{-1}$ with the conjugate $[g,h]^h$.  But the
only orbitals of $([g,h])^{-1}$ which intersect the orbitals of
$[g,h]^h$ are precisely the orbitals of $g^h$, which are disjoint from
the original orbitals of $g$.  Thus, since $[g,h]^{-1}$ has
all of the orbitals of $g$, $[g,h]^h$ does not have these orbitals,
hence the product $[g,_2h]$ has the orbital $A$.  In particular, we
have proven the first point of the lemma.

Now let $A$ be some orbital of $f=[g,_2 h]$.  $A$ is either disjoint
from the orbitals of $h$, properly contains an orbital of $h$, or is
contained in an orbital of $h$.

Suppose $A\cap Supp(h)=\emptyset$.  In this case, let us think of $f$
as the product of two elements of $H$,
\[
f = [[g,h],h]= (h^{-1})^{[g,h]}\cdot h.
\]
In the situation that $H$ is balanced and contains no transition
chains of length two, any orbital of the product of two elements in
$H$ is actually an orbital from one or the other of the two elements
of the product, or is properly contained in an orbital $B$ which is
common to the two elements of the product, by Remark
\ref{productOrbitals}.  In this case, since $A$ is disjoint from the
orbitals of $h$, $A$ must actually be an orbital of
$(h^{-1})^{[g,h]}$.  But then $A = B[g,h]$ for some orbital $B$ of
$h^{-1}$.  Now $A$ is disjoint from $B$ ($B$ is also an orbital of
$h$), so in particular, $B$ and $A$ are contained in a single orbital
$E$ of $(h^{-1})^gh$, which means $E$ is actually an orbital of
$(h^{-1})^g$.  This implies in turn that $(A\cup B) \subset Dg$, where
$D$ is another orbital of $h^{-1}$.  But now $B$ and $D$ are orbitals
of $h$ contained in an orbital $F$ of $g$, and $B$ is in the image of
the orbital $D$ of $h$ under the function $g$.  But every orbital of
$h$ in $F$ is thrown off the support of $h$ by the action of $g$ by
assumption, providing a contradiction, so that our initial assumption
that $A\cap supp(h)=\emptyset$ must be false.

Now suppose $A$ is an orbital of $f$ which properly contains an
orbital of $h$.  We see immediately then that $A$ is an orbital of
$(h^{-1})^{[g,h]}$.  But then $A$ is either an orbital of $[g,h]^{-1}$
or an orbital of $h^{-1}[g,h]$, or is contained properly in a common
orbital $B$ of both of these elements.  In the first and second cases,
$A$ would actually be an orbital of both elements, since $h^{-1}$
does not have any orbitals in $A$ near the ends of $A$.  But the
slopes of $[g,h]^{-1}$ and $h^{-1}[g,h]$ near the ends of $A$ in $A$
are inverse, so that $A$ cannot actually be an orbital of the product
$[g,h]^{-1}\cdot h^{-1}[g,h]$.  In particular, $A$ must be properly
contained in some common orbital $B$ of both $[g,h]^{-1}$ and
$h^{-1}[g,h]$.  But then $B$ is an orbital of $h^{-1}[g,h]$ induced
from an orbital of $(h^{-1})^g$ by the action of $h$, which cannot
move the ends of $B$ (h has an orbital in $A$, and therefore in $B$,
and so has no orbital containing the ends of $B$), so that $B$ is
actually an orbital of $(h^{-1})^g$.  But every orbital $(h^{-1})^g$
is disjoint from the orbitals of $h^{-1}$, so that $B$ could not
contain an orbital of $h$.

From the above, we now see that $A$ is contained in an orbital $B$ of
$h$, but we do not yet know that $B$ contains an orbital of $g$.
Suppose firstly that $B$ is disjoint from every orbital of $g$.  In
this case it is immediate that $f$ is actually trivial on $B$, so that
$A$ cannot be an orbital of $f$.

Now we know that $A$ is contained in an orbital $B$ of $H$ which
non-trivially intersects an orbital $C$ of $g$.  If $\bar{C}$ is
contained in $B$ then it is immediate that $\bar{A}\subset B$ and we
are done, so suppose instead the either $C=B$ or $\bar{B}\subset C$.

If $C=B$, then again, the slopes of $g$ and $h$ cancel near the ends
of $B$, so that $\bar{A}\subset B$, and again we are done.  Therefore
let us assume that $\bar{B}\subset C$.

In this case, given any $x\in A\subset B$, we can directly calculate the impact
of $f$ on $x$ to see that $xf = x$.  Let us write $f$ as a product to
see this.
\[
\begin{array}{l}
xf = x(g^{-1}h^{-1}gh)^{-1}h^{-1}(g^{-1}h^{-1}gh)h
=xh^{-1}\cdot(g^{-1}hg)\cdot h^{-1}\cdot(g^{-1}h^{-1}g)\cdot hh=\\
x(h^{-1}h^{-1})\cdot(g^{-1}h^{-1}g)\cdot(hh)=x(h^{-1}h^{-1})\cdot(hh) = x
\end{array}
\]
In the above the expressions with $\cdot$'s the parenthesies have been
placed suggestively to help us understand the dynamics of the point
$x$.  Any such factor in an above product involving a ``$\cdot$''
represents a factor with no net effect on the point it acts on.

Thus, if $A$ is an orbital of $f$, then $A$ is properly contained in
an orbital $B$ of $h$, and $B$ contains an orbital $C$ of $g$.
\qquad$\diamond$

\subsection{Finding infinite wreath products in groups with infinite towers}

Suppose $D = \left\{(A_i,h_i)\,|\,i\in\N\right\}$ is an exemplary
tower whose indexing respects the order of the elements so that $H =
\langle S_D\rangle$ is a balanced group that admits no transition
chains of length two, and so that whenever $B$ is an orbital of $h_i$
for some signature $h_i$ of $D$, then $B$ is contained in an orbital
$C$ of $h_{i+1}$.  We are going to find a subtower of $D$ that
satisfies a nice further property.

Suppose $B_1$ is an orbital of $h_1$.  Each signature $h_i$ of $D$ has
an orbital $B_i$ that contains $B_1$.  The orbitals $B_j$ are nested
as the index increases, but possibly not properly.  If there is an
$N_1\in\N$ so that for all $n>N$, we have $B_n = B_{n+1}$, then we
will call $B_1$ \emph{a terminal orbital of
$D$}\index{orbital!terminal}, and $(B_1,h_1)$ a terminal signed
orbital of $D$, and we will say that \emph{$B_1$ is stable after
$N_1$}.  We now extend this language to orbitals of signatures other than
$h_1$.  Given $i\in\N$, call an orbital of $h_i$ terminal in $D$ if
the orbital is terminal in the subtower of $D$ formed using only the
signed orbitals $(A_k,h_k)$ with $k\geq i$.  We will call the orbital
of any signature of the tower, where the orbital is not a terminal
orbital, a non-terminal orbital.  Observe that non-terminal signed
orbitals make good candidates for being bases of new exemplary towers.

We will rely heavily on the following technique in our proof of Lemma
\ref{tallZWreath}.

\bl[growing subtower] 
\label{growingTower}

Suppose $D = \left\{(A_i,h_i)\,|\,i\in\N\right\}$ is an exemplary
tower so that $H = \langle S_D\rangle$ is a balanced group that admits
no transition chains of length two, and so that whenever $B$ is an
orbital of $h_i$ for some signature $h_i$ of $D$, then $B$ is
contained in an orbital $C$ of $h_{i+1}$.  Then we can pass to an
infinite subtower $E$ of $D$ so that if $J$ is an orbital of any
signature $g_i$ of $E$, where $J$ is not a terminal orbital of $D$,
then there is an orbital $K$ of $g_{i+1}$ which properly contains $J$.
\el

We note by definition that the orbital $K$ will also be a non-terminal
orbital of $D$, and both will be non-terminal in $E$.

\underline{Proof of Lemma:}

We now pass repeatedly to infinite subtowers for $D$, at each stage
referring to the new tower that results as $D$, and re-indexing so
that the tower will still have the form $D =
\left\{(A_i,h_i)\,|\,i\in\N\right\}$.  Let $\mathscr{P} =
\left\{B_i\,|\,1\leq i\leq n_i, i\in\N\right\}$ represent the $n_i$
orbitals of $h_1$ that are not terminal, in left to right order.  We
improve $D$ by passing to a infinite subtower $n_i$ times.  Firstly,
for $B_1$, pass to a subtower of $D$ so that the orbitals of the $h_i$
over $B_1$ are always properly nested as we progress up the tower.
The new tower $D$ still has all the properties that we have listed for
the old $D$, but now the orbitals of the $h_i$ over $B_1$ actually
form a tower over $B_1$ when we pair them with their signatures.
Repeat this process inductively for each of the non-terminal orbitals
of $D$ in $h_1$.  Now we pass to an infinite induction, by repeating
the process again, using base signature $h_2$, so that we are
progressively improving the tower above $h_2$ so that the non-terminal
orbitals of $h_2$ are each actually the base of an infinite tower
using the signatures $h_k$ with $k>2$ paired with their appropriate
orbital containing the relevant orbital of $h_2$.  We note in passing
that the non-terminal orbitals of $h_1$ are all contained in the
non-terminal orbitals of $h_2$, and $D$ is already a perfect tower
with respect to these operations over the non-terminal orbitals of
$h_2$ which contain non-terminal orbitals of $h_1$, so that we will
only have to improve $D$ over the non-terminal orbitals of $h_2$ which
do not contain orbitals of $h_1$.  With these observations in place,
we can inductively continue this process at every level of $D$.  We
have now defined a new tower $D$.  Given any $i\in\N$, if $(A_i,h_i)$
is a signed orbital of $D$ and $B_i$ is a non-terminal orbital of
$h_i$ then for any integer $k>i$ there is an orbital $B_k$ of $h_k$ so
that $\bar{B}_i\subset B_k$.  Thus, $D_{B_i} =
\left\{(B_k,h_k)\,|\,k\geq i, k\in\N\right\}$ is itself an exemplary
tower.  \qquad$\diamond$

The following lemmas are simply restatements (with proofs) of Theorems
\ref{tallWreath} and Theorem \ref{deepWreath}.

\bl
\label{tallZWreath}

If $G$ is a subgroup of $\ploi$ and $G$ admits a tall tower, then $G$
has a subgroup of the form $(\Z\wr)^{\infty}$.

\el

pf:

We will assume that $G$ is balanced, as otherwise $G$ contains a
subgroup isomorphic to Thompson's group $F$, which itself has a
subgroup isomorphic to $(\Z\wr)^{\infty}$.  If $G$ admits transition
chains of length two then by Lemma \ref{messyWreath} $G$ admits an
exemplary bi-infinite tower $E =\left\{(A_i,g_i)\,|\,i\in\Z\right\}$
where the indexing respects the order of the signed orbitals, and
where the group $H = \langle S_E\rangle$ admits no transition chains
of length two.  In particular, we can pass to the subgroup $K$
generated by all the signatures with positive index to reproduce the
hypothesies of this lemma, with the extra condition that $G$ admits no
transition chains of length two, so let us assume that $G$ admits no
transition chains of length two.

Let $E = \left\{(A_i,g_i)\,|\,i\in\N\right\}$ be a tall tower for $G$,
where the indexing respects the order on the signed orbitals
of $E$.  Since $G$ is balanced, and contains no transition chains of
length two, we see that all towers of $G$ are exemplary,
and in particular, $E$ is exemplary.

Let $A = \cup_{i\in\N}A_i = (a,b)$.  We observe that if $B$ is an
orbital of $g_i$ for some $i$, then $B$ is disjoint from
$\left\{a,b\right\}$, and that if $B\cap A \neq\emptyset$, then
neither $a$ nor $b$ is an end of $B$.  In particular, $A$ is an
orbital of $\langle S_E \rangle$.

Now given $\epsilon>0$ so that $\epsilon< \frac{b-a}{2}$, we see that
there is an $N\in\N$ so that for all $n\in\N$ with $n\geq N$, we have
that $(a+\epsilon,b-\epsilon)\subset A_n$ since the ends of the $A_i$
must limit to the ends of $A$.  But now, we can construct a monotone
strictly increasing, order preserving function, $\phi:\N\to\N$, so
that given any $n\in\N$, all of the orbitals of $g_n$ in $A$ are
actually contained in $A_{\phi(n)}$, and since $E$ is exemplary, no
orbital of $g_n$ in $A$ actually shares an end with $A_{\phi(n)}$.  For any
$k\in\N$, let $\phi^k$ represent the product (via composition) of the
function $\phi$ with itself $k$ times in the monoid of order
preserving functions from $\N$ to $\N$.  Now define an order
preserving function $\theta:\N\to\N$, defined by the rules that
$1\mapsto 1$ and $n\mapsto \phi^{n-1}(1)$ for each
$n\in\N\backslash\left\{1\right\}$.  Replace $E$ by the exemplary
tower formed by the collection
$\left\{(A_{\theta(i)},g_{\theta(i)})\,|\,i\in\N\right\}$.  $E$ now
has the property that if $i$, $k\in\N$ with $i<k$ then all the
orbitals of $g_i$ in $A$ are actually in $A_k$, away from the ends of
$A_k$.  For each $n\in\N$, with $n>1$, let $m_n$ be an integer large
enough so that the collection of orbitals of $g_{n-1}$ inside of $A_n$
(which is all the orbitals of $g_{n-1}$ in $A$) is actually fully
contained in a single fundamental domain of $g_n^{m_n}$ in $A_n$.
Define $n_1 = 1$.  Improve $E$ by replacing each signature $g_n$ with
$g_n^{m_n}$.  Now define $H = \left<S_E\right>$.  We note in passing
that $A$ is an orbital of $H$.

We will now improve $E$ further.  Define $h_1=g_1$.  Now for each
$n\in\N$ with $n >1$, inductively define $h_n$ via the following four
step process.  

First, define $k_n = g_n^{r_n}$, where $r_n$ is a
positive integer large enough so that whenever $B$ is an orbital of
$g_n$ that is also an orbital of $h_{n-1}$, then the product
$h_{n-1}k_n$ still has orbital $B$.  

Second, define $h_n' = h_{n-1}k_n$.  Recall from Remark
\ref{productOrbitals} that any orbital of $h_{n-1}$ which properly
contains an orbital of $k_n$ will now be an orbital of $h_n'$,
and that any orbital of $k_n$ that properly contains an orbital of
$h_{n-1}$ will also be an orbital of $h_n'$.

$h_n'$ now has an orbital containing every orbital of $h_{n-1}$.

Third, choose positive integer $s_n$ large enough so that
every orbital $C$ of $h_n'^{s_n}$ which properly contains orbitals of
$h_{n-1}$ actually contains all such orbitals in a single fundamental
domain of $h_n'^{s_n}$ on $C$.  

Fourth, define $h_n = h_n'^{s_n}$.  The result is that the sequence
$(h_i)_{i\in\N}$ of signatures satisfies the following list of
properties.

\be

\item For each $n\in\N$, $A_n$ is an orbital of $h_n$.
\item For each $n\in\N$ with $n>1$, the orbitals of $h_{n-1}$ in $A$
are all contained inside a single fundamental domain of $h_n$ in
$A_n$.
\item For each $n\in\N$ with $n>1$, if $B$ is an orbital of $h_n$
which is not disjoint from the orbitals of $h_{n-1}$, then there are
two possibilities.

\be

\item $B$ is also an orbital of $h_{n-1}$.
\item $B$ properly contains a non-empty collection of orbitals of
$h_{n-1}$ in a single fundamental domain of $h_n$ on $B$.  

\ee 
\ee 

In particular, we can form the new exemplary tower $D =
\left\{(A_i,h_i)\,|\,i\in\N\right\}$.

$D$ still has the properties that $\cup_{i\in\N} A_i = A$, and that
$A$ is an orbital of the group $\langle S_D \rangle$.  Further, the
signatures satisfy all the properties of the last paragraph.

We will now improve $D$ by replacing it with the result of finding
a growing subtower, so that any non-terminal orbital of any signature
$h_i$ of $D$ is properly contained in a non-terminal orbital of a
signature with index one higher.

Our new $D$ is far superior to our old $D$, but $h_1$ will still have
terminal orbitals, if it had them to begin with.  Suppose $h_1$ does
have some terminal orbitals.  Then there is $N_1\in\N$ so that all the
terminal orbitals of $h_1$ are stable for $n\geq N_1$.  Compute a new
element $k = [h_{N_1+1},_2 h_{N_1+2}]$ (note that condition (3) above is
equivalent to saying that $h_{N_1+1}$ and $h_{N_1+2}$ satisfy the
mutual efficiency condition when all the orbitals of $h_j$ are
contained in orbitals of $h_{j+1}$ for any $j\in\N$). $k$ has the
following properties.

\be

\item The orbitals of $h_{N_1}$ which contain the terminal orbitals of
$h_1$ are not contained in the orbitals of $k$.

\item No orbital of $h_{N_1}$ which is also an orbital of $h_{N_1+1}$
is also an orbital of $k$ (these are all terminal orbitals of
$h_{N_1}$ since $D$ is the result of using a growing tower operation).

\item All the non-terminal orbitals of $h_{N_1}$ are still properly
contained in the orbitals of $k$ since $k$ contains the non-terminal
orbitals of $h_{N_1+1}$.

\ee

Now replace $k$ and $h_{N_1}$ by sufficiently high powers of
themselves so that they satify the mutual efficiency condition and let
$h = [h_{N_1},_2 k]$.  The resulting $h$ has the following properties.

\be
\item $h$ has no orbitals intersecting the terminal orbitals of $h_1$.
\item $h$ has all the non-terminal orbitals of $h_{N_1}$
\ee

Now replace $h$ and
$h_1$ by sufficiently high powers of themselves so that the satisfy
the mutual efficiency condition, and replace $h_1$ by $[h_1,_2 h]$.
Now replace $h_1$ and $h_{N_1}$ by sufficiently high powers of
themselves so they satisfy the mutual efficiency condition.  Build the
tower
\[
D' = \left\{(A_1,h_1)\right\}\cup \left\{(A_i,h_i)\,|\,i\geq N_1,
i\in\N\right\}.
\]

In this tower, $h_1$ admits only non-terminal orbitals, every
orbital of $h_1$ is properly contained in a non-terminal orbital of
$h_{N_1}$, and $h_1$ still has a copy of every non-terminal
orbital that it started with.  If we re-index the tower $D'$ and call
it $D$ again, then it satisfies all the old properties of the tower
$D$ found above, but its bottom element ($h_1$) has nice orbitals.  We
can now repeat this whole process for the subtower of $D$ starting
from level two and up, so that the new $h_2$ will admit all the
non-terminal orbitals that it started with, and other non-terminal
orbitals, and also will contain no terminal orbitals.  Inductively
proceed up the tower $D$, redefining all of the $h_i$, so that the new
tower $D$ satisfies the following properties.
\be

\item $A = \cup_{n\in\N} A_n$
\item For each $n\in\N$ with $n>1$, the orbitals of $h_{n-1}$ in $A$
are all contained inside the orbital $A_n$ of $h_n$.
\item For each $n\in\N$ with $n>1$, if $B$ is an orbital of $h_n$
which is not disjoint from the orbitals of $h_{n-1}$, then $B$
contains the closure of the union of the collection of orbitals of
$h_{n-1}$ that intersect $B$.

\ee 

Now for each index $j\in\N$, inductively replace $h_j$ and $h_{j+1}$
by sufficiently high powers of themselves so that they satisfy the
mutual efficiency condition.  (This is actually unnecessary by the
details of the proof of Lemma \ref{dcFacts}, but it hurts nothing
and explicitely guarantees that all adjacent pairs of signatures of the
tower $D$ satisfy the mutual efficiency condition.)  Note that each
signature (except $h_1$) may be replaced by progressively higher
powers of itself twice in this operation, but that once two signatures
are mutually efficient, replacing either signature by a higher power
of itself will still result in a pair that are mutually efficient.

Now any pair of adjacent signatures of the tower $D$ satisfy the
mutual efficiency condition.

For every $n\in\N$, define the subgroup $H_n=\langle
h_1,h_2,\ldots,h_n\rangle$ of $G$.  Given any two elements $f$, $g\in
H_n$, since the supports of $f$ and $g$ are contained in the support
of $h_n$, and since the support of $h_n$ in any one orbital of
$h_{n+1}$ is contained in a single fundamental domain of $h_{n+1}$ in
that orbital, we see that $f^{h_{n+1}^j}$ and $g^{h_{n+1}^k}$ have
disjoint supports and therefore commute, whenever $j\neq k$.  If
$j=k$, then the product of the conjugated $f$ and $g$ is equal to the
conjugate of the product of $f$ and $g$.  In particular, the group of
finite products of conjugates of elements of $H_n$ by $h_{n+1}$ is
isomorphic to $\sum_{j\in\Z}H_n$, where the indexing factor $j$
represents the power of $h_{n+1}$ used in the conjugation of the
element from $H_n$ under consideration.  But we can write any element
of $H_{n+1}$ as a product of an integer power of $h_{n+1}$ with a
product of conjugates of elements of $H_n$ by integer powers of
$h_{n+1}$; in short, $H_{n+1}\cong H_n\!\!\wr\!\Z$, where the $\Z$
factor is the subgroup $\langle h_{n+1}\rangle$ of $H_{n+1}$.

Now, $H_1 \cong \Z$, so $H_2\cong \Z\wr\!\Z$, $H_3\cong
(\Z\!\wr\!\Z)\!\wr\!\Z$, and etc., so that $H_n\cong
((\cdots(\Z\!\wr\!\Z)\!\wr\!\Z)\cdots\wr\!\Z$ where the finite wreath product
has $n$ factors of $\Z$.  In particular, the ascending union
$H=\langle h_1,\,h_2,\,\ldots \rangle\cong(\Z\wr)^{\infty}$.

\qquad$\diamond$

\bl
\label{deepWreathZ}

If $G$ is a subgroup of $\ploi$ and $G$ admits a deep tower, then $G$
has a subgroup of the form $(\wr\Z)^{\infty}$.

\el

pf:

We will use a similar technique to the proof of Lemma
\ref{tallZWreath}, although the analysis in this case is much simpler.

We will assume that $G$ is balanced, as otherwise $G$ contains a
subgroup isomorphic to Thompson's group $F$, which itself has a
subgroup isomorphic to $(\Z\wr)^{\infty}$.  If $G$ admits a transition
chain of length two, then by Lemma \ref{messyWreath} $G$ admits an
exemplary bi-infinite tower $E =\left\{(A_i,g_i)\,|\,i\in\Z\right\}$
where the indexing respects the order of the signed orbitals, and
where the group $H = \langle S_E\rangle$ admits no transition chains
of length two.  In particular, we can pass to the subgroup $K$
generated by all the signatures with negative index to reproduce the
hypothesies of this lemma, with the extra condition that $G$ admits no
transition chains of length two.  In particular, we can assume that
$G$ admits no transition chains of length two.

Since $G$ is balanced and admits no transition chains of length two,
any tower for $G$ is exemplary.  In particular, let $E =
\left\{(A_{-i},g_{-i})\,|\,i\in\N\right\}$ be an exemplary deep tower
for $G$ where the indexing respects the order on the elements
of the tower.  Improve $E$ by replacing the signatures of $E$ with
sufficiently high powers of themselves so that given any negative
integer $i$, then $g_{i-1}$ and $g_i$ satisfy the mutual efficiency
condition.

Let $A=A_{-1}=(a,b)$.  Since $E$ is exemplary, we see that $A$ is
actually an orbital of the subgroup $H\leq G$, where $H =
\left<S_E\right>$.  For all $i\in\N$ with $i>1$, inductively improve
$E$ (induct on increasing $i\in\N$ in the following discussion) by
replacing the signatures of $E$ according to the following three step
process.  

First, let $h_{-i} = [g_{-i},_2 g_{-i+1}]$.  

Second, define the new $g_{-i}$ to be $h_{-i}$.  

Third, replace the elements $g_{-i+1}$, $g_{-i}$, and $g_{-i-1}$
with sufficiently high powers of themselves, so that given any index
$j\in\N$, the elements $g_{-j}$ and $g_{-j-1}$ satisfy the mutual
efficiency condition (observe that if $i>3$, then $g_{-i+1}$ and
$g_{-i+2}$ will now still satisfy the mutually efficiency condition,
since we are only replacing $g_{-i+1}$ by higher powers of itself, and
these two signatures were already mutually efficient, a similar
argument shows that $g_{-i-1}$ and $g_{-i-2}$ will be mutually
efficient after this operation as well).  

Since $A_{-i}\subsetneq A_{-i+1}$ for all integers $i>1$, we see that the
resultant set of signed orbitals is still a tower (and with the same
order), so that this inductive definition will simply improve our
tower $E$.  Observe further that given any $k\in\N$, then the orbitals
of $g_{-k-1}$ are all properly contained in the orbitals of $g_{-k}$.

Define the set $\Gamma_i=\left\{g_j\,|\,j\leq i, j\in\Z\right\}$ for
each negative integer $i$.  For each negative integer $i$, define
$H_i=\langle \Gamma_i\rangle$.  For such $i$, the orbitals of $H_i$
are actually the orbitals of $h_i$, since all orbitals of the elements
$g_{k}$ with $k<i$ are contained in the orbitals of $g_i$.
Furthermore, for any such $i<-1$, the orbitals of $g_i$ are contained
in the orbitals of $g_{i+1}$ in such a way that in any individual
orbital $B$ of $g_{i+1}$, the support of $g_{i}$ in $B$ is actually
fully contained inside a single fundamental domain of $g_{i+1}$ on
$B$.  In particular, $H_i\cong H_{i-1}\!\wr\!\Z$, where the $\Z$
factor comes from the subgroup $\langle g_i\rangle$ of $H_i$.  But now
inductively, since each generator generates a group isomorphic to
$\Z$, we see that $H_1\cong(\wr\Z)^{\infty}$.

\qquad$\diamond$

The following lemma is only surprising in the sense that its proof is
somewhat involved.  It belongs in this subsection since (as will be
shown) finitely generated non-solvable subgroups of $\ploi$ always
contain infinite towers.  This lemma, together with Lemmas
\ref{deepWreathZ} and \ref{tallZWreath}, completes the proof of Theorem
\ref{fgNonSolve}.

\bl

If $H$ is a finitely generated and $H$ admits towers of arbitrary
height, then $H$ admits an infinite tower.

\el

Pf:

Suppose $H$ is a finitely generated subgroup of $\ploi$ and $H$ admits
a towers of arbitrary height.  Let $n$ be the smallest integer so that
$H$ has a generating set of size $n$.  Let $\Gamma =
\left\{g_1,g_2,\ldots,g_n\right\}$ be a minimal set of generators for $H$,
where the sum of the number of the orbitals of elements of $\Gamma$ is
also minimal amongst all generating collections for $H$ of size $n$.
Since $H$ has towers of arbitrary height, $n$ is at least two, and $H$
has at least one orbital.

If $H$ is imbalanced, then $H$ admits a bi-infinite tower by Lemma
\ref{imbalancedTower}.  So we will assume that $H$ is balanced.

If $H$ admits a transition chain of length two then $H$
admits a bi-infinite exemplary tower, so we will assume that $H$ has
no transition chains of length two.  In particular, if $h$,
$g\in H$, $A$ is an orbital of $g$ and $B$ is an orbital of $h$,
where $A\cap B\neq\emptyset$, then either $A=B$ or $\ol{A}\subset B$
or $\ol{B}\subset A$ by Remark \ref{nestedOrbitals}.

Now if $H$ has an element with a deep orbital, we can use the orbital
nesting properties of Remark \ref{nestedOrbitals} to inductively build
either a tall tower or a deep tower.  In particular, we will now
assume that no orbital of any element of $H$ has infinite depth.  In
particular, the depth of each orbital of each element of $H$ is well
defined and is finite.

We note that since $H$ admits towers of arbitrary length, it must have
elements of arbitrary finite depth.

Given any tower $T$, let us call the union of its orbitals the support
of the tower.  In our case, since there are no transition chains of
length two, and since $H$ is balanced, the group generated
by the signatures of $T$ will have an orbital equal to the support of
$T$.  If the support of a tower $T$ is contained in an orbital $A$ of
an element of $h\in H$, let us say that $A$ supports $T$ and also that
$h$ supports $T$.  If we assume that all towers of $H$ are finite,
then they are always supported by their highest signature, and their
highest orbital.

We will now inductively define a sequence of finitely generated groups
$(H_i)_{i=1}^{\infty}$ so that $H_1 = H$ and for any $i\in\N$ we have
the following properties:

\be

\item $H_{i+1}\leq H_i$.

\item There is $m_i\in\N$ and $H_i$ has a finite generating set $\Gamma_i = \left\{g_{ik}|k\in\N, 1\leq k\leq m_i\right\}$.

\item $H_i$ has an orbital $A_i$ and only the generator $g_{i1}$ in
$\Gamma_i$ realizes $A_i$.

\item $H_i$ admits a collection $X_i$ of towers of
arbitrary height for $H_i$, where all the towers of $X_i$ contain the
signed orbital $(A_i,g_{i1})$ as their largest element.

\item Every tower of $X_k$ is supported by $A_i$, whenever $k\in\N$ and $k>i$.

\item The orbital $A_{i+1}$ of $H_{i+1}$ satisfies
$\overline{A}_{i+1}\subset A_i$.

\ee

First, let $H_1 = H$.  Since $H_1$ is finitely generated, it has only
finitely many orbitals.  Since every tower for $H_1$ is supported by
one of the orbitals of $H_1$, one of the orbitals of $H_1$ supports
towers for $H_1$ of arbitrary height.  Define $A_1$ to be such an
orbital of $H_1$.  Now since $H_1$ is balanced, $H_1$ admits a
controller for $H_1$ on $A_1$.  Let $g_{11}$ be such a controller for
$H_1$ on $A_1$.  For each $k\in\N$, where $1\leq k\leq n$, the
generator $g_k$ in $\Gamma$ has a $g_{11}$-form
$g_{11}^{n_{1k}}h_{1k}$ for $A_1$, where the orbitals of $h_{1k}$ in $A_1$ are
all properly contained in $A_1$.  For each
$k\in\left\{1,2,\ldots,n\right\}$, define $g_{1(k+1)} = h_{1k}$.  Now
if we set $m_1 = n+1$ and define the set $\Gamma_1 =
\left\{g_{1k}|1\leq k\leq m_1, k\in\N\right\}$, we see that the set
$\Gamma_1$ generates $H_1$, but only $g_{11}\in\Gamma_1$ actually
realizes the orbital $A_1$.

Given any tower of $H_1$ supported by $A_1$, either the tower has a
signed orbital of the form $(A_1,f)$ as a top element (there can be
none higher) for some element $f\in H_1$, in which case we can replace
this element by $(A_1,g_{11})$, or the tower does not have such an
element, in which case we can just add the element $(A_1,g_{11})$ to
the tower.  The resultant tower from either of the two processes above
will have height unchanged, or one greater than before, and it will
have $(A_1,g_{11})$ as its top element.  In particular, we can do this
for each of the towers for $H_1$ supported by $A_1$ (which set of
towers admits towers of arbitrary height).  Now define $X_1$ to be the
set of all towers for $H_1$ which have top element $(A_1,g_{11})$.  By
our discussion, we see that the collection $X_1$ of towers of $H_1$
with top element $(A_1,g_{11})$ is a collection with towers of
arbitrary height.

We see that $H_1$ satisfies all of the requirements.

Suppose $(B,\tilde{f})$ is a signed orbital of $H_1$ of depth two.
Since $\tilde{f}\in H_1$, there is an element $f\in H_1$ so that $f$
agrees with $\tilde{f}$ on $B$, where $f$ admits a shortest length
decomposition as a product of elements of $\Gamma_1$ and their
inverses, out of all elements in $H_1$ which agree with $\tilde{f}$ on
$B$.  Let $m$ be the length of such a shortest length product, and let
$f = \prod_{i = 1}^mk_i$ be such a product.  Now $B$ is also an
orbital of $f$ and $f|_B = \tilde{f}|_B$, so we will transfer
consideration to the signed orbital $(B,f)$.

Consider the sequence of functions $(f_j)_{j=1}^m$ defined by the rule
\[
\begin{array}{l}

f_j = (k_1k_2\cdots
k_{j-1})^{-1}f(k_1k_2\cdots k_{j-1})=

\\

(k_1k_2\cdots
k_{j-1})^{-1}(k_1k_2\cdots k_m)(k_1k_2\cdots k_{j-1}) =
\\
k_jk_{j+1}\cdots k_mk_1k_2\cdots k_{j-1},
\end{array}
\]

and the corresponding sequence of signed orbitals $((B_j,f_j))_{j =
1}^{m}$ where for each $j$ the orbital $B_j$ is the signed orbital of
$f_j$ induced from $B$ by the action of $\prod_{i = 1}^{j-1}k_i$ on
$f$.  Note that the depth of each orbital $B_j$ is two, and $B_j\subset A_1$.

For each $j\in\left\{1,2,\ldots,m\right\}$, we must have that $B_j$
has nontrivial intersection with an orbital of $k_j$, or the product
decomposition could be shortened to produce a different element $g$
which is identical to $f$ (and therefore $\tilde{f}$) on $B$.

 Since $f$ does not have orbital $A_1$, we see that the number of
times that $g_{11}$ and the number of times that $g_{11}^{-1}$ appear
in the product decomposition of $f$ above must be the same, and since
$f$ is not trivial, some other elements from $\Gamma_1$ (or their
inverses) must appear in the product decomposition of $f$.  Let
$S=\left\{i\in\N|1\leq i\leq m, k_i\neq g_{11}^{\pm 1}\right\}$.  By
the previous comment, $S$ is not empty, let $t$ represent the size of
$S$.  Let $\alpha:\left\{1,2,\ldots,t\right\}\to S$ be a monic
increasing function, so that $\alpha(1)$ is the smallest element in
$S$, and $\alpha(t)$ is the largest, ie., $\alpha$ is an ordered
indexing of the elements in $S$.

Let $j\in\left\{2,3,\ldots,t\right\}$.  By the definition of an induced
orbital, we see that $B_j = Bk_1k_2\cdots k_{j-1}$, and that $B_1 =
B$.  In particular, for $i\in\left\{2,3,\ldots,t\right\}$,

\[
B_{(\alpha(i))}=B_{(\alpha(i-1))}k_{(\alpha(i-1))}k_{(\alpha(i-1)+1)}k_{(\alpha(i-1)+2)}\cdots
k_{(\alpha(i)-1)},
\]
and $B_{\alpha(1)} = Bk_1k_2\cdots k_{(\alpha(1) - 1)}$ if $\alpha(1)
\neq 1$ and$B_{\alpha(1)} = B$ if $\alpha(1) = 1$.

\vspace{.1 in}
\underline{{\bf Claim:}}
For some $i\in S$, we have that $B_{\alpha(i)}$ is
actually an orbital of $k_{\alpha(i)}$, so that $(B,f)$ is conjugate
to $(B_j,f_j)$ for some index $j$, where $B_j$ is actually an orbital
of one of the generators of $\Gamma_1$.  

Consider $f_{\alpha(1)}=k_{\alpha(1)}\cdot v_1$, where $v_1$ is the
product of the remaining terms of the product decomposition of
$f_{\alpha(1)}$ after the first term.  Now, $k_{\alpha(1)}\neq
g_{11}^{\pm 1}$ by construction.  By Lemma \ref{nestedOrbitals} and
Remark \ref{productOrbitals} we must have that $B_{\alpha(1)}$ is
either an orbital of $k_{\alpha(1)}$, in which case we have our claim,
or else it simply properly contains a smaller orbital of
$k_{\alpha(1)}$ and is in fact an orbital of $v_1$.  (Note,
$k_{\alpha(1)}$ cannot have an orbital which contains $B_{\alpha(1)}$
since $B_{\alpha(1)}$ is a depth two orbital contained only in the
depth one orbital $A_1$.)  Suppose therefore that $B_{\alpha(1)}$ is
an orbital of $v_1$ and that $B_{\alpha(1)}$ properly contains an
orbital of $k_{\alpha(1)}$.  Both of the products
\[
\begin{array}{l}
\tilde{\beta_1} = \prod_{j= \alpha(1)}^{\alpha(2)-1} k_j
\\
\textrm{and}
\\
\beta_1=\prod_{j =\alpha(1)+1}^{\alpha(2)-1}k_j
\end{array}
\]

will take $B_{\alpha(1)}$ to $B_{\alpha(2)}$, so the conjugate of $v_1$
by $\beta_1$ will produce a new function 

\[
w_2 = k_{\alpha(2)}\cdots
k_mk_1\cdots k_{(\alpha(1)-1)}k_{(\alpha(1)+1)}k_{(\alpha(1)+2)}\cdots
k_{(\alpha(2)-1)}
\]

which is like $f_{\alpha(2)}$ in that it still has orbital
$B_{\alpha(2)}$, and it is the same product as $f_{\alpha(2)}$
excepting the term $k_{\alpha(1)}$ is deleted from the product.  Now we
can write $w_2 = k_{\alpha(2)}\cdot v_2$, where $v_2$ is the product
in the definition of $w_2$, excepting the first term.  As before,
$B_{\alpha(2)}$ is either an orbital of $k_{\alpha(2)}$ or it is an
orbital of $v_2$.  If it is an orbital of $k_{\alpha(2)}$ we are
finished with the claim, otherwise we will inductively follow this
procedure, defining elements $w_{j}$ and $v_j$ along the way for each
$j\in\left\{2,3\ldots,t\right\}$, and checking whether $B_{\alpha(j)}$
is an orbital of $k_{\alpha(j)}$ for each
$j\in\left\{1,2,\ldots,t\right\}$.  If this condition is never
satisfied, then $B_{\alpha(t)}$ must be an orbital of $v_t$.  But
$v_t$ is a product purely of elements which are all $g_{11}$ and its
inverse, where the number of each type is the same, so that $v_t$ is
actually the identity!  In particular, we must have that for some
$i\in S$, $B_i$ is an orbital of $k_i$, so that $(B,f)$ is conjugate
to $(B_i,f_i)$, where $f_i = k_i\cdot u_i$, where $k_i\in\Gamma_1$, or
$k_i^{-1}\in\Gamma_1$, and $B_i$ is an orbital of depth two for $H_1$
which is also an orbital of $k_i$.  The claim is proven.
\vspace{.1 in}

With this last claim in hand, we see that the elements of $\Gamma_1$
actually contain orbitals of depth two, and there are finitely many of
them.  Further, each orbital of depth two in $A_1$ for $H_1$ is
conjugate to some element of this finite list.  In particular, every
tower for $H_1$ supported by $A_1$ can be conjugated to a tower
supported by a signed orbital $(B,g)$ of depth two where
$g\in\Gamma_1$.  Now there are towers of arbitrary height for $H_1$,
but there are only finitely many signed orbitals of the form $(A,g)$
where $g\in\Gamma_1$, so at least one element of $\Gamma_1$ has an
orbital of depth two that supports towers of arbitrary height beneath
it.  Let $(A_2,g)$ be a signed orbital of $H_1$, where $g\in\Gamma_1$,
and $A_2$ is an orbital of depth two which supports towers of
arbitrary height for $H_1$ beneath it.

Define
\[
\Lambda_2 = \left\{s\!\in\! H_1\,|\,s = g_{11}^{-k}hg_{11}^k,\,
Supp(s)\!\cap\!
A_2\neq\emptyset,\,h\!\in\!\Gamma_1\!\backslash\!\left\{\!g_{11}\!\right\},\,
k\in\Z\right\}.
\]
the set of all conjugates of elements of $\Gamma_1$ (aside from
$g_{11}$) which have resulting support intersecting nontrivially with
$A_2$.  This set is finite, since after multiple conjugates of any
element of $\Gamma_1$ (other than $g_{11}$) by $g_{11}$, the result
has all orbitals near an end of $A_1$, away from $A_2$.

We will define $H_2$ to be the group generated by $\Lambda_2$.  Note
that $A_2$ is an orbital of $H_2$ and is realized by a generator in
$\Lambda_2$.

\vspace{.1 in}

\underline{{\bf Claim:}}

The collection $\hat{X}_2$ of towers for $H_2$ all of which have
largest element of the form $(A_2,h)$ for various $h\in H_2$ contains
towers of arbitrary height for $H_2$.

\vspace{.1 in}

In order to establish this, we show that if $(A,f)$ is a signed
orbital of $H_1$, where $A\subset A_2$, then there is an element $g\in
H_2$ so that $A$ is an orbital of $g$.

Suppose therefore that $(A,f)$ is a signed orbital of $H_1$ where
$A\subset A_2$.  Let $A_f = \left\{g\in H_1\,|\,g|_A = f|_A\right\}$.
Let $g\in A_f$ be an element which can be written as a product of
elements in $\Gamma_1$ (and their inverses) which uses a minimum
number of elements distinct from $g_{11}$ or $g_{11}^{-1}$ in the
product.  Let $g = \prod_{i = 1}^sk_i$ be such a product
decomposition, and let $t$ be the number of elements distinct from
$g_{11}$ and $g_{11}^{-1}$ in the product (note that $t>0$).  Let $P =
|\left\{i\in\left\{1,2,\ldots,s\right\}\,|\,k_i = g_{11}\right\}|$ and
let $N = |\left\{i\in\left\{1,2,\ldots,s\right\},|\,k_i =
g_{11}^{-1}\right\}|$.  Since $g$ must have orbital $A$, and
$\bar{A}\subset A_1$, we see that $P = N$.  In particular, by
inserting an equal number of elements $g_{11}$ and $g_{11}^{-1}$ in
appropriate places, we can think of $g$ as actually being a product of
conjugates of elements of $\Gamma_1\backslash\left\{g_{11}\right\}$
(and their inverses) by elements of the form $g_{11}^j$, where
$j\in\Z$.  Let $g = \prod_{i=1}^tv_i$ be such an expression for $g$,
formed as described.  

If $t = 1$, then we note that $g = v_1$, and that $v_1$ is a conjugate
of an element, $h$, of $\Gamma_1\backslash\left\{g_{11}\right\}$ (or
an inverse of one such), which conjugate has orbital $A\subset A_2$.
But now $A_2$ is an orbital of depth two, and $A$ is of depth at least
two, so that $v_1\in\Lambda_2$ or $v_1^{-1}\in\Lambda_2$.  In this
case, we see that $g\in\ H_2$.

Now let us assume that $t>1$.  Let $w = \prod_{i = 2}^tv_i$.  Now $g =
v_1w$.  Suppose $v_1$ has no orbital intersecting $A$, in this case
$w|_A = g|_A$, so that $w$ would be an element of $A_f$ which has a
product decomposition using elements of $\Gamma_1$ with fewer uses of
elements distinct from $g_{11}$ and $g_{11}^{-1}$, which is
impossible, by the construction of $g$.  So $v_1$ has an orbital $B_1$
so that $B_1\cap A\neq \emptyset$.  Now, by the definition of the
elements $v_i$, we know that the depth of $B_1$ is at least two in
$H_1$, so that $B_1\subset A_2$.  Define $C_1 = A$, and for each
integer $i$ where $2\leq i \leq t$, inductively define
$C_i=C_{i-1}v_{i-1}$.  Suppose that for some
$i\in\left\{1,2,\ldots,t\right\}$, $C_i\cap Supp(v_i) = \emptyset$.
In this case, the product $v_1v_2\cdots v_{i-1}v_{i+1}\cdots v_t$,
which is the product for $g$ without the $v_i$ term, will be a shorter
product which equals $f$ on $A$, which violates our construction.  In
particular, each $C_i$ has non-trivial intersection with an orbital
$B_i$ of $v_i$.  Now $B_1$ is an orbital of $v_1$ that has non-trivial
intersection with $A=C_1$, and so $B_1\subset A_2$, but now this
implies that $C_2\subset A_2$, and so $B_2\subset A_2$, This same
argument repeated $t-2$ more times shows that for each
$i\in\left\{1,2,\ldots,t\right\}$, we have $C_i\subset A_2$, and so
$B_i \subset A_2$.  In particular, each $v_i\in\Lambda_2$, or
$v_i^{-1}\in\Lambda_2$.  But now we have shown that $g$ is a product
of elements of $H_2$, so our claim is proven.
\vspace{.1 in}

Now define $g_{21}$ to be a controller for $H_2$ on $A_2$, and further
define $X_2$ to be the set of towers for $H_2$ with largest element
$(A_2,g_{21})$.  From the claim just proven (that $H_2$
admits towers of arbitrary height with largest elements of the form
$(A_2,h)$ for various $h\in H_2$), we see that $X_2$ has towers of
arbitrary height.

Now define the following set:
\[
\Gamma_2 = \left\{g\in H_2\,|\,(g_{21}^kg) = \lambda,
\lambda\in\Lambda_2, A_2\nsubseteq Supp(g),
k\in\Z\right\}\cup\left\{g_{21}\right\}
\]
This is essentially looking at the interior parts of the generators of
$H_2$ under the orbital $A_2$ of $H_2$ using the controller $g_{21}$.

It is immediate from the definition of controller that $H_2$ is
generated by the finite set $\Gamma_2$.  Suppose that the order of
$\Gamma_2$ is $m_2$, and extend the indexing of $g_{21}$ to the other
elements of $\Gamma_2$, so that $\Gamma_2 = \left\{g_{2i}\,|\, 1\leq i
\leq m_2, i\in\N\right\}$, where $g_{21}$ is the only element of
$\Gamma_2$ which realizes the orbital $A_2$.

We now have that $\Gamma_2$ is a finite set of generators for $H_2$
which satisfies conditions (2) and (3).

Now by construction, $H_2$, $m_2$, $\Gamma_2$, $A_2$, and $X_2$ satisfy the six
constraints.

But now, by temporarily relabeling all of the items mentioned in the
six conditions for the series of groups $(H_i)_{i=1}^{\infty}$, with
index one (the ones that all started with index two), we can repeat
the argument that generated $H_2$ et al from $H_1$, to create a new
subgroup $H_2$, and so forth.  Remembering our original labelling, we
see our new group and objects are $H_3$, et al, which again satisfy
conditions (1--4), (6).  But now that $A_3$ satisfies condition (6),
and $A_2$ satisfies condition (6) vis-a-vis $A_1$, we see that $X_3$
satisifies condition (5).  We can induct on this argument to now build
the sequence of groups $(H_i)_{i=1}^{\infty}$ as promised, that
satisfy all of the six conditions.

  Let $T =\left\{(A_i,g_{i1})|i\in\N\right\}$.  By construction,
$g_{i1}\in H_i\leq H$ has orbital $A_i$ and $\overline{A}_{i+1}\subset
A_i$ for all $i\in\N$.  In particular, $T$ is an infinite tower for
$H$.  \qquad$\diamond$

\subsection{$W$ and other groups}
In this subsection, we will discuss how $W$ relates to the various
groups that we have found in the non-solvable subgroups of $\ploi$.
This will result in proofs of the Lemmas \ref{WInStuff} and \ref{stuffInW}.

We first pick representations of the groups $(\wr\Z\wr)^{\infty}$,
$(\wr\Z)^{\infty}$, and $(\Z\wr)^{\infty}$ in $\ploi$.  Our
presentations of these groups in Thompson's group $F$ will be more
explicit than that given in Brin's \cite{BrinEG}, and his proof that
the resulting groups really are isomorphic with the wreath products
named carries through with no difficulties, although in our concrete
situation we can also use the outline of the proof in the introduction
for at least the group $(\wr\Z)^{\infty}$.

Namely, define $\alpha\in\ploi$ to be the element so that given any
$x\in I$, we have
\[
x\alpha = \left\{\begin{array}{ll}
\frac{1}{4}x & 0\leq x < \frac{1}{4},
\\
x-\frac{3}{16} & \frac{1}{4} \leq x<\frac{7}{16},
\\
4x - \frac{3}{2} & \frac{7}{16}\leq x<\frac{9}{16},
\\
x+\frac{3}{16} & \frac{9}{16}\leq x<\frac{3}{4},
\\
\frac{1}{4}x + \frac{3}{4}& \frac{3}{4}\leq x\leq 1,
\end{array}\right.
\]
and define $\beta_0\in\ploi$ to be the element so that given any $x\in
I$, we have
\[
x\beta_0 =\left\{\begin{array}{ll}
x&0\leq x<\frac{7}{16},
\\
2x-\frac{7}{16}&\frac{7}{16} \leq x< \frac{15}{32},
\\
x+\frac{1}{32} & \frac{15}{32}\leq x<\frac{1}{2},
\\
\frac{1}{2}x + \frac{9}{32} & \frac{1}{2} \leq x < \frac{9}{16},
\\
x & \frac{9}{16}\leq x\leq 1.
\end{array}\right.
\]
The graphs of these elements (superimposed) are given below.  

\begin{center}
\psfrag{a}[c]{$\alpha$}
\psfrag{b1}[c]{$\beta_0$}
\includegraphics[height=340pt,width=340 pt]{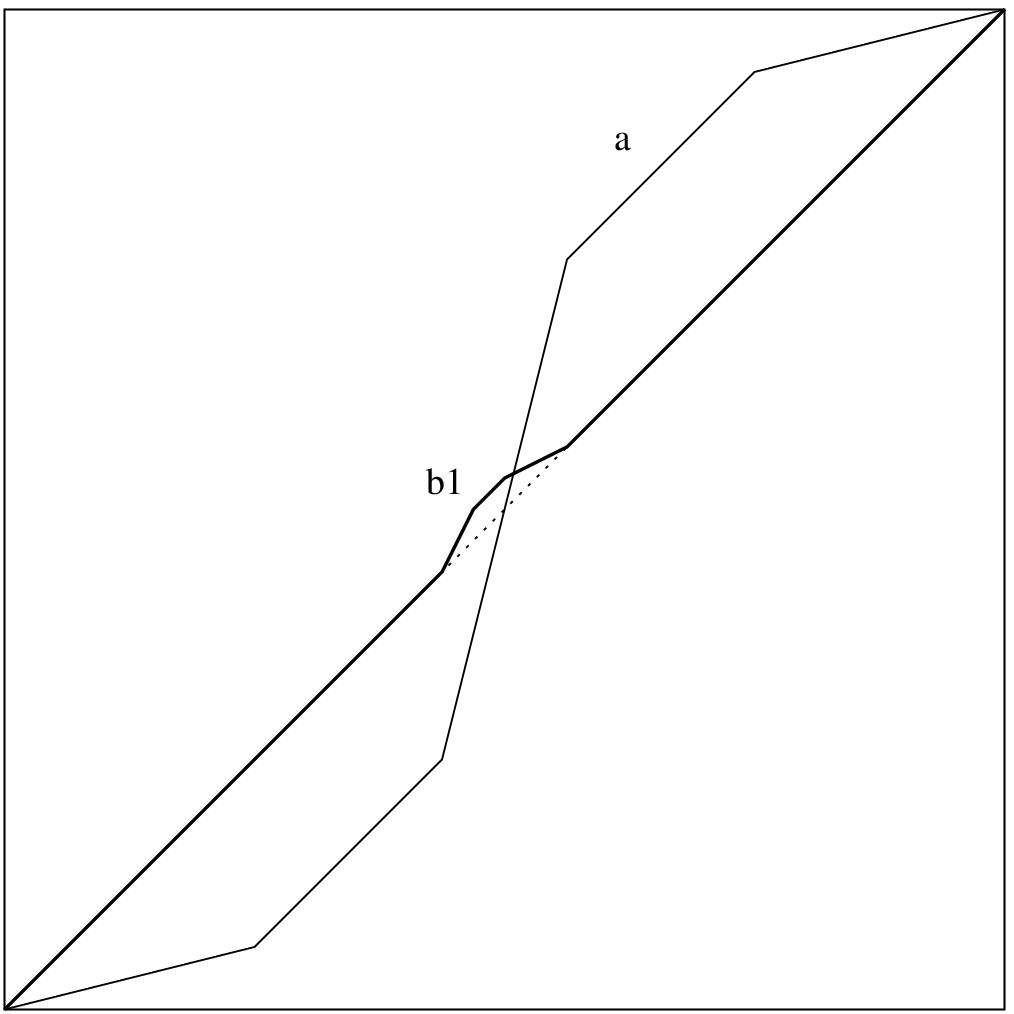}
\end{center}

Now define $\beta_k =
\beta_0^{\alpha^k}$ for each integer $k$.  In particular, $\beta_{1}$
is given by the rule
\[
x\beta_1=x\alpha^{-1}\beta_0\alpha = \left\{\begin{array}{ll}
x &0\leq x<\frac{1}{4},
\\
2x-\frac{1}{4} &\frac{1}{4}\leq x < \frac{3}{8},
\\
x+\frac{1}{8} &\frac{3}{8}\leq x < \frac{1}{2},
\\
\frac{1}{2}x+\frac{3}{8} & \frac{1}{2}\leq x<\frac{3}{4},
\\
x&\frac{3}{4}\leq x\leq 1.
\end{array}\right.
\]
Observe that the support of $\beta_0$ is
$(\frac{7}{16},\frac{9}{16})$, and since $\frac{7}{16}\beta_1 =
\frac{9}{16}$, the support of $\beta_0$ is contained in a single
fundamental domain of $\beta_1$.  In particular, given any $i\in\Z$,
we see that the support of $\beta_{i-1}$ is contained in a
single fundamental domain of the support of $\beta_i$, since these two
elements are conjugates of $\beta_0$ and $\beta_1$.  Now following
Brin in \cite{BrinEG}, the group $(\wr\Z\wr)^{\infty}$ is the group
generated by the full collection of the $\beta_i$, while
$(\wr\Z)^{\infty}$ is the group generated by the $\beta_i$ with $i$ a
negative index, while $(\Z\wr)^{\infty}$ is the group generated by the
$\beta_i$ where we allow only the non-negative indices.  The point is
that whenever $i<j$ are indices, the support of $\beta_i$ is fully
contained inside a single fundamental domain of $\beta_j$.  The reader
may have realized that the group generated by the signatures of the
tower found in Lemma \ref{messyWreath} is isomorphic to
$(\wr\Z\wr)^{\infty}$.  Note that all of the $\beta_i$ are found in
Thompson's group $F$.

We will call $\beta_{-1}$ the top generator of $(\wr\Z)^{\infty}$, and
$\beta_0$ the bottom generator of $(\Z\wr)^{\infty}$.

\bl 
$W$ embeds in both $(\wr\Z)^{\infty}$ and
$(\Z\wr)^{\infty}$.  
\el

pf:
We first embed $W$ in $(\wr\Z)^{\infty}$.  For each $i\in\N$, define  
\[
\Gamma_i = \left\{\gamma^i_j\,|\, \gamma^i_j = \beta_j^{\beta_{i+1}},
1\leq j\leq i, j\in\N\right\}
\]
Note that each collection $\Gamma_i$ generates a group isomorphic to
$W_i$, by the argument given in the introduction after the discussion
of $W_2$, or also by the details of Brin in \cite{BrinEG}.  Further,
the supports of the generators in $\Gamma_i$ are all disjoint from the
supports of the generators in $\Gamma_j$ whenever $i$, $j\in\N$ with
$i\neq j$, so that the elements of $\Gamma_i$ found in
$(\Z\wr)^{\infty}$ above commute with the elements of $\Gamma_j$ in
this case.  Hence, the set
\[
\Gamma = \cup_{i\in\N}\Gamma_i
\]
generates a group
\[
\langle \Gamma \rangle\cong \bigoplus_{i\in\N}W_i \cong W.\]

(Note: We will use this realization of $W$ throughout the rest of this
subsection when we refer to our realization of $W$ in $\ploi$.  When
we refer to ``the first $n$ summands of $W$'' we will mean the $n$
subgroups (each isomorphic to a $W_i$) of $W$ corresponding to the
individual groups generated by the sets $\Gamma_i$ for $1\leq i\leq
n$, $i\in\N$.)

We now embed $W$ in $(\wr\Z)^{\infty}$ in a similar fashion, finding
copies of each $W_i$ in $(\wr\Z)^{\infty}$, all of which occur with
mutually disjoint supports in $I$, the union of their generators will
then generate a group isomorphic to $W$. Let $i\in\N$, and define
\[
\Upsilon_i=\left\{\theta^i_j\,|\, \theta^i_j =
\beta_{-i+j-2}^{\beta_{-1}^i}, 1\leq j\leq i, j\in\N\right\}.
\]
so that $\Upsilon_i$ is the collection of the $i$ generators beneath
$\beta_{-1}$ of the $i$'th conjugate of the generators of
$(\wr\Z)^{\infty}$ by $\beta_{-1}$.  Each collection $\Upsilon_i$
therefore generates a group isomorphic to $W_i$, while if $i\neq j$,
any generator in $\Upsilon_i$ has disjoint support from the generators
of $\Upsilon_j$, so that the union 
\[
\Upsilon =
\cup_{i\in\N}\Upsilon_i
\]
 has the property that
\[
\langle \Upsilon\rangle\cong \bigoplus_{i\in\N}W_i\cong W.
\]
\qquad$\diamond$

\bl 

Neither $(\wr\Z)^{\infty}$, nor $(\Z\wr)^{\infty}$, nor
$(\wr\Z\wr)^{\infty}$ embed in $W$.  

\el 

pf:

If we show that one of $(\wr\Z)^{\infty}$ or $(\Z\wr)^{\infty}$ fails
to embed in $W$, then we will have shown that $(\wr\Z\wr)^{\infty}$ fails
to embed in $W$, so we will say nothing more about
$(\wr\Z\wr)^{\infty}$.

From the details of the proof of Lemma \ref{shortTower}, we see that
the orbital of a signed orbital in a tall tower of a subgroup $G$ of
$\ploi$ survives to be an orbital of an element in the derived group
$G'$.  If we let $G^{(n)}$ temporarily represent the $n$'th derived
group of $(\Z\wr)^{\infty}$, as realized above, then each of the
orbitals of the generators of $(\Z\wr)^{\infty}$ will be orbitals of
some element of $G^{(n)}$.  None of these elements of $G^{(n)}$ with
orbitals shared with a generator $\beta_{k}$ (with $k>0$, $k\in\N$)
will commute with the generator $\beta_0$ in $(\Z\wr)^{\infty}$.  Now
suppose we have found an embedded copy of $(\Z\wr)^{\infty}$ in our
realization of $W$ above.  There is $N\in\N$ so that the orbitals of
the image of the bottom generator of $(\Z\wr)^{\infty}$ in $W$ is
contained in $\bigoplus_{1\leq i\leq N} W_i$, the first $N$ summands used
in the definition of $W$, since elements of $\ploi$ have finitely many
orbitals.  Now the $N$'th derived group of $W$ is trivial over the
domain of these summands, by the argument of Lemma
\ref{solvableClassification}, so that the image of the bottom
generator of $(\Z\wr)^{\infty}$ in $W$ must commute with all of the
elements of the $M$'th derived group of $(\Z\wr)^{\infty}$ in $W$.
Therefore, $(\Z\wr)^{\infty}$ cannot embed in $W$.

To see that $(\wr\Z)^{\infty}$ does not embed in $W$.  Let $G^{(n)}$
now represent the $n$'th derived group of $(\wr\Z)^{\infty}$.  Suppose
we have an embedding of $(\wr\Z)^{\infty}$ in $W$, and let
$\alpha_{-1}$ represent the image of the top generator $\beta_{-1}$ of
$(\wr\Z)^{\infty}$ in $W$.  Note that for any positive integer $m$,
$\beta_{-1}$ does not commute with any of the non-trivial elements of
the $G^{(m)}$ in $(\wr\Z)^{\infty}$ and that this collection is
non-empty.  This implies that $\alpha_{-1}$ must not commute with some
elements of the image of $G^{(m)}$ in $W$.  Observe that this
image is contained in $W^{(m)}$.  Now suppose $n$ is the largest index
where the summand $W_n$ of $W$ shares some support with
$\alpha_{-1}$. $\alpha_{-1}$ must commute with every element of
$W^{(n)}$ in $W$, since the group $W^{(n)}$ has no support in the
supports of the first $n$ summands of $W$, by Lemma
\ref{solvableClassification}, but this contradicts our statement that
the subgroup $G^{(n)}$ of $W^{(n)}\leq W$ has elements that fail to
commute with $\beta_{-1}$ for any $n>0$.  \qquad$\diamond$

\subsection{$W$ in arbitrary non-solvable subgroups of $\ploi$}
The lemma below is a restatement of Lemma \ref{arbTowersW} and
completes the proof of Theorem \ref{nonSolveClassification}.

\bl

If $H$ is a non-solvable subgroup of $\ploi$ then $H$ contains a subgroup isomorphic to $W$.

\el
pf:

Suppose $H$ is a non-solvable subgroup of $\ploi$.  By Lemma
\ref{WInStuff} we know that $W$ embeds in $(\wr\Z)^{\infty}$
and $(\Z\wr)^{\infty}$, and therefore also into $(\wr\Z\wr)^{\infty}$.
Therefore, if $H$ admits infinite towers then we already have the
result, so let us assume that $H$ does not admit infinite towers.  In
particular, $H$ admits towers of arbitrary finite height, $H$ is balanced,
and $H$ admits no transition chains of length two.

Since $H$ does not admit infinite towers, the depth of any signed
orbital of $H$ is well defined and finite.  Since $H$ is not the
trivial group, $H$ has a non-empty collection of orbitals.  The
analysis now breaks into two cases.

\underline{{\bf Case 1:}}

Suppose $H$ admits no orbital that supports towers of arbitrary
height.  In this case the depth of any orbital of $H$ is well defined,
every orbital of $H$ has finite depth, and given any $n\in\N$, $H$ has
orbitals with depth greater than $n$.

Now, pick an element $g^1_1$ of $H$ so that $\hat{T}_1 =
\left\{(B_1^1,g_1^1)\right\}$ is a tower of height one for $H$.
$g^1_1$ will be our generator for $W_1$.  $g_1^1$ has finitely many
orbitals, and so there is a maximum depth $j_1$ of the orbitals of $H$
that are not disjoint from the support of $g_1^1$.  We will now pick
our remaining generators from the group $H^{(j_1)}$, the $j_1$'st
derived subgroup of $H$.  We note that no element in $H^{(j_1)}$ can
have support intersecting $g_1^1$, since $H^{(j_1)}$ has trivial
support over the orbitals of $H$ of depth less than or equal to $j_1$
as a consequence of the details of the proof of Lemma
\ref{solvableClassification}.  We also observe that $H^{(j_1)}$ still
admits towers of arbitrary height, and infinitely many orbitals, of
arbitrary finite depth.  We now find a tower $\hat{T}_2 =
\left\{(B_1^2,g_1^2),(B_2^2,g_2^2)\right\}$ for $H^{(j_1)}\leq H$ of
height two.  Now the signatures of $\hat{T}_2$ admit a finite total
number of orbitals, and therefore the union of this collection of
element orbitals is contained in the union of the collection of
orbitals of $H$ of depth less than some integer $j_2>j_1$.  We
therefore will pick a tower $\hat{T}_3=\left\{(B_1^3,g_1^3),
(B_2^3,g_2^3), (B_3^3,g_3^3)\right\}$ for $H^{(j_2)}$ which has
signatures whose supports must be disjoint from the supports of the
signatures of the first two towers $\hat{T}_1$ and $\hat{T}_2$.  We
can continue in this fashion to inductively define towers $\hat{T}_k$
and integers $j_{k-1}$ for each positive integer $k$ so that the
integers $j_{k}$ are always getting larger, and so that the towers
$\hat{T}_k$ always have height $k$ and have signatures which are
disjoint in support from the signatures of the previous towers.  (This
last follows since $\hat{T}_k$ has all of the orbitals of its
signatures in orbitals of $H^{(j_{k-1})}$ whose supports are away from
the orbitals of depth less than $j_{k-1}$ of $H$, which orbitals
contain all of the supports of the signatures of all of the towers
with smaller index).

Let $k\in\N$.  Let $\hat{G}_k$ represent the group generated by the
signatures of $\hat{T}_k$.  We can use the techniques of the proof of
Lemma \ref{deepWreathZ} to replace $\hat{T}_k$ with a new tower $T_k$
supported by a subset of orbitals of $H$ that support $\hat{T}_k$, so
that the signatures of $T_k$ generate a group $G_k$ isomorphic to
$W_k$.  Do this for all $k\in\N$.

Now the union of all the signatures of all of the towers $T_k$ forms a
collection of generators of a group isomorphic to $W$.

\underline{{\bf Case 2:}}

Suppose now that $H$ admits an orbital $A$ that supports towers of
arbitrary height.

If $A$ is not an orbital of any element of $H$ then $A$ can be written
as a union of an infinite collection of nested element orbitals of
$H$, so that $H$ would then admit an infinite tower, therefore there
is an element $d$ of $H$ so that $(A,d)$ is a signed orbital of depth
one for $H$.

We will now restrict our attention to a special subgroup $H_d$ of $H$
which is directed by the element $d$, in a sense that will be made
clear.  Given any element $h\in H$, let $k_h$ and $j_h$ represent the
smallest positive integers so that $h^{k_h}$ and $d^{j_h}$ satisfy the
mutual efficiency condition.  Let
\[
\Gamma_d = \left\{[h^{k_h},_2 d^{j_h}]\,|\, h\in H\right\}\cup\left\{d\right\}.
\]
The elements of $\Gamma_d$ have all of their orbitals properly
contained inside the orbitals of $d$, and since the orbital $A$ of $H$
admits towers of arbitrary height, and any element orbital $B$ which
is properly contained inside $A$ will be realized as an orbital of
some element $g$ of $\Gamma_d$ (note that it does not matter that we
passed to high powers to guarantee the mutual efficiency condition),
we see that the group $H_d = \langle \Gamma_d\rangle$ admits towers of
arbitrary height. We now observe that given any finite set $X$ of
elements of $H_d$ that do not support any signed orbitals of depth one
for $H_d$, and a finite tower $T$ for $H_d$ which also contains no
signed orbital of depth one, we can find a minimal power $k_X$ of $d$
so that so that the tower $T^{d^{k_X}}$ for $H_d$ induced from $T$ via
conjugation of the signatures of $T$ by $d^{k_X}$ will have all of its
signatures with disjoint support from the signatures of $X$, since we can conjugate the tower to be arbitrarily near to an end of an orbital of $d$.
Therefore, for each positive integer $n$, let $\tilde{T}_n$ be a tower
for $H_d$ of height $n$.  Now inductively define towers $\hat{T}_n$
which are towers induced from the $\tilde{T}_n$ by conjugation by
powers of $d$ so that given any positive integer $k$, the tower
$\hat{T}_k$ has signatures whose supports are all disjoint from the
signatures of the towers $\hat{T}_j$ whenever $j<k$ is a positve
integer.

Now follow the procedure at the end of the previous case to improve
the towers $\hat{T}_n$ to new towers $T_n$ so that for each
positive integer $n$, the signatures of the tower $T_n$ generate a
group isomorphic with $W_n$, while preserving the conditions that the
signatures of distinct towers $T_k$ and $T_j$ have disjoint supports
from each other, so that the union of all of the signatures of all of
the towers $T_k$ is a set of generators of a group isomorphic with $W$
in $H$.

\qquad$\diamond$